\documentclass[12pt]{article}        %4/8/02

\usepackage{color,graphicx}
\usepackage{latexsym,amssymb,amsmath}
\numberwithin{equation}{section}

\newtheorem{lemma}{Lemma}
\newcommand{\bsquare}{\hbox{\rule{6pt}{6pt}}}

\newtheorem{theorem}{Theorem}%[section]
\newtheorem{remark}{{\it{Remark}}}
\begin{document}

\title{Non-local Markovian symmetric forms on infinite dimensional spaces\\
{\small{I. The closability  {%\textcolor{red}
{and quasi-regularity}} \\
}}}
%\TitleHead{}          %optional
%\dedicatory{and Here is a Dedication}           %optional
\author{
{{Sergio
 \textsc{Albeverio}}}
          \footnote{Inst. Angewandte Mathematik, 
and HCM, Univ. Bonn, Germany, 
email \, :albeverio@iam.uni-bonn.de},
\quad 
{{Toshinao \textsc{Kagawa}}} \footnote{Dept. Information Systems Kanagawa Univ., Yokohama, Japan }, 
\quad
 {{Yumi \textsc{Yahagi}}} \footnote{Dept. Mathematical information Tokyo Univ. of Information, Chiba, Japan}, 
\\
 and \quad
{{Minoru \textsc{W. Yoshida}}} \footnote{Dept. Information Systems Kanagawa Univ., Yokohama, Japan,  
 email:\, washizuminoru@hotmail.com }
 }

\maketitle

\begin{abstract}
General theorems on the closability and quasi-regularity of non-local Markovian symmetric forms on probability spaces $(S, {\cal B}(S), \mu)$, with $S$ 
Fr{\'e}chet spaces such that 
$S \subset {\mathbb R}^{\mathbb N}$, ${\cal B}(S)$ is the Borel $\sigma$-field of $S$, and $\mu$ is a Borel probability measure on $S$, are introduced. 
Firstly, a family of non-local Markovian symmetric forms ${\cal E}_{(\alpha)}$,
$0 < \alpha < 2$, acting in each given $L^2(S; \mu)$ is 
defined, 
the index $\alpha$ characterizing the order of the non-locality.
Then, it is shown that all the forms ${\cal E}_{(\alpha)}$ defined on $\bigcup_{n \in {\mathbb N}} C^{\infty}_0({\mathbb R}^n)$ are closable in $L^2(S;\mu)$. 
Moreover, sufficient conditions under which the 
closure of the closable 
forms, that are Dirichlet forms,  become strictly quasi-regular, are given. 
Finally, an existence theorem for Hunt processes properly associated to the Dirichlet forms 
is given.  
The application of the above 
theorems to the problem of stochastic quantizations of Euclidean $\Phi^4_d$ fields, for $d =2, 3$, 
by means of these Hunt processes  is indicated.

\medskip
\noindent
%Insert your abstract here. Include keywords, PACS and mathematical
%subject classification numbers as needed.
{\bf{Key words.}} \, \, Non-local Dirichlet forms, \, \,  Fr{\'e}chet spaces, \, \,  Stochastic quantization, \, \,  Euclidean $\Phi^4_3$ field
\end{abstract}

\section{Introduction}
\label{intro}
We consider a space $S$ that is either a real 
 Banach space  $l^p$, $1 \leq p \leq \infty$, 
with suitable weights, 
or the 
 direct product space ${\mathbb R}^{\mathbb N}$ (with ${\mathbb R}$ and  ${\mathbb N}$ the spaces of real numbers and  natural numbers, respectively).
Both will be looked upon  as Fr{\'e}chet spaces. Let $\mu$ be a Borel probability measure on $S$. On the {\it{real}} $L^2(S; \mu)$ space, for each $0 < \alpha < 2$, we give an explicit formulation of $\alpha$-stable type non-local quasi-regular (cf. section IV-3 of [M,R 92]) Dirichlet forms $({\cal E}_{(\alpha)}, {\cal D}({\cal E}_{(\alpha)}))$ 
(with a domain ${\cal D}({\cal E}_{(\alpha)})$), and show the existence of $S$-valued Hunt processes properly associated to $({\cal E}_{(\alpha)}, {\cal D}({\cal E}_{(\alpha)}))$.
$\alpha$-stable is understood in analogy with the $\alpha$-stable Dirichlet forms defined on 
$L^2({\mathbb R}^d)$, for $d \in {\mathbb N}$, e.g., in 
 [Fukushima,Uemura 2012], section 5.

As an application of the above general results, in Example 1 and 2 in section 5, we consider the problem of stochastic quantization of the Euclidean free  field over ${\mathbb R}^d$, the $\Phi^4_2$ and $\Phi^4_3$ fields over ${\mathbb R}^2$ 
and  ${\mathbb R}^3$, i.e., fields with 
no (self) interaction respectively, 
(self) interaction of the $4$-th power.  By using the property that, for example,  the support of the 
Euclidean $\Phi^4_3$ field measure $\mu$  is in some {\it{real}} Hilbert space ${\cal H}_{-3}$ (cf. (5.9) for the explicit definition), which is a subspace of the  Schwartz space of {\it{real}} tempered distributions ${\cal S}'({\mathbb R}^3 \to {\mathbb R})$,  we define 
 an isometric isomorphism $\tau_{-3}$  from ${\cal H}_{-3}$ to "some weighted $l^2$ space"\, (cf. (5.22) for the explicit definition). By making use of $\tau_{-3}$, we then apply the above general theorems formulated on the abstract $L^2(S; \mu)$ space to 
the case of the space $L^2({\cal H}_{-3}; \mu)$ for the 
Euclidean $\Phi^4_3$ field, 
and 
for each $0 < \alpha \leq 1$ we show the existence of an ${\cal H}_{-3}$-valued Hunt process $(Y_t)_{t \geq 0}$ the invariant measure of which is $\mu$.

$(Y_t)_{t \geq 0}$ can be understood as a stochastic quantization of the Euclidean ${\Phi}^4_3$ field 
 realized by a Hunt process through the {\it{non-local}} Dirichlet form $({\cal E}_{(\alpha)}, {\cal D}({\cal E}_{(\alpha)}))$ for $0 < \alpha \leq 1$,
in the sense that ${\Phi}^4_3$-measure $\mu$ is the invariant measure for $(Y_t)_{t \geq 0}$.

1) \quad As far as we know, there has been no explicit
proposal of a general 
  formulation of {\it{non-local}} quasi-regular Dirichlet forms on infinite dimensional (Fr{\'e}chet)  topological vector spaces,
which admit interpretations as Dirichlet forms 
on
 concrete random fields on Fr{\'e}chet spaces.
This is different from the situation associated with the local case, 
 i.e., the case where the associated Markov processes are (continuous) diffusions, where much has been developed and known (cf. a short review given below).

2) \quad Although there have been derived several results on the existence of (continuous) diffusions (i.e., roughly speaking,  that are associated to quadratic forms and generators of local type) corresponding with 
stochastic quantizations of, e.g., ${\Phi}^4_2$ and  
 ${\Phi}^4_3$ Euclidean fields (see below for references), as far as we know, there 
do not 
exist no explicit corresponding considerations for {\it{non-local}} type Markov processes, 
associated with {\it{non-local}}  Dirichlet forms related to such fields.
The only examples so far provided have been obtained by  subordination, 
for the Euclidean free field (on ${\cal S}'({\mathbb R}^d), \, d \geq1$) and for the ${\Phi}^4_2$-Euclidean field with space cut-off (interaction in a bounded euclidean space region) a non-local type stochastic quantization procedure has been discussed in 
 [A,R{\"u}diger 2003] 
 (see (5.29), (5.62) and Remark 12).

The present paper is a first development that gives answers to the 
above mentioned 
open problems 1), i.e.,  to give an abstract formulation in this concern,  and 2), i.e.,  to perform a corresponding consideration in the case of ${\Phi}^4_3$-model.

Before giving an explanation of the contents of the present paper, we give a brief review of the 
theory of  Dirichlet forms and its 
applications to 
stochastic quantization having correspondences with the present considerations. 

In the case where the state space $E$ is finite dimensional, much 
 has  been done on the considerations of both {\it{local}} and {\it{non-local}} Dirichlet 
forms and semi-Dirichlet forms (non-symmetric Dirichlet forms).
The natural setting is the one where the Hilbert space, where the Dirichlet forms are defined 
(as quadratic forms), is  
 $L^2(E;m)$, 
 with $E$ a general locally compact separable metric space (when $E$ is a topological vector space, the dimension of the space is thus finite) and $m$ a positive Radon measure on it (cf., e.g., [Fukushima 80], [F,Oshima,Takeda 2011], [F,Uemura 2012], [Hoh,Jacob 96], [Masamune,Uemura,Wang 2012], [Schlling,Wang 2014], [Shiozawa,Uemura 2014], [A],
[A,Song 93], [A,Ugolini 2015] 
  and references therein).  Also, many results have been developed on the theory of 
general 
(non-symmetric) local Dirichlet forms  defined on $L^2(E;m)$, with general topological spaces including the case of some infinite dimensional topological vector spaces, and $m$ some  Radon measures on them (cf., e.g., [A 2003], [A,H-K 76], [A,H-K1 77], [A,H-K2 77], [A,H-K,S 77],  
[A,Ma,R 2015], [A,Ma,R1 93], [A,Ma,R2 93], [A,R 89, 90,91], [Kusuoka-sige 82], [M,R 92], [A,R{\"u}diger 2003], [Schumuland 90], [Shigekawa 2000] and references therein). 
In the general abstract framework, in particular [M,R 92], [A 2003], the Dirichlet forms need not be local ones, however all examples 
except those considered through the framework of subordination given in [A,R{\"u}diger] (cf. Remark 12 in section 5) 
 treated in above references, treating the case where $E$ is {\it{infinite}} dimensional, 
concern  Dirichlet forms theories that are local ones (i.e., associated with differential operators).

One of the most important applications of the theory of Dirichlet forms on $L^2(E;m)$, $E$ an infinite dimensional space, is to construct an $E$-valued Markov process, the invariant measure of which is $m$ (assuming $m$ to be a probability measure). In particular, if $L^2(E;m)$ is a 
space associated to the 
 Euclidean quantum field (for the construction of $m$ cf., e.g., [Glim,Jaffe 87], [Simon 74]), the constructed Markov process is referred to as stochastic quantization of the field.  In the  literature 
quoted above concerning local Dirichlet forms 
 we find several applications to stochastic quantizations, where the Markov processes constructed are (continuous) diffusion processes associated to {\it{local}} Dirichlet forms.  
 For the stochastic quantizations of Euclidean $P(\phi)_2$ fields on ${\cal S}'({\mathbb R}^2 \to {\mathbb R})$, i.e., fields with polynomial interactions, fields with trigonometric interactions and exponential interactions 
on ${\cal S}'({\mathbb R}^2 \to {\mathbb R})$, these considerations by means of the {\it{local}} Dirichlet form arguments were completed by [A,H-K 76], [A,H-K1 77], [A,H-K2 77], [A,H-K,S 77], 
[Bo,Ch,Mi 88],
 [A,Hida,Po,R,Str1 89], [A,Hida,Po,R,Str2 89],
[A,Ma,R 2015], [A,Ma,R1 93], [A,Ma,R2 93], [A,R 89, 90, 91], 
[Hida,Kuo,Po,Str 93].
In this direction of the application of {\it{local}} Dirichlet forms, there also are corresponding considerations 
for measures $m$ describing  
 infinite particle systems (cf., e.g., [Y 96], [Osada 96], [A,Kondratiev,R{\"o}ckner 98], [A,Ugolini 2015] and references therein).
For works on stochastic quantization of ${\Phi}^4_2$ using other methods see 
 [Jo-La,Mi 85],
[Da Prato,Debussche 03], 
[Da Prato,Tubaro 2000] 
 and references in [Albeverio, Ma, R{\"o}ckner], [A, Kusuoka-sei 2017].
 For other models in two or less Euclidean space dimensions where stochastic quantization methods have been applied see the recent work [A,DeVe,Gu1 2020], [A,DeVe,Gu2 2020] 
and references therein.

The problem of stochastic quantization of the $\Phi^4_3$ Euclidean field on the $3$-dimensional torus ${\mathbb T}^3$ was  solved firstly by [Hairer 2014], by not passing through the arguments by means of the Dirichlet form theory, but through the stochastic  partial differential equation (SPDE in short) arguments with the theory of regularity structures developed there.  A little more precisely, 
in [Hairer 2014] and subsequent works (see below) 
 the existence of a solution, that is a (continuous) diffusion, of an SPDE, which were expected to be satisfied by a solution of the stochastic quantization of an Euclidean  $\Phi^4_3$ model, restricted to ${\mathbb T}^3$, 
is shown 
(cf., e.g., Theorem 1.1 of [A,Liang,Zegarlinski 2006]).  Since the appearance of [Hairer 2014], there have been  developed several results on (continuous) diffusion type stochastic quantization of the $\Phi^4_3$ Euclidean field (cf., e.g., [Cat,Chouk 2018], [Hairer,Mattingly 2016], 
[Mourrat,Weber 2017]
[A,Kusuoka-sei 2020], 
[A,Kusuoka-sei 2021], 
[Gu,Ho 2019], [Gu,Ho 2021] and references therein). [Zu,Zu1 2018], [Zu,Zu2 2018] considered the existence of a local Dirichlet form corresponding to the Euclidean $\Phi^4_3$ field on ${\mathbb T}^3$ by making use of the result on existence of the solution of the SPDE by [Hairer 2014].

Now, let us give a brief summary of the contents of the single sections of this work.
Theorems 1, 2, 3 and 4 are statements corresponding to the {\it{non-local}} Dirichlet forms constructed explicitly on $L^2(S;\mu)$, where $S$ denotes  the weighted $l^p$, $1 \leq p \leq \infty$, spaces and subspace  of the direct product ${\mathbb R}^{\mathbb N}$, respectively,
and $\mu$ is a Borel probability measure on $S$.  Theorem 5 is a restatement of the Bochner-Minlos theorem in this framework. Theorem 6 
contains an application to 
 a {\it{non-local}} type stochastic quantization of the Euclidean $\Phi^4_d$, $d =2, 3$, theory, derived by applying the above general theorems to the $\Phi^4_d$ field by making use of an isometric isomorphism between a Hilbert space, taken as  support of the $\Phi^4_d$ field measure, and a weighted $l^2$ space on which Theorems 1, 2 and 4 are considered.

In section 2, for each $0 < \alpha <2$ we define the symmetric Markovian form ${\cal E}_{(\alpha)}$ on $L^2(S; \mu)$ (cf. (2.8), (2.9) and (2.10)), which is given by using the conditional probability of $\mu$, 
where the index $\alpha$ characterizes the order of the non-locality. 
The definition is a natural analogue of the one for  $\alpha$-stable type 
{\it{non local}} 
Dirichlet forms on ${\mathbb R}^d$, $d< \infty$ (cf. (2.13) in Remark 2 of the present paper and (5.3), (1.4)  of  [Fukushima,Uemura 2012]). 
It has a natural correspondence to the one for local classical Dirichlet forms 
on infinite dimensional topological vector spaces 
(cf. [A,R 89, 90, 91] and [Kusuoka-sige 82]), which is carried out by making use of directional derivatives.
 Theorem 1 shows that the non-local symmetric forms ${\cal E}_{(\alpha)}$ with the domain 
${\cal F}C^{\infty}_0$, the space of regular cylindrical functions (cf. (2.7)), are closable in $L^2(S;\mu)$, hence their closures are {\it{non-local}} Dirichlet forms denoted by $({\cal E}_{(\alpha)}, {\cal D}({\cal E}_{(\alpha)}))$.

We give a proof of Theorem 1 in section 3, by considering separately the cases $0 < \alpha \leq 1$ and $1 < \alpha <2$.

Section 4 discusses  {\it{strictly}} quasi-regularity (cf., e.g., section IV-3 and section V-2 of [M,,R 92]) of the {\it{non-local}} Dirichlet form $({\cal E}_{(\alpha)}, {\cal D}({\cal E}_{(\alpha)}))$ on $L^2(S; \mu)$, by which the existence of Hunt processes properly associated with $({\cal E}_{(\alpha)}, {\cal D}({\cal E}_{(\alpha)}))$ is guaranteed.  Theorem 2 and  Theorem 3 give sufficient conditions under which $({\cal E}_{(\alpha)}, {\cal D}({\cal E}_{(\alpha)}))$ 
for $0 < \alpha \leq 1$ and for $1 < \alpha <2$ become {\it{strictly}} quasi-regular, respectively.   In the case where $S$ is a weighted $l^p$, $1 \leq p \leq \infty$, space, the proofs are carried out by making 
an efficient 
use of the structures of  weighted spaces.  In short,  the concept of the quasi-regularity on ${\cal D}({\cal E}_{(\alpha)})$ is 
a requirement that ${\cal D}({\cal E}_{(\alpha)})$ contains  a sequence of subsets which is dense in 
 ${\cal D}({\cal E}_{(\alpha)})$,
the elements of which are functions on $S$ with compact supports.  By making use of the weight of the $l^p$ spaces, compact sets in $S$ can be given explicitly (cf. (4.15)), and then 
a corresponding 
 sequence of subsets in 
 ${\cal D}({\cal E}_{(\alpha)})$ can be constructed.

Theorem 4 in section 5 gives the statement of the existence of $S$-valued Hunt processes 
that are 
 properly associated to the {\it{non-local}} strictly quasi-regular Dirichlet forms 
 $({\cal E}_{(\alpha)}, {\cal D}({\cal E}_{(\alpha)}))$ on  $L^2(S; \mu)$ given by Theorems 2 and 3.  In the same section, we recall the Bochner-Minlos theorem (cf. Theorem 5), and then by making use of this theorem and Theorem 4 we solve the problems of stochastic quantization of the Euclidean free field on 
${\cal S}'({\mathbb R}^d \to {\mathbb R})$, $d \in {\mathbb N}$, and  the Euclidean $\Phi^4_d$ field on 
${\cal S}'({\mathbb R}^d \to {\mathbb R})$, $d =2, 3$ in Example 1, and in  Example 2, respectively.  As a consequence, the stochastic quantizations are realized by Hunt processes properly associated to the {\it{non-local}} 
Dirichlet forms 
$({\cal E}_{(\alpha)}, {\cal D}({\cal E}_{(\alpha)}))$ for $0 < \alpha \leq 1$.
In these examples, in order to apply Theorem 4 to the Euclidean free field measure and 
to the 
$\Phi^4_2$ and $\Phi^4_3$ field measures, respectively, on the Schwartz space of tempered distributions, we firstly certify that these measures have support in the Hilbert spaces ${\cal H}_{-2}$ and  ${\cal H}_{-3}$ (cf. (5.11)), respectively,  and define 
the isometric isomorphism 
$$\tau_{-2} \, : {\cal H}_{-2} \to {\mbox{a weighted $l^2$ space, $l^2_{(\lambda_i^4)}$ \, (cf. (5.27))}}$$
and, 
$$\tau_{-3} \, : {\cal H}_{-3} \to {\mbox{a weighted $l^2$ space, $l^2_{(\lambda_i^6)}$ \, (cf. (5.56))}},$$
and we then identify ${\cal H}_{-2}$ with $l^2_{(\lambda_i^4)}$ and 
${\cal H}_{-3}$ with $l^2_{(\lambda_i^6)}$, respectively.  For the weighted $l^2$ spaces  $l^2_{(\lambda_i^4)}$ and 
 $l^2_{(\lambda_i^6)}$ we can apply Theorems 1, 2, 4 and have 
$l^2_{(\lambda_i^4)}$ and 
 $l^2_{(\lambda_i^6)}$ valued Hunt processes $(X_t)_{t \geq 0}$ properly associated to corresponding Dirichlet forms $({\cal E}_{(\alpha)}, {\cal D}({\cal E}_{(\alpha)}))$, respectively (in case of $\alpha = 1$, for the Euclidean free field on ${\cal S}'({\mathbb R}^d \to {\mathbb R})$, $d \in {\mathbb N}$ an $l^2_{(\lambda^4_i)}$-valued Hunt process is defined, and generally, for $0 < \alpha \leq 1$,  $l^2_{(\lambda^6_i)}$-valued Hunt processes are defined 
for both the Euclidean free field on ${\cal S}'({\mathbb R}^d \to {\mathbb R})$, $d \in {\mathbb N}$ and 
the Euclidean $\Phi^4_d$ field on 
${\cal S}'({\mathbb R}^d \to {\mathbb R})$, $d =2, 3$).
  Finally, we define the corresponding ${\cal H}_{-2}$ and  ${\cal H}_{-3}$ valued Hunt processes $(Y_t)_{t \geq 0}$ which are inverse images of $(X_t)_{t \geq 0}$ through $\tau_{-2}$ and  $\tau_{-3}$, respectively (see Theorem 6).

In section 6, we 
give a short outlook to 
 future developments, 
in the line of 
 the present formulations and results.

%Your text comes here. Separate text sections with
\section{Markovian symmetric forms individually adapted to each measure space}
%{\textcolor{blue}{
The state space $S$, on which we define the Markovian symmetric forms, is 
one of the following 
Fr{\'e}chet spaces (i.e.,   complete 
infinite dimensional topological vector spaces with a system of countable semi-norms): \\
A weighted $l^p$ space, denoted by  $l^p_{(\beta_i)}$, 
such that, for some $p \in [1, \infty)$ and a weight $(\beta_i)_{i \in {\mathbb N}}$ with  $\beta_i \geq 0, 
i \in {\mathbb N}$,  
\begin{equation}
S = l^p_{(\beta_i)} \equiv \bigl\{  {\mathbf x} = (x_1, x_2, \dots) \in {\mathbb R}^{\mathbb N} \, : \, 
\| {\mathbf x} \|_{l^p_{(\beta_i)}} \equiv (  \sum_{i=1}^{\infty} {\beta}_i |x_i|^p )^{\frac1p} < \infty    \bigr\},
\end{equation}
or 
a weighted $l^{\infty}$ space, denoted by  $l^{\infty}_{(\beta_i)}$, 
such that for  a weight $(\beta_i)_{i \in {\mathbb N}}$ with $\beta_i \geq 0, 
i \in {\mathbb N}$,  
\begin{equation}
S = l^{\infty}_{(\beta_i)} \equiv \bigl\{  {\mathbf x} = (x_1, x_2, \dots) \in {\mathbb R}^{\mathbb N} \, : \, 
\| {\mathbf x} \|_{l^{\infty}_{(\beta_i)}} \equiv \sup_{i \in {\mathbb N}} \beta_i |x_i| < \infty    \bigr\},
\end{equation}
or 
\begin{equation}
S = {\mathbb R}^{\mathbb N}, \, \,  {\mbox{the direct product space with the metric 
$d(\cdot, \cdot)$}} 
\end{equation}
such that for ${\mathbf x}, \,{\mathbf x}' \in {\mathbb R}^{\mathbb N}$,
$$d(\mathbf x, {\mathbf x}') \equiv \sum_{k=1}^{\infty} (\frac12)^k 
\frac{\|\mathbf x- {\mathbf x}'\|_k}{\|\mathbf x- {\mathbf x}'\|_k +1},$$
$$ {\mbox{with}} \qquad 
\| \mathbf x\|_k = ( \sum_{i=1}^k (x_i)^2)^{\frac12}, \, \,  {\mathbf x} =(x_1, x_2, \dots) \in {\mathbb R}^{\mathbb N}.
 $$
{\it{
By making use of the concrete expressions of $S$ as above, i.e., their expressions by means of a subspaces  of ${\mathbb R}^{\mathbb N}$,  abstract discussions on  Dirichlet forms 
(i.e., closed Markovian symmetric forms,  on 
$L^2(S; \mu)$) can be made more concrete.   Moreover, 
 our choice of $S$ permits 
 an effective  application of  the theory of Dirichlet forms to 
certain problems related to 
the stochastic quantization of Euclidean quantum fields  (cf. Remark in the last section).}}
%}}

We denote 
by ${\cal B}(S)$ the 
 Borel $\sigma$-field of $S$. 
Suppose that we are given a Borel probability measure $\mu$ on $(S, {\cal {B}}(S))$. 
For each $i \in {\mathbb N}$, 
%{\textcolor{red}{
let  $\sigma_{i^c}$  be the 
sub $\sigma$-field
%}}
 of 
${\cal B}(S)$ that is generated by the Borel sets 
\begin{equation}
B 
%{\textcolor{blue}{
\equiv
%}}
 \left\{ {\mathbf x} \in S \, \, \, \big| \, x_{j_1} \in B_1, \dots x_{j_n} \in B_n \right\}, \quad 
 j_k \ne i, \, \, B_k \in {\cal B}^1, \, \, k=1, \dots, n, \, \, n \in {\mathbb N}, \, \,
\end{equation}
where ${\cal B}^1$ denotes the Borel $\sigma$-field of ${\mathbb R}^1$, i.e.,  $\sigma_{i^c}$  is the smallest $\sigma$-field that includes every $B$ given by 
(2.4). 
{\it{Namely, ${\sigma}_{i^c}$ is the sub $\sigma$-field of ${\cal B}(S)$ generated by the variables ${\mathbf x} \setminus x_i$, i.e., the variables except of the $i$-th variable $x_i$}}.
For each $i \in {\mathbb N}$, let $\mu(\cdot \, \big| \, \sigma_{i^c})$ be the conditional 
probability,
a one-dimensional probability distribution-valued  $\sigma_{i^c}$ measurable function 
%{\textcolor{blue}{
(i.e., the probability distribution for the $i$-th component $x_i$),
%}}
 that is characterized by (cf. (2.4) of [A,R 91])
\begin{equation}
\mu \big( \{ {\mathbf x} \, \, : \, x_i \in A \} \cap B \big)
= \int_B \mu(A \, \big| \, {\sigma}_{i^c})   \, \mu(d {\mathbf x}),  \quad \forall A \in {\cal B}^1, \, \,
\forall B \in {\sigma}_{i^c}.
\end{equation}
Define
\begin{equation} L^2(S; \mu) \equiv \left\{ f \, \, \,  \Big| \, f : S \to {\mathbb R}, \, {\mbox{measurable and }} \, 
\|f \|_{L^2} = \Bigl( \int_S |f({\mathbf x})|^2 \mu(d {\mathbf x})  \Bigr)^{\frac12} < \infty \right\},
\end{equation}
and
%{\textcolor{blue}{
\begin{equation}
{\cal F}C^{\infty}_0 \equiv 
{\mbox{the $\mu$ equivalence class of}} 
\left\{ f   \, \Big| \, \exists n \in {\mathbb N}, \, f  \in C^{\infty}_0({\mathbb R}^n \to {\mathbb R})  \right\} \subset 
L^2(S; \mu),
\end{equation}
%}}
where $C^{\infty}_0({\mathbb R}^n \to {\mathbb R})$ denotes the space of {\it{real valued}} infinitely differentiable functions on ${\mathbb R}^n$ with compact supports 
(cf. (4.29) below).

%{\textcolor{blue}{
\begin{remark}{}\quad i) \quad For the subsequent discussions, in particular 
those concerning 
the closability of quadratic forms on $L^2(S;\mu)$, we first have to certify that ${\cal F}C^{\infty}_0$ defined by (2.7) is dense in $L^2(S;\mu)$. This can be seen through the same argument performed in section II-3-a) of [M,R 92] for a {\it{Souslin space}},
where the corresponding discussion is carried out for ${\cal F}C_b^{\infty}$ that is defined by
substituting $C^{\infty}_0({\mathbb R}^n \to {\mathbb R})$ in (2.7) by 
$C^{\infty}_b({\mathbb R}^n \to {\mathbb R})$ (see (3.1) in section II-3 of [M,R 92]).
But in the present 
formulation where the state space $S$ is given by (2.1), (2.2) or (2.3), the fact that 
 ${\cal F}C^{\infty}_0$ is dense in $L^2(S;\mu)$ can be seen directly as follows: \, The Borel $\sigma$-field ${\cal B}(S)$ of $S$ is generated by the 
collection of 
finite direct products of the 
sets such that $\{x_i \, :\, a_i \leq x_i \leq b_i\}$,  $(i \in {\mathbb N})$, and for each indicator function $\displaystyle{I_{[a_i,b_i]}(x_i)}$ there exists a sequence of functions $\{f_n(x_i)\}_{n \in {\mathbb N}}$, $f_n(\cdot) \in C^{\infty}_0({\mathbb R}), \, n \in {\mathbb N}$, that has a uniform bound, such that it converges to $\displaystyle{I_{[a_i,b_i]}(x_i)}$ in $L^2(S;\mu)$ by the Lebesgue's bounded convergence theorem.\\
ii) \quad Let $C(S \to {\mathbb R})$ be the space of real valued continuous functions on $S$. For each $S$ that is given by (2.1), (2.2) or (2.3), the fact that ${\cal F}C_0^{\infty} \subset C(S \to {\mathbb R})$ can be also seen 
directly. But, since $S= {\mathbb R}^{\mathbb N}$ defined by (2.3) is equipped the weakest topology, and clearly  ${\cal F}C_0^{\infty} \subset C({\mathbb R}^{\mathbb N} \to {\mathbb R})$,  the same is also true for $S$ as given by (2.1) and (2.2).\\
iii) \quad The completeness of the metric space $S$ is not necessary for the discussions in this section, but it shall be used for the considerations of the {\it{quasi regularity}} of symmetric forms performed in section 4.
%\qquad \qquad \qquad \qquad \qquad \qquad \qquad \qquad \qquad \qquad  \qquad \qquad
%\qquad  \qquad \qquad  \qquad \qquad \quad 
%\bsquare
\end{remark}
%}}

On $L^2(S;\mu)$, 
for any $0 < \alpha <2$, let us define 
the Markovian symmetric form ${\cal E}_{(\alpha)}$  called  {\it{individually adapted Markovian symmetric form 
{%\textcolor{red}
{of index $\alpha$}} to the measure $\mu$}}, 
the definition of which 
is a natural analogue of  the one for  $\alpha$-stable type 
({\it{non local}})  
Dirichlet form on ${\mathbb R}^d$, $d< \infty$ (cf. Remark 2 given below and  (5.3), (1.4)  of  [Fukushima,Uemura 2012]). 
The concept extends to the non local case the one of  local classical Dirichlet forms 
on infinite dimensional topological vector spaces 
(cf. [A,R 89, 90, 91] and [Kusuoka-sige 82]), which is carried out by making use of directional derivatives.
Our  definition is 
as follows:
Firstly, 
for each $0 < \alpha <2$ and  $i \in {\mathbb N}$,
and
for 
the variables \,
$y_i, \, y'_i \in {\mathbb R}^1$, ${\mathbf x}  = (x_1, \dots , x_{i-1},x_i, x_{i+1}, \dots) \in S$ and  \,
${\mathbf x} \setminus x_i \equiv (x_1, \dots, x_{i-1},x_{i+1}, \dots)$, let 
\begin{eqnarray} 
\lefteqn{
\Phi_{\alpha}(u,v;y_i,y'_i,{\mathbf x} \setminus x_i)
} \nonumber \\
 &&\equiv 
\frac{1}{|y_i- y'_i|^{\alpha +1}} 
\times 
\big\{
u(x_1, \dots, x_{i-1}, y_i, x_{i +1}, \dots) - u(x_1, \dots, x_{i-1}, y'_i, x_{i +1}, \dots) \big\}
 \nonumber \\
&&\times 
\big\{
v(x_1, \dots, x_{i-1}, y_i, x_{i +1}, \dots) - v(x_1, \dots, x_{i-1}, y'_i, x_{i +1}, \dots)
\big\}.
\end{eqnarray}

For each $0 < \alpha \leq 1$ and $i \in {\mathbb N}$, define
\begin{equation}
{\cal E}_{(\alpha)}^{(i)}(u,v) \equiv 
\int_S
\Big\{
\int_{\mathbb R} I_{\{y\,:\, y \ne x_i\}}(y_i) \,
\Phi_{\alpha}(u,v;y_i,x_i, {\mathbf x}\setminus x_i)  \, {\mu \big( dy_i \, \big| \, \sigma_{i^c} \big) }
\Big\} \mu(d{\mathbf x}),
\end{equation}
%{\textcolor{blue}{
for any $u, v$ such that the right hand side of (2.9) is finite,
where for a set $A$ and a variable $y$, 
$I_A(y)$ denotes the indicator function, and in the sequel, to simplify the notations, we  denote 
$I_{\{y\,:\, y \ne x_i\}}(y_i)$ by, e.g.,  $I_{\{y_i \ne x_i\}}(y_i)$ or $I_{\{ y_i \ne x_i\}}$.

By ${\cal D}_i$, we denote the subset of the space of real valued ${\cal B}(S)$-measurable functions such that 
the right hand side of (2.9) is finite for any $u, \, v \in {\cal D}_i$.
Let us call $({\cal E}^{(i)}_{(\alpha)}, {\cal D}_i)$ this form, 
${\cal D}_i$ being its domain, 
and then 
 define 
\begin{equation}
{\cal E}_{(\alpha)}(u,v) \equiv \sum_{i \in {\mathbb N}} {\cal E}_{(\alpha)}^{(i)}(u,v), 
\qquad \forall u, \, v \in \bigcap_{i \in {\mathbb N}} {\cal D}_i, 
\end{equation}
For $y_i \ne y'_i$, 
(2.8) is well defined for any real valued ${\cal B}(S)$-measurable functions $u$ and $v$. 
Moreover, 
for {\it{the Lipschiz continuous}} functions 
${\tilde{u}} \in C^{\infty}_0({\mathbb R}^n \to {\mathbb R}) \subset {\cal F}C^{\infty}_0$ and
${\tilde{v}} \in C^{\infty}_0({\mathbb R}^m \to {\mathbb R}) \subset {\cal F}C^{\infty}_0$,  $n,m \in {\mathbb N}$ with compact supports, which are 
%{\textcolor{blue}{
representatives
%}}
of 
$u \in {\cal F}C^{\infty}_0$ and
$v \in {\cal F}C^{\infty}_0$, respectively (see (2.7)), (2.9) and (2.10)  are well defined (the right hand side of (2.10) has only a finite number of sums).  In Theorem 1 given below we will see that (2.9) and (2.10) are well defined for ${\cal F}C^{\infty}_0$, the space of $\mu$-equivalent class, and in particular  
${\cal F}C^{\infty}_0 \subset (\cap_{i \in {\mathbb N}} {\cal D}_i)$.
%}}

For  $1 < \alpha <2$,  we suppose that 
for each $i \in {\mathbb N}$, the conditional distribution 
 $\mu(\cdot \, \big| \, \sigma_{i^c})$ can be expressed by 
%{\textcolor{blue}{
a {\it{locally bounded
%}}
 probability density function}} $\rho(\cdot \, \big| \, \sigma_{i^c})$, 
$\mu$-a.e.,
i.e., $\mu(\cdot | {\sigma}_{i^c}) = \rho (\cdot|{\sigma}_{i^c}) dx$.
Precisely (cf. (2.5) of [A,R91]), there exists a $\sigma_{i^c}$-measurable function 
$0 \leq \rho(\cdot \, \big| \, \sigma_{i^c})$ on ${\mathbb R}^1$ 
and
\begin{equation}
\mu(dy \, \big| \, \sigma_{i^c}) = 
\rho(y \, \big| \, \sigma_{i^c}) \, dy, \quad \mu-a.e.,
\end{equation}
holds, with  $\rho(\cdot \, \big| \, \sigma_{i^c})$ a function such that
for each $ i \in {\mathbb N}$
and
 for any compact 
$K \subset {\mathbb R}$ 
there exists a {\it{ bound}} $L_{K,i} <\infty$, which may depend on $i$ such that 
%{\textcolor{blue}{
\begin{equation}
{\rm{ess}} \, \sup_{y \in 
{K}}  \rho(y \, \big| \, \sigma_{i^c}) \leq L_{K,i}, \qquad  \mu-a.e.,
\end{equation}
where ${\rm{ess}}\sup_{y \in 
{K}}$ is taken with respect to the Lebesgue measure on ${\mathbb R}^1$.
%}}
%}}
%$L_{K_{M,i},i}$
Then define 
the non-local form ${\cal E}_{(\alpha)}(u,v)$, for  $1 < \alpha <2$, by the same formulas as 
(2.9) and (2.10), for all $u, v \in \cap_{i \in {\mathbb N}} {\cal D}_i$,
where
 ${\cal D}_i$ denotes the subset of the space of real valued ${\cal B}(S)$-measurable functions such that 
the right hand side of (2.9), with a given  $\alpha$ such that $1 < \alpha <2$, is finite for any $u, \, v \in {\cal D}_i$.

%{\textcolor{blue}{
\begin{remark} {}\qquad \\
i) \quad
For the ${\cal B}(S)$ measurable function 
$\int_{\mathbb R} I_{\{y_i \ne x_i \}} \Phi_{\alpha}(u,v;y_i,x_i, {\mathbf x} \setminus x_i)
\mu(dy_i \, | \, \sigma_{i^c})$ 
by taking the expectation conditioned by the sub $\sigma$-field $\sigma_{i^c}$, 
through the definition of the conditional probability measures, we can give the  equivalent expressions 
 for  ${\cal E}^{\alpha}_i(u,v)$ defined by (2.9), 
$u, v  \in \cap_{i \in {\mathbb N}} {\cal D}_i$,
as follows:
\begin{eqnarray}
{\cal E}_{(\alpha)}^{(i)}(u,v) &\equiv& 
\int_S
\Big\{
\int_{{\mathbb R} \setminus \{x_i\}} 
%I_{\{y_i \ne x_i\}}\,
\Phi_{\alpha}(u,v;y_i,x_i, {\mathbf x \setminus x_i})  \, {\mu \big( dy_i \, \big| \, \sigma_{i^c} \big)} 
\Big\} \mu(d{\mathbf x}) \nonumber \\
&=&
\int_S
\int_{\mathbb R} \Big\{
\int_{{\mathbb R} \setminus \{x_i\}} 
%I_{\{y_i \ne x_i\}}\,
\Phi_{\alpha}(u,v;y_i,x_i, {\mathbf x \setminus x_i})  \, {\mu \big( dy_i \, \big| \, \sigma_{i^c} \big)}
\Big\}
 {\mu \big( dx_i\, \big| \, \sigma_{i^c} \big) }
\, \mu(d{\mathbf x}) \nonumber \\
&=&
\int_S
\Big\{
\int_{\{y_i \ne y'_i\}} 
\Phi_{\alpha}(u,v;y_i,y'_i, {\mathbf x \setminus x_i})  \, {\mu \big( dy_i \, \big| \, \sigma_{i^c} \big)}\,
 {\mu \big( dy'_i \, \big| \, \sigma_{i^c} \big) }
\Big\} \mu(d{\mathbf x})   \\
&=&
\int_{S \setminus x_i}
\Big\{
\int_{\{y_i \ne y'_i\}} 
\Phi_{\alpha}(u,v;y_i,y'_i, {\mathbf x} \setminus x_i)  \, {\mu \big( dy_i \, \big| \, \sigma_{i^c} \big)}\,
 {\mu \big( dy'_i \, \big| \, \sigma_{i^c} \big) }
\Big\} \mu(d({\mathbf x} \setminus x_i)) \nonumber
\end{eqnarray}
where
$\displaystyle{ S \setminus x_i \equiv \{ (x_1, \dots, x_{i-1}, x_{i +1}, \dots)  :  (x_1, \dots) \in S\} }$  and 
$\mu(d({\mathbf x} \setminus x_i))$ is the 
marginal probability distribution of the variable 
${\mathbf x} \setminus x_i$ (cf. the notation used in (2.8)), i.e., for any $A \in \sigma_{i^c}$, \\
$\displaystyle{\int_A \mu(d({\mathbf x} \setminus x_i)) = \int_S I_{\mathbb {R}}(x_i) \, I_A({\mathbf x} \setminus x_i) \, \mu (d {\mathbf x})}$. 
The third and fourth formulas give the more symmetric definitions for ${\cal E}^{\alpha}_i(u,v)$ with respect to the variables $y_i$ and $x_i$ (cf. (1.2.1) of [Fukushima 80]), which will be used in section 4.

Under the adequate assumptions on the density functions $\rho$ (cf. (2.11)), by multiplying adequate 
weighted functions to $\Phi_{\alpha}$,  and then by taking the limit
 $\alpha \uparrow 2$ 
for the 
{\bf non-local}
symmetric form ${\cal E}_{(\alpha)}^{(i)}$ given  by (2.9), it is possible to define 
various 
{\it \bf local} Dirichlet forms 
 (cf., e.g., section II-2 of [M,R 92]). 
The corresponding  considerations 
will be announced in forthcoming works.
%(also for the infinite dimensional cases) 
%for the {\bf infinite} dimensional situation corresponding to the above {\bf finite} dimensional 
%corresponding to the above observation will be carried out in forthcoming work.
%\\
%\quad \quad \quad \qquad \qquad \qquad \qquad \qquad 
%\qquad \qquad  \qquad 
%\qquad
%\qquad  
%\qquad \qquad  \qquad 
\qquad 
\bsquare
\end{remark}
%}}

\medskip
The following is the main theorem of the closability part of this paper.

%{\textcolor{blue}{
\begin{theorem}
For the symmetric non-local forms ${\cal E}_{(\alpha)}$, $0 <\alpha <2$ given by (2.10) 
(for $1 < \alpha <2$ with the additional assumption (2.11) 
with (2.12)
) the following hold
(cf. Remark 1-i),ii)):
\\
i) \quad \,  ${\cal E}_{(\alpha)}$ is well-defined on ${\cal F}C^{\infty}_0$; \\
ii) \quad  
$({\cal E}_{(\alpha)},{\cal F}C^{\infty}_0)$ is 
closable in $L^2(S;\mu)$; \\
iii) \quad 
$({\cal E}_{(\alpha)},{\cal F}C^{\infty}_0)$ is 
Markovian.\\
Thus, 
for each $0 < \alpha <2$, the closed extension of $({\cal E}_{(\alpha)},{\cal F}C^{\infty}_0)$ denoted by 
$({\cal E}_{(\alpha)}, {\cal D}({\cal E}_{(\alpha)}))$ with the domain ${\cal D}({\cal E}_{(\alpha)})$,
 is a non-local Dirichlet form on $L^2(S;\mu)$.\\
%\\
%\qquad \qquad \qquad \qquad \qquad \qquad \qquad \qquad \qquad \qquad  \qquad \qquad
%\qquad  \qquad \qquad  \qquad \qquad \quad 
%\bsquare
\end{theorem}
%}}
The proof of Theorem 1 is presented in the next section.

\section{Proof of Theorem 1.}
{\bf{We start the proof of Theorem 1 for $0 < \alpha \leq 1$}}\\
For the statement i),  we have to show that \\
i-1) \quad 
for any real valued ${\cal B}(S)$-measurable function $u$ on $S$,  such that\\
 $u=0, \, \mu-a.e.$, it holds that\quad 
$\displaystyle{{\cal E}_{(\alpha)}(u,u) =0}$
(cf. (3.8) given below), and \\
%}}
i-2) \quad for any $u, v \in {\cal F}C^{\infty}_0$, \, there corresponds only one value ${\cal E}_{(\alpha)}(u,v) \in {\mathbb R}$, \\
For the statement ii), we have to show the following: \, For a sequence $\{u_n\}_{n \in {\mathbb N}}$, 
$u_n \in {\cal F}C^{\infty}_0$, $n \in {\mathbb N}$, if 
\begin{equation}
\lim_{n \to \infty} \| u_n \|_{L^2(S; \mu)} = 0,
\end{equation}
and
\begin{equation}
\lim_{n,m \to \infty} {\cal E}_{(\alpha)}(u_n -u_m, u_n -u_m) =0,
\end{equation}
then
\begin{equation}
\lim_{n \to \infty} {\cal E}_{(\alpha)} (u_n,u_n) =0.
\end{equation}
For the statement iii), we have to show that (cf. [Fukushima 80] and Proposition I-4.10 of [M,R 92]) for any $\epsilon >0$ there exists a real function $\varphi_{\epsilon}(t)$, $ -\infty < t < \infty$, such that ${\varphi}_{\epsilon}(t) = t, \, \forall t \in [0,1]$, 
$-\epsilon \leq {\varphi}_{\epsilon}(t) \leq 1 + \epsilon, \, \forall t \in (-\infty, \infty)$, and
$0 \leq {\varphi}_{\epsilon}(t')-{\varphi}_{\epsilon}(t) \leq t'-t$ for $t < t'$, such that for any $u \in {\cal F}C^{\infty}_0$ it holds that ${\varphi}_{\epsilon}(u) \in {\cal F}C^{\infty}_0$ and
\begin{equation}
{\cal E}_{(\alpha)}({\varphi}_{\epsilon}(u), {\varphi}_{\epsilon}(u)) \leq {\cal E}_{(\alpha)}(u,u).
\end{equation}
%To prove iv) we use Lemma 1 prepared in the next section, a general statement on the elements of 
%domain ${\cal D}({\cal E}_{(\alpha)})$ (see below). 

i-1) can be seen as follows: \\
%{\textcolor{blue}{
For each $i \in {\mathbb N}$ and any
real valued ${\cal B}(S)$-measurable function $u$, note that for each $\epsilon >0$, 
%and any compact subset $K$ of ${\mathbb R}$, 
$$  I_{\{\epsilon < |x_i-y_i|\} }(y_i) \,
%\, I_K(y_i)
\Phi_{\alpha}(u,u;y_i, x_i, {\mathbf x} \setminus x_i)$$
defines a ${\cal B}(S \times {\mathbb R})$-measurable function.
The function 
$\Phi_{\alpha}(u,u;y_i,x_i, {\mathbf x} \setminus x_i)$, is defined by setting 
$v =u, \, x =x_i$, 
in (2.8).   
 ${\cal B}(S \times {\mathbb R})$ is the Borel $\sigma$-field of $S \times {\mathbb R}$.
${\mathbf x} = (x_i, i \in {\mathbb N}) \in S$ and $y_i \in {\mathbb R}$.
Then, 
%for any compact subset $K$ of ${\mathbb R}$, 
$0 \leq I_{\{\epsilon < |x_i-y_i|\} }(y_i) \,
%\, I_K(y_i)
\Phi_{\alpha}(u,u;y_i, x_i, {\mathbf x} \setminus x_i)$ converges monotonically to
$ I_{\{y_i \ne x_i\}}(y_i) \,
\Phi_{\alpha}(u,u;y_i,x_i, {\mathbf x} \setminus x_i)$ as 
%$ K \uparrow {\mathbb R}$ and 
$\epsilon \downarrow 0$, for every $y_i \in {\mathbb R}$, \, ${\mathbf x} \in S$, 
and by 
Fatou's Lemma, 
we have 
\begin{eqnarray}
\lefteqn{
\int_S  
\left\{
\int_{\mathbb R} I_{\{y_i \ne x_i\}}(y_i) \,
\Phi_{\alpha}(u,u;y_i,x_i, {\mathbf x} \setminus x_i)  \, {\mu \big( dy_i \, \big| \, \sigma_{i^c} \big) }
\right\} \mu(d{\mathbf x})
}
  \\
&& = 
\int_S  
%\liminf_{K \uparrow {\mathbb R}} 
\liminf_{\epsilon \downarrow 0} 
 \left\{ 
 \int_{\mathbb R}  I_{\{\epsilon < |x_i-y_i|\} }(y_i) 
%\, I_K(y_i)
\,
\Phi_{\alpha}(u,u;y_i, x_i, {\mathbf x} \setminus x_i)  \, {\mu \big( dy_i \, \big| \, \sigma_{i^c} \big) }
\right\} \mu(d{\mathbf x}) 
\nonumber \\
&& \leq 
%\liminf_{K \uparrow {\mathbb R}} 
\liminf_{\epsilon \downarrow 0}
\int_S   
 \left\{ 
 \int_{\mathbb R}  I_{\{\epsilon < |x_i-y_i| \}}(y_i) 
%\, I_K(y_i) 
\,
\Phi_{\alpha}(u,u;y_i, x_i, {\mathbf x} \setminus x_i)  \, {\mu \big( dy_i \, \big| \, \sigma_{i^c} \big) }
\right\} \mu(d{\mathbf x}). \nonumber
\end{eqnarray}
%where $K$ denotes a compact set of  ${\mathbb R}$ 
%and $I_K$ the corresponding indicator function.
 Through the definition of the conditional probability distributions and conditional expectations, we see that, 
for any $\epsilon >0$,  
\begin{eqnarray}
\lefteqn{
\int_S \left\{ \int_{\mathbb R}  I_{\{\epsilon < |x_i-y_i| \}}(y_i)  
%\, I_K(y_i) 
\, 
\frac{1}{|y_i -x_i|^{ \alpha + 1}} \big(u(x_1, \dots, x_{i-1}, y_i, x_{i +1}, \dots) \big)^2 
\, {\mu \big( dy_i \, \big| \, \sigma_{i^c} \big) }
\right\} \mu(d{\mathbf x})
}  \nonumber \\
&&\leq 
\frac{1}{{\epsilon}^{ \alpha +1}} 
\int_S \left\{ \int_{\mathbb R}  I_{\{\epsilon < |x_i-y_i| \}}(y_i) 
%\, I_K(y_i) 
\, 
 \big(u(x_1, \dots, x_{i-1}, y_i, x_{i +1}, \dots) \big)^2 
\, {\mu \big( dy_i \, \big| \, \sigma_{i^c} \big) }
\right\} \mu(d{\mathbf x}) 
 \nonumber \\
&&\leq 
\frac{1}{{\epsilon}^{ \alpha +1}} 
\int_S \left\{ \int_{\mathbb R}  
 \big(u(x_1, \dots, x_{i-1}, y_i, x_{i +1}, \dots) \big)^2 
\, {\mu \big( dy_i \, \big| \,  \sigma_{i^c} \big) }
\right\} \mu(d{\mathbf x}) 
 \nonumber \\
&&=
\frac{1}{{\epsilon}^{ \alpha +1}} 
\int_S 
 \big(u(x_1, \dots, x_{i-1}, x_i, x_{i +1}, \dots) \big)^2 
 \mu(d{\mathbf x}),
\end{eqnarray}
and 
\begin{equation}
\int_S 
\big(u(x_1, \dots) \big)^2 
\left\{ \int_{\mathbb R}  I_{\{\epsilon < |x_i-y_i| \}}(y_i) 
%\, I_K(y_i) 
\, 
\frac{1}{|y_i -x_i|^{ \alpha + 1}} 
\, {\mu \big( dy_i \, \big| \, \sigma_{i^c} \big) }
\right\} \mu(d{\mathbf x})
 \nonumber 
\end{equation}
\begin{equation}
\leq
\frac{1}{{\epsilon}^{ \alpha +1}} 
\int_S 
 \big(u(x_1,\dots) \big)^2 
 \mu(d{\mathbf x}).
\end{equation}
From (3.6), by making use of the Cauchy Schwarz inequality
\begin{eqnarray}
\lefteqn{ 
\Big|
\int_S 
u(x_1, \dots) 
\Big\{ \int_{\mathbb R}  I_{\{\epsilon < |x_i-y_i| \}}(y_i) 
%\, I_K(y_i) 
\, 
\frac{1}{|y_i -x_i|^{ \alpha + 1}} 
} \nonumber \\
&&\times 
u(x_1, \dots, x_{i-1}, y_i, x_{i +1}, \dots)  
\, {\mu \big( dy_i \, \big| \, \sigma_{i^c} \big) }
\Big\} \mu(d{\mathbf x}) \Big| 
\nonumber \\
&&\leq
\frac{1}{{\epsilon}^{\alpha +1}} 
\int_S 
 \big(u(x_1,\dots) \big)^2 
 \mu(d{\mathbf x}). \nonumber
\end{eqnarray}
%{\textcolor{blue}{
By making use of 
this,  (3.6) and (3.7), from (3.5) 
it follows the property  i-1): \,
for all $i \in {\mathbb N}$,
%\begin{equation*}
${\cal E}_{(\alpha)}^{(i)} (u,u) = 0$, %\quad \forall i \in {\mathbb N}$, 
%\end{equation*}
for any real valued ${\cal B}(S)$-measurable function $u$ such that $u=0$, 
\, $\mu$-a.e..
%\begin{equation}
 %{\mbox{for any real valued ${\cal B}(S)$-measurable function $u$ such that $u=0$, 
%\, $\mu$-a.e.}}.
%\end{equation}
%}}
This then implies 
\begin{equation}
{\cal E}_{(\alpha)} (u,u) =0 \qquad {\mbox{ for all such $u$}}.
\end{equation}
%}}

In order to prove i-2), 
 {%\textcolor{red}
{for $0 < \alpha \leq 1$}},  
take any 
{\it{representative}} 
${\tilde{u}} \in C^{\infty}_0({\mathbb R}^n)$ of $u \in {\cal F}C^{\infty}_0$, $n \in {\mathbb N}$ (precisely, ${\tilde{u}}$ is defined as (3.15)) . 
 {%\textcolor{red}
Since ${\tilde{u}}$ is a {\it{Lipschitz}} function, it is easy to see  from the definition (2.8) that
(cf. the formula on $|\eta_{M,i}(x) -\eta_{M,i}(y)|^2$ given above (4.33), in the next section)
there exists an $M < \infty$ depending on ${\tilde{u}}$  such that
$$
0 \leq \Phi_{\alpha}({\tilde{u}},{\tilde{u}};y_i,y'_i,{\mathbf x} \setminus x_i) \leq M,
 \quad {\mbox{$ \forall {\mathbf x} \in S$,  and $\forall y_i, \, y'_i \in {\mathbb R}$}}.
$$
By this, we have that 
\begin{equation}
{\cal E}_{(\alpha)}^{(i)}({\tilde{u}},{\tilde{u}}) \in {\mathbb R}
\end{equation}
(in fact, for only a finite number of $i \in {\mathbb N}$. 
we have 
${\cal E}_{(\alpha)}^{(i)}({\tilde{u}},{\tilde{u}})  \ne 0$, cf. also (2.10)). 
Since, $u = {\tilde u} + {\overline{0}}$ for some 
real valued ${\cal B}(S)$-measurable function ${\overline{0}}$ such that ${\overline{0}}=0$, $\mu$-a.e., 
by (3.9) together with i-1) (cf. (3.8)) and  the Cauchy-Schwarz inequality, for $u \in {\cal F}C^{\infty}_0$, 
${\cal E}_{(\alpha)}(u,u)
= {\cal E}_{(\alpha)}({\tilde{u}},{\tilde{u}}) 
 \in {\mathbb R}$, $0 < 
\alpha \leq 1$. 
Thus ${\cal E}_{(\alpha)}(u,u)$ is {\it{uniquely defined}}.
Then by the Cauchy-Schwarz inequality i-2) follows.

ii) can be proved as follows (cf. section 1 of [Fukushima 80]): \, Suppose that a sequence $\{u_n \}_{n \in {\mathbb N}}$ satisfies (3.1) and (3.2).  Then, by (3.1) there exists a measurable set ${\cal N} \in {\cal B}(S)$ and a 
{%\textcolor{blue}
{
sub sequence $\{ u_{n_k} \}$ of $\{u_n\}$ such that
$$
\mu({\cal N}) = 0, \qquad   \lim_{{n_k} \to \infty} u_{n_k} (\mathbf x) = 0, \qquad \forall {\mathbf x} \in S \setminus {\cal N}.
$$
Define
$$
{\tilde{u}}_{n_k} ({\mathbf x}) = u_{n_k}({\mathbf x}) \quad {\mbox{for ${\mathbf x} \in S \setminus {\cal N}$, \qquad and}} \qquad 
{\tilde{u}}_{n_k} ({\mathbf x}) = 0 \quad {\mbox{for ${\mathbf x} \in  {\cal N}$}}.
$$
Then, 
\begin{equation} {\tilde{u}}_{n_k} (\mathbf x) = u_{n_k}(\mathbf x),  \, \, \mu-a.e., 
\qquad  \lim_{n_k \to \infty} {\tilde{u}}_{n_k} (\mathbf x) = 0, \quad \forall {\mathbf x} \in S.
\end{equation}
}}
By the fact  i-1), 
precisely by  (3.8) 
 shown above and (3.10), for each $i$, we see that
\begin{eqnarray}
\lefteqn{
\int_S  
\Big\{
 \int_{\mathbb R} I_{\{y_i \ne x_i\}}(y_i) \,
\Phi_{\alpha}(u_n,u_n;y_i,x_i, {\mathbf x} \setminus x_i)  \, {\mu \big( dy_i \, \big| \, \sigma_{i^c} \big) }
\Big\} \mu(d{\mathbf x})
}
 \nonumber \\
&&=
\int_S  
 \Big\{ 
\int_{\mathbb R} I_{\{y_i \ne x_i\}}(y_i) \,
\lim_{n_k \to \infty}
\Phi_{\alpha}(u_n - {\tilde{u}}_{n_k},u_n- {\tilde{u}}_{n_k};y_i,x_i, {\mathbf x} \setminus x_i)  \, {\mu \big( dy_i \, \big| \, \sigma_{i^c} \big) }
\Big\} \mu(d{\mathbf x})
 \nonumber \\
&&\leq
\liminf_{n_k \to \infty}
\int_S 
\Big\{ 
\int_{\mathbb R} I_{\{y_i \ne x_i\}} \,
\Phi_{\alpha}(u_n - {\tilde{u}}_{n_k},u_n- {\tilde{u}}_{n_k};y_i, x_i, {\mathbf x} \setminus x_i)  \, {\mu \big( dy_i \, \big| \, \sigma_{i^c} \big) }
\Big\} \mu(d{\mathbf x})
 \nonumber \\
&&=
\liminf_{n_k \to \infty}
\int_S  \, 
 \Big\{ 
 \int_{\mathbb R} I_{\{y_i \ne x_i\}}\,
\Phi_{\alpha}(u_n - {{u}}_{n_k},u_n- {{u}}_{n_k};y_i, x_i, {\mathbf x} \setminus x_i)  \, {\mu \big( dy_i \, \big| \, \sigma_{i^c} \big) }
\Big\} \mu(d{\mathbf x})
 \nonumber \\
&&\equiv
\liminf_{n_k \to \infty}
{\cal E}^{(i)}_{(\alpha)}(u_n -u_{n_k}, u_n -u_{n_k}).
\end{eqnarray}
Now, by applying the assumption (3.2) to the right hand side of (3.11), we get 
\begin{equation}
\lim_{n \to \infty} {\cal E}_{(\alpha)}^{(i)}(u_n,u_n) = 0, \qquad \forall i \in {\mathbb N}.
\end{equation}
%{\textcolor{blue}{
(3.12) together with i) shows that for each $i \in {\mathbb N}$, ${\cal E}^{(i)}_{(\alpha)}$ with the domain ${\cal F}C_0^{\infty}$ is closable in $L^2(S;\mu)$.
Since, ${\cal E}_{(\alpha)} \equiv \sum_{i \in {\mathbb N}} {\cal E}^{(i)}_{(\alpha)}$, by using Fatou's Lemma, from (3.12) and the assumption (3.2) we see that 
$$
{\cal E}_{(\alpha)}(u_n,u_n) = \sum_{i \in {\mathbb N}} \lim_{m \to \infty} {\cal E}^{(i)}_{(\alpha)}(u_n -u_m, u_n -u_m) 
\leq \liminf_{m \to \infty}  {\cal E}_{(\alpha)}(u_n -u_m, u_n -u_m)  \to 0,
$$
as $n \to 
\infty$.
This proves  (3.3) 
(cf.  Proposition I-3.7 of [M,R 92] for a general argument). The proof of ii) is thus complete.
%}}

iii) can be shown as follows: \, 
For each $\epsilon >0$, take a smooth function ${\varphi}_{\epsilon}(t)$ such that
$\varphi_{\epsilon}(t) = t$ for $t \in [- \frac{\epsilon}2, 1+ \frac{\epsilon}2]$,  
$\varphi_{\epsilon}(t) = - \epsilon$ for $t \leq - 2 \epsilon$,
$\varphi_{\epsilon}(t) = 1 + \epsilon$ for $t \geq 1+ 2 \epsilon$,
and it satisfies $0 \leq {\varphi}_{\epsilon}(t') - \varphi_{\epsilon}(t) \leq t'-t$ for $t < t'$).
Then, for $u \in {\cal F}C_0^{\infty} \subset {\cal D}({\cal E}_{(\alpha)})$, it holds that 
${\varphi}_{\epsilon}(u({\mathbf x}))  \in {\cal F}C_0^{\infty} \subset {\cal D}({\cal E}_{(\alpha)})$,
and  (3.4) is satisfied (cf. section 1 of Fukushima 80]).

For $0 < \alpha \leq 1$ the proof of Theorem 1 is therefore complete.

{\bf{The proof of Theorem 1,  for $1< \alpha < 2$}}

The proof of i-1), ii) and  iii)   can be carried out by  
 exactly the 
same manner as the previous proof 
we have provided 
for the case $0 < \alpha \leq 1$. 
We  only  show that  i-2), i.e., ${\cal E}_{(\alpha)}(u,u) < \infty, \, \forall u \in {\cal F}C^{\infty}_0$ 
also holds when we make use of 
%{\textcolor{blue}{
 the additional assumption (2.11) with (2.12), 
i.e., of
 existence of 
a {\it{locally bounded probability density}} (cf. (2.12)),
$\rho(y \, \big| \, \sigma_{i^c})$ for 
$\mu(dy \, \big| \, \sigma_{i^c})$, i.e., 
 $\mu(dy \, \big| \, \sigma_{i^c}) = 
\rho(y \, \big| \, \sigma_{i^c}) \, dy, \, \, \mu-a.e.$.
%}}
Under this assumption by applying  Young's inequality, we  derive i-2) 
(cf. the formula given above (3.9), also cf.  (3.6), (3.7) and (3.8)) as follows: 
For $1 < \alpha <2$, by the definition (2.8), we note that  
\begin{equation}
\Phi_{\alpha}(u,u; y_i, x_i,{\mathbf x} \setminus x_i) = \Phi_{1}(u,u; y_i, x_i, {\mathbf x} \setminus x_i) \cdot \frac{1}{|x_i- y_i|^{ \alpha -1}}.
\end{equation}
%{\textcolor{blue}
%{
Take a
{\it{representative}} 
${\tilde{u}} \in C^{\infty}_0({\mathbb R}^n)$ of $u \in {\cal F}C^{\infty}_0$, given by  (3.15) below, and 
%For each $u \in {\cal F}C^{\infty}_0$, 
let
% $\displaystyle{\frac{I_{K_i}}{|x|^{2 \alpha -1}}  \in L^1({\mathbb R}^1)}$ holds, where  
$I_{K_i}$ 
be
 the indicator function for a compact set $K_i$ such that
\begin{equation}
K_i \equiv \{ x \in {\mathbb R} \, : \, x = x' - x'', \, \, x',x'' \in U_i\},
\end{equation}
with
\begin{equation*}
U_i \equiv {\mbox{the closure of}} \, \{ x_i \in {\mathbb R} \, : \, f(x_1, \dots, x_{i-1},x_i, x_{i+1}, \dots x_n) \ne 0 \}, \, i=1, \dots, n,
\end{equation*}
where
\begin{equation}
f \in C^{\infty}_0({\mathbb R}^n \to {\mathbb R}), \quad {\mbox{such that}} \quad 
{\tilde{u}}({\mathbf{x}}) = f \cdot \prod_{j \in {\mathbb N}} I_{\mathbb R}(x_j)  \, \in {\cal F}C^{\infty}_0.
\end{equation}
Also by the assumption (2.12), the probability density satisfies the following:
\, for each $i \in {\mathbb N}$ and the compact $U_i$,  there exists an $0< L_{U_i,i }<\infty$, 
and for $a.e. \, y \in {\mathbb R}$, with respect to the Lebesgue measure, 
$$
I_{U_i,i}(y)  \rho(y \, \big| \, \sigma_{i^c}) \leq L_{U_i,i}, \quad  \mu-a.e..$$
Hence 
by  Young's inequality (for the 
convolutions) we have the following bound: \,
for \, $ i=1, \dots, n$, 
%\begin{eqnarray}
%\lefteqn{
\begin{equation}
%{\textcolor{red}{
\sup_{x \in {\mathbb R}} | \int_{\mathbb R} \frac{I_{K_i}(x-y)}{|x-y|^{\alpha-1}} I_{U_i}(y) \rho(y |{\sigma}^c_i) dy |
%}}
 %\|\frac{I_{K_i}}{|\cdot|^{ \alpha -1}} \ast \rho(\cdot \, \big| \, \sigma_{i^c})\|_{L^{\infty}({\mathbb R})} 
%} 
%\\
%&& \leq  \| \frac{I_{K_i}}{|\cdot|^{ \alpha -1}}\|_{L^1({\mathbb R})} \, 
%\|I_{K_i}(\cdot) \, \rho(\cdot \, \big| \, \sigma_{i^c}) \|_{L^{\infty}({\mathbb R})}  
= L_{U_i,i} \, \| \frac{I_{K_i}}{|\cdot|^{ \alpha -1}}\|_{L^1({\mathbb R})}, \quad  \mu-a.e..
% \nonumber
%\end{eqnarray}
\end{equation}
For  notational simplicity, let  $f(x_1, \dots, x_{i-1},x_i, x_{i+1}, \dots, x_n)$ be denoted by $f(\cdot, x_i, \cdot)$.
Since,   $f(\cdot, x_i, \cdot) = f(\cdot, y_i, \cdot)=0$ \, for $x_i \in U^c_i$ and $y_i \in U^c_i$, 
we note that the following is true:
$$(f(\cdot, x_i, \cdot) - f(\cdot, y_i, \cdot)) I_{U_i}^c(x_i)  I_{U_i}^c(y_i)
 = 0.$$
We then have the following evaluation for $x_i \ne y_i$:
\begin{eqnarray}
\lefteqn{
I_{K_i}(x_i -y_i) \frac{(f(\cdot, x_i, \cdot) - f(\cdot, y_i, \cdot))^2}{|x_i -y_i|^{\alpha +1} }
} \nonumber \\
&&= I_{K_i}(x_i -y_i) \frac{(f(\cdot, x_i, \cdot) - f(\cdot, y_i, \cdot))^2}{|x_i -y_i|^{\alpha +1} }
\nonumber \\
&&\times \{ I_{U_i}(x_i)I_{U_i}(y_i) + I_{U_i^c}(x_i)I_{U_i}(y_i) +I_{U_i}(x_i)I_{U_i^c}(y_i) +
I_{U_i^c}(x_i)I_{U_i^c}(y_i) \}
\nonumber \\
&& = I_{K_i}(x_i -y_i) \frac{(f(\cdot, x_i, \cdot) - f(\cdot, y_i, \cdot))^2}{|x_i -y_i|^{2} }
\frac{1}{{|x_i -y_i|^{\alpha -1} }}
\nonumber \\
&&\times \{ I_{U_i}(x_i)I_{U_i}(y_i) + I_{U_i^c}(x_i)I_{U_i}(y_i) +I_{U_i}(x_i)I_{U_i^c}(y_i) 
\}
\nonumber \\
&& \leq (\sup_{0 < \theta <1} |f_{x_i}(\cdot, x_i + \theta (y_i - x_i), \cdot)|^2) 
\frac{I_{K_i}(x_i -y_i)}{{|x_i -y_i|^{\alpha -1} }}
\nonumber \\
&&\times \{ I_{U_i}(x_i)I_{U_i}(y_i) + I_{U_i^c}(x_i)I_{U_i}(y_i) +I_{U_i}(x_i)I_{U_i^c}(y_i) \},
\end{eqnarray}
in the above and below we denoted $\frac{\partial f}{\partial x_i} = f_{x_i}$.
From (3.17), by making use of (3.16) we have that
\begin{eqnarray}
\lefteqn{
\int_S \{
\int_{x_i \ne y_i} 
I_{K_i}(x_i -y_i) \frac{(f(\cdot, x_i, \cdot) - f(\cdot, y_i, \cdot))^2}{|x_i -y_i|^{\alpha +1} }
\rho(x_i|{\sigma}^c_i) \rho(y_i|{\sigma}^c_i) dx_i dy_i \}
\mu(d{\mathbf x})
}
\nonumber \\
&& \leq 
 3(\sup_{0 < \theta <1} |f_{x_i}(\cdot, x_i + \theta (y_i - x_i), \cdot)|^2) 
\{ 
\int_S \{
\int_{\mathbb R} 
L_{U_i,i} \, \| \frac{I_{K_i}}{|\cdot|^{ \alpha -1}}\|_{L^1({\mathbb R})}
\rho(x_i|{\sigma}^c_i) dx_i \}
\mu(d {\mathbf x}) \}
\nonumber \\
&&= 3(\sup_{0 < \theta <1} |f_{x_i}(\cdot, x_i + \theta (y_i - x_i), \cdot)|^2) 
L_{U_i,i} \, \| \frac{I_{K_i}}{|\cdot|^{ \alpha -1}}\|_{L^1({\mathbb R})}.
\end{eqnarray}

Moreover, 
for ${\tilde{u}} \in {\cal F}C^{\infty}_0$ given by (3.15) with $f$ satisfying $f \ne 0$,  
since $f \in C^{\infty}_0({\mathbb R}^n \to {\mathbb R})$ for some $n \in {\mathbb N}$, 
for each 
$i=1, \dots, n$, 
there exists $m_i > 0$ and it holds that 
${\displaystyle{K_i^c \subset \{ x \, : \, |x | \geq m_i \} \subset {\mathbb R}}}$.
Hence for each $1 < \alpha <2$ and $i = 1, \dots, n$, there exists an $M'_i < \infty$ such that 
\begin{equation}
\frac{1}{|x_i -y_i|^{ \alpha -1}} \, I_{K_i^c}(x_i -y_i) \leq M'_i,  \qquad \forall x_i, \, y_i \in {\mathbb R}.
\end{equation}
%there exists an $M'_i < \infty$ such that 
%$$
%\frac{1}{|x_i -y_i|^{ \alpha -1}} I_{K_i^c}(x_i -y_i) I_{U_i}(x_i) I_{U_i}(y_i) \leq M'_i, 
%$$
%$$
%\frac{1}{|x_i -y_i|^{\alpha -1}} I_{K_i^c}(x_i -y_i) I_{U_i}(x_i) I_{U_i^c}(y_i) \leq M'_i
%$$
%and 
%\begin{equation}
%\frac{1}{|x_i -y_i|^{ \alpha -1}} I_{K_i^c}(x_i -y_i) I_{U_i^c}(x_i) I_{U_i^c}(y_i) = 0, \qquad 
%\forall x_i, \, y_i \in {\mathbb R}, \quad i =1, \dots, n,
%\end{equation}
In  addition, by the definition (2.8),  for each ${\tilde{u}} \in {\cal F}C^{\infty}_0$ given by (3.15),  there exists an $M< \infty$ (cf. the formula given above (3.9)) and 
\begin{equation*}
0 \leq 
\Phi_{1}({\tilde{u}},{\tilde{u}}; y_i, y'_i, {\mathbf x} \setminus x_i) \leq M, 
\quad {\mbox{$ \forall {\mathbf x} \in S$,  and $\forall y_i, \, y'_i \in {\mathbb R}$, $i =1, \dots, n$}}.
\end{equation*}
By making use of  this and (3.19),  with (3.13),  
we also have that (cf. (3.18)) 
\begin{eqnarray}
\int_S \{
\int_{x_i \ne y_i} 
I_{K_i^c}(x_i \! &-& \!y_i) \frac{(f(\cdot, x_i, \cdot) - f(\cdot, y_i, \cdot))^2}{|x_i -y_i|^{\alpha +1} }
\rho(x_i|{\sigma}^c_i) \rho(y_i|{\sigma}^c_i) dx_i dy_i \}
\mu(d{\mathbf x})
\nonumber
\\
%\nonumber \\
& \leq &
 M \cdot M'_i.
\end{eqnarray}
From (3.20) and (3.18), 
%{\textcolor{red}{
by the symmetric expression (2.13)
%}}
we see that
${\cal E}^{(i)}_{(\alpha)}({\tilde{u}},{\tilde{u}}) < \infty$, $\forall i \in {\mathbb N}$,   
(cf. (4.64) in section 4, for similar detailed evaluation)
and 
${\cal E}_{(\alpha)}({\tilde{u}},{\tilde{u}}) < \infty$ 
 (in fact, for only a finite number of $i \in {\mathbb N}$, 
${\cal E}^{(i)}_{(\alpha)}({\tilde{u}},{\tilde{u}})  \ne 0$, cf. also (2.10))
for any ${\tilde{u}} \in {\cal F}C^{\infty}_0$ given by (3.15). 
Since i-1) holds for $1 < \alpha <2$, 
 passing through the same argument as for the case where    $0 < \alpha \leq 1$,  we have i-2) for $1 < \alpha <2$.
%}} 
\\
\qquad \qquad \qquad \qquad \qquad 
\quad
\qquad 
\quad 
\qquad \qquad 
\qquad \qquad  
\qquad \qquad
\qquad  \qquad \qquad 
% \qquad \qquad \quad 
\bsquare

\begin{remark}{}\qquad
%i) \quad In Theorem 1, for the case where $1 < \alpha <2$,  a regularity assumption of the probability measure,(2.11) with (2.12), is assumed. By observing the proof of 
%i-2) for $1< \alpha <2$ (cf. (3.16)), we see that the assumption (2.11) with (2.12) can be substituted by the following abstract one, under which the statement of Theorem 1 for 
%$1< \alpha < 2$ is also true: \, for each compact $K \subset {\mathbb R}$, there exists an 
%$L_i < \infty$ and 
%\begin{equation}
% \sup_{y \in {\mathbb R}} I_K(y) \int_{{\mathbb R}} \frac{I_K(y')}{|y -y'|^{ \alpha -1}} \mu(dy'\,|\, 
%\sigma_{i^c})  \, \leq L_i, \quad \mu-a.e..
%\end{equation}
% \\
i) \quad The condition 
%(3.18), and hence 
 (2.12) can be substituted by a general abstract condition that is an analogue of (2.1) of [Fukushima,Uemura 2012] given for the finite dimensional cases.
\\
ii) \quad For each $0 \leq \alpha <2$, by using the formulas (2.9), if we define a quadratic form 
${\overline{\cal E}}_{(\alpha)}$ on $L^2(S;\mu)$ (cf. (2.10)) such that
\begin{equation}
{\overline{\cal E}}_{(\alpha)} (u,v) \equiv \sum_{i \in {\mathbb N}}{\cal E}^{(i)}_{(\alpha)}(u,v), \qquad 
u,v \in {\cal D}({\overline{\cal E}}_{(\alpha)}),
\end{equation}
where ${\cal D}({\overline{\cal E}}_{(\alpha)})$ is the domain of ${\overline{\cal E}}_{(\alpha)}$ defined by 
\begin{equation}
{\cal D}({\overline{\cal E}}_{(\alpha)}) \equiv \{ u \in L^2(S; \mu) \, : \, {\overline{\cal E}}_{(\alpha)}(u,u) < \infty \},
\end{equation}
then, by the same arguments as for the proof of Theorem 1, it is possible to see that 
${\overline{\cal E}}_{(\alpha)}$ with the domain ${\cal D}({\overline{\cal E}}_{(\alpha)})$ on $L^2(S; \mu)$ is a closed form. But it is not guaranteed that  ${\cal D}({\overline{\cal E}}_{(\alpha)})$ includes 
${\cal F}C^{\infty}_0$ as its dense subset, and hence 
for ${\overline{\cal E}}_{(\alpha)}$,
the Markovian property (cf. Theorem 1-ii) and  (3.1)) 
and the quasi regularity (cf. next section) can not be discussed 
through a standard argument (cf. e.g., section 1.2 of [Fukushima 80]).
Moreover it may happen that ${\cal D}({\overline{\cal E}}_{(\alpha)}) = \{0\}$, where $0$ is the $\mu$ equivalent class of measurable functions $f = 0, \, \mu-a.e.$.\\
\end{remark}

\section{Strict Quasi-regularity}
In this section, we give sufficient conditions (cf. Theorem 2 and Theorem 3 below) under which the Dirichlet forms (i.e. the closed Markovian symmetric forms)  defined in the previous section are {\it{strictly quasi-regular}} (cf., [A,R 98, 90, 91] and section IV-3 of [M,R 92]).

Denote by $({\cal E}_{(\alpha)}, {\cal D}({\cal E}_{(\alpha)}))$ the Dirichlet form on $L^2(S;\mu)$, 
with the domain ${\cal D}({\cal E}_{(\alpha)}))$ defined through Theorem 1 in the previous 
section,  obtained as  the closed extension of 
the closable Markovian symmetric form ${\cal E}_{(\alpha)}$ 
of its restriction to 
 ${\cal F}C^{\infty}_0$.   We shall use the same notation ${\cal E}_{(\alpha)}$  for the closable form and the closed form.

Recall that the state space $S$ is 
taken to be a  
Fr{\'e}chet space, either 
a weighted $l^p$ space, $l^p_{(\beta_i)}$ defined by (2.1), or  a weighted $l^{\infty}$ space, 
$l^{\infty}_{(\beta_i)}$ defined by (2.2), or  the direct product space ${\mathbb R}^{\mathbb N}$ defined by (2.3).

For each $i \in {\mathbb N}$, we denote by $X_i$ the random variable (i.e., measurable function) on $(S, {\cal B}(S), \mu)$ , that represents the coordinate $x_i$ of  ${\mathbf x} = (x_1, x_2, \dots)$, precisely,
\begin{equation}
X_i \, : S \ni {\mathbf x} \longmapsto x_i \in {\mathbb R}.
\end{equation}
By making use of the random variable $X_i$, we have the following probabilistic expression:
\begin{equation}
\int_S I_{B}(x_i) \, \mu(d {\mathbf x}) = \mu (X_i \in B), \qquad {\mbox{for \, \, $B \in 
{\cal B}(S)$}}.
\end{equation}
\begin{theorem} {}\qquad 
Let \,  $0 < \alpha \leq 1$, and 
$({\cal E}_{(\alpha)}, {\cal D}({\cal E}_{(\alpha)}))$ be the closed Markovian symmetric form defined  through Theorem 1.\\
i) \quad In the case where $S= l^p_{(\beta_i)}$, 
%{\textcolor{blue}{
for some $ 1 \leq p < \infty$, 
as defined by (2.1), 
if there exists a 
positive sequence $\{ \gamma_i \}_{i \in {\mathbb N}}$ such that 
$\sum_{i=1}^{\infty} \gamma_i^{-1} < \infty$ (i.e., $\{\gamma^{-\frac{1}p}_i \}_{i \in {\mathbb N}}$ is a positive $l^p$ sequence),
and 
 an  $M_0 \in (0,\infty)$,
%}}
  and both 
\begin{equation}
\sum_{i=1}^{\infty} 
(\beta_i \gamma_i)^{\frac{\alpha +1}p} \cdot
\mu \Big(\beta_i^{\frac{1}p} |X_i| > M_0 \cdot  \gamma_i^{-\frac{1}p}\Big) < \infty,
\end{equation}
\begin{equation}
\mu \Big( \bigcup_{M \in {\mathbb N}} \big\{ |X_i| \leq M \cdot \beta_i^{-\frac{1}{p}}\, \gamma_i^{-\frac{1}p}, \, \forall i \in  {\mathbb N} \big\}  \Big) =1,
\end{equation}
hold, then $({\cal E}_{(\alpha)}, {\cal D}({\cal E}_{(\alpha)}))$ is a  quasi-regular Dirichlet form. 
 \\
ii) \quad 
 In the case where $S= l^{\infty}_{(\beta_i)}$ as defined by (2.2), if there exists a  
sequence $\{ \gamma_i \}_{i \in {\mathbb N}}$ such that $0 < \gamma_1  \leq \gamma_2 \leq \cdots \to \infty$, and an   $M_0 \in (0,\infty)$,  and  both 
\begin{equation}
 \sum_{i=1}^{\infty} 
(\beta_i {\gamma}_i)^{\alpha +1} \cdot
\mu \Big(\beta_i |X_i| > M_0 \cdot \gamma_i^{-1} \Big) < \infty,
\end{equation}
%{\textcolor{blue}{
\begin{equation}
\mu \Big( \bigcup_{M \in {\mathbb N}} \big\{  |X_i| \leq M \cdot \beta_i^{-1}\, \gamma_i^{-1}, \, \forall i \in  {\mathbb N} \big\}  \Big) =1,
\end{equation}
%}}
hold, then $({\cal E}_{(\alpha)}, {\cal D}({\cal E}_{(\alpha)}))$ is a  quasi-regular Dirichlet form. \\
iii) \quad 
 In the case where $S= {\mathbb R}^{\mathbb N}$ as defined by (2.3), 
 $({\cal E}_{(\alpha)}, {\cal D}({\cal E}_{(\alpha)}))$ is a  quasi-regular Dirichlet form. \\
%{\textcolor{blue}{
iv) \quad 
The forms $({\cal E}_{(\alpha)}, {\cal D}({\cal E}_{(\alpha)}))$ in the statements i), ii), iii) are strictly quasi-regular Dirichlet forms.
%}}
%\\
%\qquad \qquad \qquad \qquad \qquad \qquad \qquad \qquad \qquad \qquad  \qquad \qquad
%\qquad  \qquad \qquad  \qquad \qquad \quad 
%\bsquare
\end{theorem}

\begin{remark}{}\quad i) \,  
Generally, 
for a real valued random variable $X$ on a probability space $(\Omega, {\cal F}, P)$ (denoting the expectations with respect to $P$ by $E[\cdot]$), 
for $1 \leq r < \infty$, the following Cebysev type inequality holds:
\begin{equation}
P \left( |X| > K) \right) < \frac{E[|X|^r]}{K^r}
\end{equation}
Thus, by denoting the expectations (i.e. the integrations) with respect to the measure $\mu$ by   $E_{\mu}[ \cdot]$, we see that  the following inequality is a sufficient condition for (4.3):
\begin{equation}
\sum_{i=1}^{\infty} E_{\mu}[|X_i|^2] \cdot (\beta_i \, {\gamma}_i)^{\frac{2(\alpha +1)}p} 
< \infty.
\end{equation}
%\begin{equation}
%\sum_{i=1}^{\infty} E_{\mu}[|X_i|] \cdot (\beta_i^{\frac{1}{p}}\, {\gamma}_i^{\frac{1}{p}} )^3  < \infty.
%\end{equation}
Similarly, by (4.7), 
we see that the following inequality is a sufficient condition for (4.5):
\begin{equation}
\sum_{i=1}^{\infty} E_{\mu}[|X_i|^2] \, (\beta_i \, \gamma_i)^{2(\alpha +1)} < \infty.
\end{equation}
%\begin{equation}
%\sum_{i=1}^{\infty} E_{\mu}[|X_i|] \, \beta_i^3 \, \gamma_i^3 < \infty.
%\end{equation}
ii) \, 
In the case where $S= l^2_{({\beta}_i)}$,
if the probability mesure $\mu$ is constructed through the Bochner-Minlos theorem, then  the support property (4.4) of $\mu$ can be discussed as  part of the statements of the same theorem (cf. Theorem 5 and Examples in section 5).\\
iii) \, 
To prove Theorem 2, for each $f \in C^{\infty}_0({\mathbb R}^n \to {\mathbb R})$ (for some $n \in {\mathbb N}$), the functions 
$u(\mathbf x) \equiv f(x_1, \dots, x_n) \prod_{i \geq 1} I_{\mathbb R}(x_i)$, 
%{\textcolor{blue}{
$f_M \in {\cal D}({\cal E}_{(\alpha)})$ (see (4.43) below), which is  
a function having a compact support in $S$
defined by (4.27) below, 
and $f_{M,k}$ defined by (4.36) below, which is a continuos function on $S$, play the crucial roles. For the infinite dimensional topological vector spaces $S$, any compact set in $S$ can not have an interior (cf., e.g., 
Theorem 9.2 in [Tr{\`e}ves 67]), 
and hence 
$f_M$ is not a continuous function on $S$, but under the condition (4.4) or (4.6) it will be 
 proven  that 
$f_M \in {\cal D}({\cal E}_{(\alpha)})$
%}}
 and $f_M$ approximates $u$ arbitraly 
with respect to the norm 
$\|\cdot\|_{L^2(S; \mu)} + \sqrt{{\cal E}_{(\alpha)}}$  
as $M \to \infty$ 
(see (4.43) and (4.46)).
%\qquad \qquad \qquad \qquad \qquad \qquad \qquad \qquad \qquad \qquad  \qquad \qquad
%\qquad  \qquad \qquad  \qquad \qquad \quad 
%\bsquare
\end{remark}

We prepare the following Lemma 1, that is the Lemma 2.12 in section I-2 of [M,R 92],  by which 
the proof of Theorem 2 will follow. Here we quote it in a simplified way that is adapted to the present paper:
%{\textcolor{blue}{
\begin{lemma}
For the closed form $({\cal E}_{(\alpha)}, {\cal D}({\cal E}_{(\alpha)})$ with the domain ${\cal D}({\cal E}_{(\alpha)})$ that is the closure of ${\cal F}C^{\infty}_0$, a dense subset of $L^2(S;\mu)$  (cf. Remark 1), defined by Theorem 1, the following holds: \, If a sequence $\{u_n\}_{n \in {\mathbb N}}$, $u_n  \in 
{\cal D}({\cal E}_{(\alpha)}), \, n \in {\mathbb N}$ satisfies 
\begin{equation*}
\sup_{n \in {\mathbb N}} {\cal E}_{(\alpha)}(u_n,u_n) < \infty,
\qquad {\mbox{and}} \qquad \lim_{n \to \infty} u_n = u, \, \, {\mbox{in $L^2(S;\mu)$}},
\end{equation*}
then 
\begin{equation}
 u \in {\cal D}({\cal E}_{(\alpha)}) \qquad 
{\mbox{and}} \qquad 
{\cal E}_{(\alpha)}(u,u) \leq \liminf_{n \to \infty} {\cal E}_{(\alpha)}(u_n,u_n).
\end{equation}
Moreover, there exists a subsequence $\{u_{n_k}\}_{k \in {\mathbb N}}$ of  $\{u_n\}_{n \in {\mathbb N}}$ such that its Ces{\`a}ro mean 
\begin{equation}
w_n \equiv   \frac{1}n \sum_{k=1}^n u_{n_k}  \to u  \quad  {\mbox{in}}   \,\,  {\cal D}({\cal E}_{(\alpha)}) 
\end{equation}
$$
{\mbox{with respect to the norm}} \,\, 
\|\cdot\|_{L^2(S; \mu)} + \sqrt{{\cal E}_{(\alpha)}} \, \, {\mbox{as $ n \to \infty$}}.
$$
%\qquad \qquad \qquad \qquad \qquad \qquad \qquad \qquad \qquad \qquad  \qquad \qquad
%\qquad  \qquad \qquad  \qquad \qquad \quad 
%\bsquare
\end{lemma}
%}}

{\it{\bf{Proof of Theorem 2.}}} \quad We have to verify that the Dirichlet forms 
$({\cal E}_{(\alpha)}, {\cal D}({\cal E}_{(\alpha)}))$ satisfy the definition of  quasi-regularity given by Definition 3.1 in section IV-3 of [M,R 92]. Namely, by using the same notions as in [M,R 92], we have to certify that the following i), ii) and iii) are satisfied by 
$({\cal E}_{(\alpha)}, {\cal D}({\cal E}_{(\alpha)}))$: \\
i) \quad There exists an ${\cal E}_{(\alpha)}$-nest $(D_M)_{M \in {\mathbb N}}$ consisting of compact sets. \\
ii) \quad There exists a subset of ${\cal D}({\cal E}_{(\alpha)})$, that is dense with respect to the norm $\displaystyle{\| \cdot \|_{L^2(S;\mu)} + \sqrt{{\cal E}_{(\alpha)}}}$,
and the elements of the subset have ${\cal E}_{(\alpha)}$-quasi continuous versions. \\
iii) \quad There exists $u_n \in {\cal D}({\cal E}_{(\alpha)})$, $n \in {\mathbb N}$, having ${\cal E}_{(\alpha)}$-quasi continuous $\mu$-versions ${\tilde{u}}_n$, $n \in {\mathbb N}$, and an ${\cal E}_{(\alpha)}$-exceptional set ${\cal N} \subset S$ such that $\{{\tilde{u}}_n \, :\, 
n \in {\mathbb N} \}$ separates the points of $S \setminus {\cal N}$.\\
Also, the fact that the quasi-regular Dirichlet form $({\cal E}_{(\alpha)}, {\cal D}({\cal E}_{(\alpha)}))$ is looked upon as a  {\it{strictly}} quasi-regular Dirichlet form can be guaranteed by showing the following (see Proposition V-2.15 of [M,R 92]):\\
iv) \quad $1 \in {\cal D}({\cal E}_{(\alpha)})$.

\bigskip

By Remark 1 and Theorem 1 in section 2, the above ii) and iii) hold for $({\cal E}_{(\alpha)}, {\cal D}({\cal E}_{(\alpha)}))$: \, since ${\cal F}C^{\infty}_0 \subset C(S \to {\mathbb R})$ by Remark 1-ii) and ${\cal D}({\cal E}_{(\alpha)})$ is the closure of ${\cal F}C^{\infty}_0$ by Theorem 1, we can take ${\cal F}C^{\infty}_0$ as the subset of ${\cal D}({\cal E}_{(\alpha)})$ mentioned in the above ii), also since  ${\cal F}C^{\infty}_0$ separates the points $S$,  we see that the above iii) holds.

Hence, we have only to show that the above i) and iv)  hold for $({\cal E}_{(\alpha)}, {\cal D}({\cal E}_{(\alpha)}))$.

We start the proof of i). 
 Equivalently (cf. Definition 2.1. in section III-2 of [M,R 92]), we have  to show that there exists an increasing sequence $(D_M)_{M \in {\mathbb N}}$ of compact subsets of $S$ such that $\displaystyle{
\cup_{m \geq 1} {\cal D}({\cal E}_{(\alpha)})_{D_M}}$ is dense in ${\cal D}({\cal E}_{(\alpha)})$ 
(with respect to the norm  $\displaystyle{\| \cdot \|_{L^2(S;\mu)} + \sqrt{{\cal E}_{(\alpha)}}})$, where ${\cal D}({\cal E}_{(\alpha)})_{D_M}$ is the subspace of ${\cal D}({\cal E}_{(\alpha)})$ the elements of which are functions with supports contained in  $D_M$.  For this, 
by Theorem 1, since ${\cal D}({\cal E}_{(\alpha)})$ is the closure of ${\cal F}C^{\infty}_0$, it suffices to show the following: \,
there exists a sequence of compact sets 
\begin{equation}
D_M \subset S, \quad M \in {\mathbb N}
\end{equation}
and a subset \, ${\tilde{\cal D}}({\cal E}_{(\alpha)}) \subset L^2(S;\mu)$ \, 
that satisfies
\begin{equation}
{\tilde{\cal D}}({\cal E}_{(\alpha)}) \subset
\bigcup_{M \geq 1} {\cal D}({\cal E}_{(\alpha)})_{D_M};
\end{equation}
and, moreover,
for any $u \in {\cal F}C^{\infty}_0$ there exists a sequence $\{ u_n \}_{n \in {\mathbb N}}$, $u_n \in {\tilde{\cal D}}({\cal E}_{(\alpha)})$, $n \in {\mathbb N}$, such that
\begin{equation}
\lim_{n \to \infty} u_n = u 
\quad {\mbox{in}} \quad   {\cal D}({\cal E}_{(\alpha)}),
\end{equation}
with respect to the norm  \,
$\| \cdot \|_{L^2(S;\mu)} + \sqrt{{\cal E}_{(\alpha)}}$.

{\it{\bf{Proof of (4.13) and (4.14) for $S = l^p_{(\beta_i)}, \, 1 \leq p <\infty.$}}} \quad We start from the proof of (4.13) and (4.14) by
considering
a 
suitable 
system 
$D_M$, $M \in {\mathbb N}$, 
of  compact sets 
  Let $1 \leq p < \infty$ be fixed.   For each $M \in {\mathbb N}$, define 
\begin{equation}
D_M \equiv \left\{ {\mathbf x} \in l^p_{(\beta_i)} \, : \, \, {\beta}_i^{\frac{1}{p}} | x_i|  \leq 
M \cdot {\gamma}_i^{-\frac{1}p}, \,  i \in {\mathbb N} 
 \right\},
\end{equation}
then $D_M$ is a compact set in $S = l^p_{(\beta_i)}$. This 
is proven 
 through a standard argument as follows: \, Since $\{ \gamma_i^{-1} \}_{i \in {\mathbb N}}$ is assumed to be a positive $l^1$  sequence, 
for any $\epsilon >0$, there exists an $N_{M,\epsilon} \in {\mathbb N}$ and 
\begin{equation}
\sum_{i = N_{M, \epsilon} +1}^{\infty} \beta_i |x_i|^p \leq \sum_{i = N_{M, \epsilon} +1}^{\infty} M^p \cdot 
{\gamma}_i^{-1} \leq \left(\frac{\epsilon}3 \right)^p, \qquad \forall {\mathbf x} \in D_M.
\end{equation}
Also, for any ${\mathbf x} \in D_M$, it holds that (without loss of generality 
assuming that $ 0 < {\gamma}_i^{-1} \leq 1$)
\begin{equation}
\beta_i^{\frac{1}p} |x_i| \leq M, \qquad \forall i =1,\dots, N_{M, \epsilon},
\end{equation}
In view of the 
 above, we can construct an $\epsilon$-net of $D_M$ as follows: 
Take $\displaystyle{ \epsilon' \equiv \frac{\epsilon}3 \cdot (N_{M, \epsilon})^{- \frac1{p}} }$, and for each $i \leq N_{M, \epsilon}$, set  $x_{i,j}$ in order that
\begin{equation}
\beta_i^{\frac{1}p} x_{i,j} \equiv -M + \epsilon' \cdot j, \qquad j=0, 1, \dots, [ {2M} \cdot {\epsilon'}^{-1}] +1,
\end{equation}
where $[y]$ denotes the greatest integer that is not greater than $y \in {\mathbb R}$.
Define  finite elements ${\mathbf x}_{j_1, j_2, \dots, j_{N_{M, \epsilon}}}$ in   $l^p_{(\beta_i)}$ as follows:
\begin{equation}
{\mathbf x}_{j_1, j_2, \dots, j_{N_{M, \epsilon}}} \equiv (x_{1,j_1}, x_{2, j_2}, \dots x_{N_{M, \epsilon}, j_{N_{m, \epsilon}}}, 0, 0, \dots) \in l^p_{(\beta_i)}, 
\end{equation}
$$
 j_l =0, 1, \dots, [2M\cdot {\epsilon'}^{-1}] +1,  \qquad l=1, \dots, N_{M,\epsilon}.$$
On the other hand, by (4.17) and (4.18), for any 
\begin{equation}
{\mathbf x} = (x_1, x_2, \dots, x_{N_{M,\epsilon}}, x_{N_{M,\epsilon} +1}, \dots) \in D_M
\end{equation}
there exists a vector 
$ (x_{1,j_1}, x_{2, j_2}, \dots x_{N_{M, \epsilon}, j_{N_{m, \epsilon}}}) \in {\mathbb R}^{N_{M,\epsilon}}$ such that 
\begin{equation}
|\beta^{\frac{1}p}_i x_i - \beta^{\frac{1}p}_i x_{i,j} | < \epsilon', \quad \forall i = 1, \dots, N_{M, \epsilon},
%{\textcolor{blue}{
\quad {\mbox{$x_i$ as in (4.17), $x_{i,j}$ as in (4.18).}}
%}}
\end{equation}
By combining (4.16) with (4.21), for any ${\mathbf x} \in D_M$ with the expression (4.20), there exists an 
${\mathbf x}_{j_1, j_2, \dots, j_{N_{M, \epsilon}}} \in l^p_{(\beta_i)}$ defined by (4.19), 
such that the following holds for any $\epsilon >0$:
\begin{eqnarray}
\| {\mathbf x} - {\mathbf x}_{j_1, j_2, \dots, j_{N_{M, \epsilon}}} \|_{l^p_{(\beta_i)}}^p 
& < & \sum_{i=1}^{N_{M, \epsilon}} \beta_i | x_i - x_{i,j} |^p + 
\sum_{i=N_{M,\epsilon} +1}^{\infty} \beta_i |x_i|^p  \\
& \leq & \sum_{i=1}^{N_{M, \epsilon}} (\epsilon')^p + \left(\frac{\epsilon}3\right)^p  \leq 
 \left(\frac{\epsilon}3\right)^p +  \left(\frac{\epsilon}3\right)^p < \epsilon^p.
\nonumber
\end{eqnarray}
(4.22) shows that 
for any $\epsilon > 0$ there exists a finite (cf. (4.19)) open covering of $D_M$ such that
\begin{equation}
D_M \subset \bigcup_{j_1, \dots, j_{N_{M, \epsilon}}} \left\{
{\mathbf x}' \in l^p_{(\beta_i)} \, : \, \|{\mathbf x}' - {\mathbf x}_{j_1, \dots, j_{N_{m,\epsilon}}} \|_{l^p_{(\beta_i)}} < \epsilon \right\}.
\end{equation}
Hence, for the subset 
 $D_M$ of $S= l^p_{(\beta_i)}$, there exists an $\epsilon$ net and it is totally bounded. Since, obviously, $D_M$ is a closed set and since $l^p_{(\beta_i)}$ is a complete metric space (cf. Remark 1-iii)), 
by  Fr{\'e}chet's compactness criterion for complete metric spaces, 
we see that for each $M \in {\mathbb N}$, $D_M$ is a compact subset of 
$l^p_{(\beta_i)}$.

%Next, we proceed to define ${\tilde{\cal D}}({\cal E}_{(\alpha)})$ for which (4.13) and (4.14) hold. 
%Let \, $\eta(\cdot) \in C^{\infty}_0({\mathbb R})$ be a function such that \, $\eta(x) \geq 0$, \, 
%$\displaystyle{|\frac{d}{dx} \eta(x)| \leq 1}$, \, $\forall x \in {\mathbb R}$ \, and 
Next, we proceed to define ${\tilde{\cal D}}({\cal E}_{(\alpha)})$ for which (4.13) and (4.14) hold. 
Let \, $\eta(\cdot) \in C^{\infty}_0({\mathbb R})$ be a function such that \, 
%{\textcolor{blue}{
$0 \leq \eta(x) \leq  1$,
%}} 
\, 
$\displaystyle{|\frac{d}{dx} \eta(x)| \leq 1}$, \, $\forall x \in {\mathbb R}$ \, and 
\begin{equation}
\eta(x) = \left\{ 
           \begin{array}{ll}
             1, & \qquad \qquad  |x| \leq 1; \\ [0.2cm]
             0, & \qquad \qquad  |x| \geq 3.
           \end{array} \right.
\end{equation}
For each $M \in {\mathbb N}$\, and \,$i \in {\mathbb N}$, let 
\begin{equation}
\eta_{M,i}(x) \equiv \eta \left(M^{-1} \cdot {\gamma}_i^{\frac{1}p} \beta_i^{\frac{1}p}\cdot x \right), \qquad x \in {\mathbb R},
\end{equation}
then, by (4.15), (4.24) and (4.25), we see that  
%$\displaystyle{\prod_{i \geq 1} \eta_{M,i} \in l^p_{(\beta_i)}}$, and
\begin{equation}
{\rm{supp}} \, \big[ \prod_{i \geq 1}\eta_{M,i} \big] \subset D_{3M}, \qquad M \in {\mathbb N}.
\end{equation}
For each $f \in C^{\infty}_0({\mathbb R}^n \to {\mathbb R})$, \, $n \in {\mathbb N}$, define
\begin{equation}
f_M(x_1, \dots, x_n, x_{n+1}, \dots)  \equiv f(x_1, \dots, x_n) \cdot \prod_{i \geq 1} \eta_{M,i}(x_i), \quad {\mathbf x}=(x_1, x_2, \dots) \in {\mathbb R}^{\mathbb N},
\end{equation}
and then define a subspace ${\tilde{\cal D}}({\cal E}_{(\alpha)}) \subset L^2(S; \mu)$, that is the {\it{linear span}} of the family of  ${f_M}'s$,  $M \in {\mathbb N}$,  $f \in C^{\infty}_0({\mathbb R}^n \to {\mathbb R})$, $n \in {\mathbb N}$, defined by (4.27), i.e., the {\it{space of finite linear combinations}} of ${f_M}'s$ defined by (4.27):
\begin{equation}
{\tilde{\cal D}}({\cal E}_{(\alpha)}) \equiv 
{\mbox{the linear span of}} \,
\Big\{ f_M, \, \, M \in {\mathbb N}\, \, : \, 
f \in C^{\infty}_0({\mathbb R}^n \to {\mathbb R}), \, n \in {\mathbb N}  \,  \Big\}.
\end{equation}
Note that ${\cal F}C^{\infty}_0$ defined by (2.7) is expressed  in a similar way as (4.27):
\begin{equation}
{\cal F}C^{\infty}_0 = \Big\{ f(x_1, \dots, x_n) \cdot \prod_{i \geq 1}I_{\mathbb R}(x_i) \, \, : \, f \in C^{\infty}_0({\mathbb R}^n \to {\mathbb R}), \, n \in {\mathbb N} \, \Big\},
\end{equation}
where $I_A(x)$, $A \in {\cal B}({\mathbb R})$, denotes the indicator function.
Then, for any $f_M$, defined by (4.27), 
%there exists an $M' \geq 3M$ such that 
\, ${\rm{supp}}[f_M] \subset D_{M'}$ \, (cf. (4.26))  holds 
for all $M' \geq 3M$.   Thus, if we see that 
\begin{equation}
{\tilde{\cal D}}({\cal E}_{(\alpha)}) \subset {\cal D}({\cal E}_{(\alpha)}),
\end{equation}
and for any $\displaystyle{u \equiv f \cdot \prod_{i \geq 1} I_{\mathbb R}(x_i) 
\in {\cal F}C^{\infty}_0 \subset {\cal D}({\cal E}_{(\alpha)})}$ 
(cf. (4.29)),
 there exists a sequence $\{w_m\}_{m \in {\mathbb N}}$, $ w_m \in {\tilde{\cal D}}({\cal E}_{(\alpha)})$, $m \in {\mathbb N}$, such that 
\begin{equation}
\lim_{m \to \infty} w_m = u \equiv f \cdot \prod_{i \geq 1} I_{\mathbb R}(x_i), \qquad {\mbox{in \quad ${\cal D}({\cal E}_{(\alpha)})$}},
\end{equation}
then  (4.13) and (4.14) are verified.

Let us show that (4.30) and (4.31) hold.  For each $M \in {\mathbb N}$.
 Denote 
(cf. (4.15) and (4.26))
\begin{equation}
a_{M,i} \equiv M {\gamma}_i^{- \frac{1}p} \beta_i^{- \frac1p}, \qquad 
b_{M,i} \equiv 3M {\gamma}_i^{- \frac{1}p} \beta_i^{- \frac1p}, \qquad i \in {\mathbb N}.
\end{equation}
Below, we use simplified notations for the indicator functions such that, e.g., 
$I_{[-a_{M,i},a_{M,i}]} (y) = I_{|y| \leq a_{M,i}}$ and 
$I_{(a_{M,i}, b_{M,i}]} (y) = 1 _{a_{M,i} < y \leq b_{M,i}}$ and so on.
From the definition of $\eta$ (cf. (4.24) and (4.25)), 
since $0 \leq \eta_{M,i}(x ) \leq 1$ ($ \forall x \in {\mathbb R}$), and
by the  {\it{mean value theorem}}, 
since for any $x < y$, there exists a $\theta \in (0,1)$ such that 
$ \eta_{M,i}(x ) - \eta_{M,i}(y) = \eta'_{M_i}(x + \theta \cdot (y-x)) \cdot (y-x)$, 
 we have the following bound for $0 < \alpha \leq 1$:
\begin{eqnarray*}
|\eta_{M,i}(x ) - \eta_{M,i}(y) |^2 & =& |\eta_{M,i}(x ) - \eta_{M,i}(y) |^{\alpha + 1} |\eta_{M,i}(x ) - \eta_{M,i}(y) |^{1 - \alpha}
\\
& \leq & 2^{1-\alpha}  |\eta_{M,i}(x ) - \eta_{M,i}(y) |^{\alpha + 1} \\
& \leq & 2^{1-\alpha} |\eta'_{M_i}(x + \theta \cdot (y-x))|^{\alpha +1} \cdot |x-y|^{\alpha +1}.
\end{eqnarray*}
By making use of this bound and 
 noting that $\eta_{M,i}(x) = \eta_{M,i}(y) = 1$ for $|x|, \, |y| \leq a_{M,i}$, from 
the definition (4.25) for $\eta_{M_i}(x)$, we then have  
the following evaluation:
\begin{eqnarray}
\lefteqn{\frac{(\eta_{M,i}(y ) - \eta_{M,i}(x))^2}{|y -x |^{ \alpha +1}} I_{|y| \leq a_{M,i}}  I_{x \ne y} } \nonumber \\
&& = \frac{(\eta_{M,i}(y ) - \eta_{M,i}(x))^2}{|y -x |^{\alpha +1}} \cdot 
\big( I_{|y| \leq a_{M,i}} I_{|x| \leq a_{M,i}} + I_{|y| \leq a_{M,i}} I_{-b_{M,i} \leq x < - a_{M,i}} 
\nonumber \\
&&+ I_{|y| \leq a_{M,i}} I_{a_{M,i} < x \leq b_{M,i}} 
+ I_{|y| \leq a_{M,i}} I_{|x| > b_{M,i}} \big)I_{x \ne y}
\nonumber \\
&& \leq 
2^{1- \alpha} \big(\sup_{t \in {\mathbb R}} | \frac{d}{dt} \eta_{M,i}(t)|^{\alpha +1} \big) 
\cdot 
\big( I_{|y| \leq a_{M,i}} I_{-b_{M,i} \leq x < - a_{M,i}} 
\nonumber \\
&&+ I_{|y| \leq a_{M,i}} I_{a_{M,i} < x \leq b_{M,i}} 
+ I_{|y| \leq a_{M,i}} I_{|x| > b_{M,i}} \big)I_{x \ne y}
\nonumber \\
&& \leq 
2^{1- \alpha} \big(\sup_{t \in {\mathbb R}} | \frac{d}{dt} \eta_{M,i}(t)|^{\alpha +1} \big) 
\cdot 
\big( I_{x < - a_{M,i}} + I_{a_{M,i} < x} + I_{|x| > b_{M,i}} \big) \, I_{x \ne y}
\nonumber \\
&& \leq 
2^{1- \alpha} 3 \big( M^{-1} {\gamma}_i^{\frac{1}p} \beta_i^{\frac1p} \big)^{\alpha +1} I_{|x| > a_{M,i}} I_{x \ne y}.
\end{eqnarray}
Similarly,
\begin{eqnarray}
\lefteqn{\frac{(\eta_{M,i}(y ) - \eta_{M,i}(x))^2}{|y -x |^{ \alpha +1}} 
\big(
I_{-b_{M,i} \leq y < - a_{M,i}} + I_{a_{M,i} < y \leq b_{M,i}} \big) I_{x \ne y}
}  \nonumber \\
&& \leq 
2 \cdot 
2^{1- \alpha} \big(\sup_{t \in {\mathbb R}} | \frac{d}{dt} \eta_{M,i}(t)|^{\alpha +1} \big) 
\cdot 
I_{|y| >a_{M,i}} I_{x \ne y} 
\nonumber \\
 &&\leq \, 
2^{2-\alpha} \big( M^{-1} {\gamma}_i^{\frac{1}p} \beta_i^{\frac1p} \big)^{\alpha + 1}
I_{|y| >a_{M,i}} I_{x \ne y},  
\end{eqnarray}
and 
\begin{eqnarray}
\lefteqn{\frac{(\eta_{M,i}(y ) - \eta_{M,i}(x))^2}{|y -x |^{ \alpha +1}} 
I_{|y| > b_{M,i}}  I_{x \ne y}
} 
\nonumber \\
&& \leq 
2^{1- \alpha} \big(\sup_{t \in {\mathbb R}} | \frac{d}{dt} \eta_{M,i}(t)|^{\alpha +1} \big)  
\cdot 
I_{|y| > b_{M,i}}  I_{x \ne y} 
\nonumber \\
 &&\leq \, 
2^{1- \alpha} \big( M^{-1} {\gamma}_i^{\frac{1}p} \beta_i^{\frac1p} \big)^{\alpha +1} 
I_{|y| >a_{M,i}} I_{x \ne y}. 
\end{eqnarray}
Next, 
for each $f \in C^{\infty}_0({\mathbb R}^n \to {\mathbb R})$, $n \in {\mathbb N}$, define $f_{M,k} \in{\cal F}C^{\infty}_0$, $k \geq n$, (cf. (4.27) and (4.29)) as follows:
\begin{equation}
f_{M,k}(x_1, \dots, x_n, x_{n+1}, \dots)  \equiv f(x_1, \dots, x_n) \cdot \prod_{i =1}^k \eta_{M,i}(x_i) \prod_{j \geq k +1} I_{{\mathbb R}}(x_j), \qquad k \geq n.
\end{equation}
In the sequel 
for $f_{M,k}$ defined through $f \in C^{\infty}_0({\mathbb R}^n \to {\mathbb R})$, 
we assume,  {\it{without mentioning it explicitly}},  that 
 $M$ 
is taken 
sufficiently large in order that the following holds:
$${\rm{supp}}\, [f] \subset [-M \gamma_1^{- \frac1p} \beta_1^{- \frac1p}, +M \gamma_1^{- \frac1p} \beta_1^{- \frac1p}]
\times \cdots \times  [-M \gamma_n^{- \frac1p} \beta_n^{- \frac1p}, +M \gamma_n^{- \frac1p} \beta_n^{- \frac1p}],
$$
%i.e., by (4.25), we assume that  $M$ is taken  sufficiently large in order that
%$
%\prod_{i =1}^n \eta_{M,i} (x_i)  = 1$ holds for \, $(x_1, \dots, x_n) \in {\rm{supp}}\, [f]$.
%Through the above discussion, to define  
%$f_{M,k}$ by (4.36), 
%we 
%always 
%assume that $M< \infty$ is taken 
%sufficiently large in order that the following holds:
%\begin{equation}
% f(x_1, \dots, x_n) \cdot \prod_{i =1}^n \eta_{M,i} (x_i) = f(x_1, \dots, x_n), \qquad \forall (x_1, \dots, x_n) \in {\mathbb R}^n. 
%\end{equation}
%Let    (cf. (4.29), and (4.36)) 
  %$$u \equiv f \cdot \prod_{i \geq n +1} I_{{\mathbb R}}(x_i) \in {\cal F}C^{\infty}_0.$$ 
%By making use of  (4.33), (4.34) and (4.35), we  have the following evaluation for  the quadratic form 
%${\cal E}_{(\alpha)}
%\equiv \sum_{i \in {\mathbb N}} {\cal E}^{(i)}_{(\alpha)}$, (cf. (2.10)),
% defined through ${\cal E}^{(i)}_{(\alpha)}$ given by {\it{(2.13) in Remark 2}}, that is equivalent to the formula (2.9): \, 
% Let $f_{M,k}$, be the function  defined by (4.36) satisfying  (4.37), also 
% let 
%$k \geq n$.  \\
%For \, $1 \leq i \leq n$, since 
%${\displaystyle{\Phi_{\alpha}(f_{M,k},f_{M,k}; y_i,y'_i, {\mathbf x}\setminus x_i) = \Phi_{\alpha}(u,u; y_i,y'_i, {\mathbf x}\setminus x_i)}}$, 
%for all $y, \, y', \, {\mathbf x}$, by (4.37), it holds that,  for $1 \leq i \leq n$
%\begin{equation}
%  {\cal E}^{(i)}_{(\alpha)} (f_{M,k}, f_{M,k}) = {\cal E}^{(i)}_{(\alpha)}(u,u),
%\end{equation}
i.e., by (4.25), we assume that  $M$ is taken  sufficiently large in order that
$
\prod_{i =1}^n \eta_{M,i} (x_i)  = 1$ holds for \, $(x_1, \dots, x_n) \in {\rm{supp}}\, [f]$.
Thus, to define  
$f_{M,k}$ by (4.36), 
assume  that $M< \infty$ is taken 
sufficiently large and the following holds:
\begin{equation*}
 f(x_1, \dots, x_n) \cdot \prod_{i =1}^n \eta_{M,i} (x_i) = f(x_1, \dots, x_n), \qquad \forall (x_1, \dots, x_n) \in {\mathbb R}^n. 
\end{equation*}
Let    (cf. (4.29), and (4.36)) 
  $$u \equiv f \cdot \prod_{i \geq n +1} I_{{\mathbb R}}(x_i) \in {\cal F}C^{\infty}_0,$$ 
%{\textcolor{blue}{
then, by this and (4.36) with  the above equality, which holds for sufficiently large $M < \infty$, for 
 $ n \leq k$  we have
%}}
%{\textcolor{blue}{
\begin{eqnarray}
  \lefteqn{ (f_{M,k}(x_1, \dots, x_{i-1}, y_i, x_{i +1}, \dots, x_n, \dots) 
              - f_{M,k}(x_1, \dots, x_{i-1}, y'_i, x_{i +1}, \dots, x_n, \dots) )^2} 
    \nonumber \\
   &&\leq (u(x_1, \dots, x_{i-1}, y_i, x_{i +1}, \dots, x_n, \dots) 
              - u(x_1, \dots, x_{i-1}, y'_i, x_{i +1}, \dots, x_n, \dots) )^2,
     \nonumber \\
&&{}
\qquad \qquad \forall {\mathbf x} \in {\mathbb R}^{\mathbb N}, \quad 
\forall y_i,\, y'_i \in {\mathbb R}, \quad 1 \leq i \leq n. 
%\qquad \qquad
%\qquad \qquad\qquad \quad (4.37)
\end{eqnarray}
%}}
By making use of  (4.33), (4.34) and (4.35), we  have the following evaluation for  the quadratic form 
${\cal E}_{(\alpha)}
\equiv \sum_{i \in {\mathbb N}} {\cal E}^{(i)}_{(\alpha)}$, (cf. (2.10)),
 defined through ${\cal E}^{(i)}_{(\alpha)}$ given by {\it{(2.13) in Remark 2}}, that is equivalent to the formula (2.9). \, 
%{\textcolor{blue}{
For  $k \geq n$, 
%}}
 let $f_{M,k}$, be the function  defined by (4.36) satisfying  (4.37). Then
for \, $1 \leq i \leq n$, since 
%{\textcolor{blue}{
${\displaystyle{\Phi_{\alpha}(f_{M,k},f_{M,k}; y_i,y'_i, {\mathbf x}\setminus x_i) \leq \Phi_{\alpha}(u,u; y_i,y'_i, {\mathbf x}\setminus x_i)}}$, 
%}}
for all $y, \, y', \, {\mathbf x}$, it holds that,  for $1 \leq i \leq n$,
%{\textcolor{blue}{
\begin{equation}
%{}\qquad \qquad\qquad \quad \qquad
  {\cal E}^{(i)}_{(\alpha)} (f_{M,k}, f_{M,k}) \leq {\cal E}^{(i)}_{(\alpha)}(u,u), 
 %\qquad\qquad\qquad\qquad \quad(4.38)
\end{equation}
%}}
also  for \, $n+1 \leq i \leq k$,
by (4.33), (4.34) and (4.35),  with $C_{\alpha} \equiv 6 \cdot 2^{1- \alpha}$, 
\begin{eqnarray}
\lefteqn{
{\cal E}^{(i)}_{(\alpha)}(f_{M,k},f_{M,k}) 
} \nonumber \\
&& =
\int_S (f(x_1, \dots, x_n))^2 \cdot \prod_{j \geq n+1, j \ne i}^k \big( \eta_{M,j}(x_j) \big)^2 
\prod_{l \geq k+1} \big( I_{\mathbb R}(x_l) \big)^2 
\nonumber \\
&& \times
\Big\{
\int_{{\mathbb R}^2} I_{y_i \ne y'_i} 
\frac{1}{|y_i -y'_i|^{ \alpha + 1}} ( \eta_{M,i}(y_i) - \eta_{M,i}(y'_i))^2 \mu(dy_i\,|\, \sigma_{i^c}) \mu(dy'_i\,|\, \sigma_{i^c}) \Big\} \, \mu(d{\mathbf x})
\nonumber \\
&& \leq 
C_{\alpha} \big( M^{-1} {\gamma}_i^{\frac{1}p} \beta_i^{\frac1p} \big)^{\alpha +1} \int_S (f(x_1, \dots, x_n))^2
\prod_{j \ne i} I_{\mathbb R}(x_j) \big\{ \int_{\mathbb R} I_{|y'_i| > a_{M,i}} \mu(dy'_i \, | \, \sigma_{i^c}) \big\} \mu(d{\mathbf x}) 
\nonumber \\
&& \leq 
C_{\alpha} \big( M^{-1} {\gamma}_i^{\frac{1}p} \beta_i^{\frac1p} \big)^{\alpha +1} \, \|f\|_{L^{\infty}}^2
\int_S \big\{ \int_{\mathbb R} I_{|y'_i| > a_{M,i}} 
\mu(dy'_i \, | \, \sigma_{i^c}) \big\} \mu(d{\mathbf x}) 
\nonumber \\
&& \leq 
C_{\alpha} \big( M^{-1} {\gamma}_i^{\frac{1}p} \beta_i^{\frac1p} \big)^{\alpha +1} \, \|f\|_{L^{\infty}}^2
\, \mu(|x_i| > a_{M,i}),
\end{eqnarray}
and
  for\, $ k+1 \leq i$, 
since 
${\displaystyle{\Phi_{\alpha}(f_{M,k},f_{M,k}; y_i,y'_i, {\mathbf x}\setminus x_i) = 0}}$, 
for all $y, \, y', \, {\mathbf x}$, by (4.36),  for\, $ k+1 \leq i$, 
\begin{equation}
{\cal E}^{(i)}_{(\alpha)} (f_{M,k},f_{M,k}) = 0.
\end{equation}
By combining  (4.38), (4.39) and (4.40), from (4.32), the definition of $a_{M,i}$, 
for $f_{M,k}$ defined through $f \in C^{\infty}_0({\mathbb R}^n \to {\mathbb R})$, $n \in {\mathbb N}$
we have 
\begin{equation*}
{\cal E}^{(i)}_{(\alpha)} (f_{M,k},f_{M,k}) \leq {\cal E}^{(i)}_{(\alpha)} (u,u) + C_{\alpha} \big(
M^{-1} {\gamma}_i^{\frac{1}p} \beta_i^{\frac1p}\big)^{\alpha +1} \, \| f\|_{L^{\infty}}^2 \cdot \mu \left(
|X_i| > M {\gamma}_i^{- \frac{1}{p}} \beta_i^{- \frac1p} \right),
\end{equation*}
\begin{equation}
{\mbox{for any}} 
\quad  k \geq n.
\end{equation}
By (4.41), for $f_{M,k}$ defined through $f \in C^{\infty}_0({\mathbb R}^n \to {\mathbb R})$, $n \in {\mathbb N}$, in order to make the discussions simpler by 
{\it{taking $M$ satisfying both (4.37) and $M \geq M_0$ (cf. (4.3))}}, 
since 
$\mu (
|X_i| > M {\gamma}_i^{- \frac{1}{p}} \beta_i^{- \frac1p} )$
is a decreasing function of $M >0$, 
we see that if the condition (4.3) of Theorem 2-i) is satisfied,  then the following holds: \, 
%{\textcolor{blue}{
There exists a constant $C_{\alpha} < \infty$ that does not depend on $M \geq M_0$ and 
\begin{equation}
\sup_{k \geq 1} {\cal E}_{(\alpha)} (f_{M,k}, f_{M,k}) = \sup_{k \geq 1} \sum_{i \in {\mathbb N}}
{\cal E}^{(i)}_{(\alpha)}(f_{M,k}, f_{M,k}) \leq 
{\cal E}_{(\alpha)}(u,u) + C_{\alpha} M^{-(\alpha +1)} \|f \|_{L^{\infty}}^2,
\end{equation}
%}}
%\begin{equation}
%{\mbox{ for any}}
% \quad  M \geq M_0.
%\end{equation}
for any sufficiently large $M$ by which satisfies  (4.37)  holds with   $M \geq M_0$.\\
Since, by the definitions (4.27) and (4.36), 
$\lim_{k \to \infty} f_{M,k} =f_M$, $\mu-a.e.$ , and by  Lebesgue's bounded convergence theorem \, $\lim_{k \to \infty} \| f_{M,k} - f_M \|_{L^2(S;\mu)} = 0$, (4.42) shows that the sequence $\{f_{M,k} \}_{k \in {\mathbb N}}$, $f_{M,k} \in {\cal F}C^{\infty}_0 \subset {\cal D}({\cal E}_{(\alpha)})$ satisfies the condition for (4.10) of Lemma 1.
Hence, we conclude that
for \, $f_M$ \, defined by (4.27) through $f \in C^{\infty}_0({\mathbb R}^n \to {\mathbb R}), \, n \in {\mathbb N}$, and the subspace  ${\tilde{\cal D}}({\cal E}_{(\alpha)}) \subset L^2(S;\mu)$ defined by (4.28) 
the following is true (cf. (4.10)):
%{\textcolor{blue}{
\begin{equation}
 f_M  \in {\cal D}({\cal E}_{(\alpha)}),  \qquad {\mbox{and}} \qquad {\tilde{\cal D}}({\cal E}_{(\alpha)}) \subset {\cal D}({\cal E}_{(\alpha)}),
\end{equation}
also by (4.42) (cf. (4.10)),
\begin{equation}
{\cal E}_{(\alpha)} (f_{M}, f_{M}) \leq 
\liminf_{k \to \infty} {\cal E}_{(\alpha)} (f_{M,k}, f_{M,k}) 
 \leq 
{\cal E}_{(\alpha)}(u,u) + C_{\alpha} M^{-(\alpha + 1)} \|f \|_{L^{\infty}}^2, 
\end{equation}
%$$
%{\mbox{
for any sufficiently large $M$ by which satisfies  (4.37)  holds with   $M \geq M_0$.\\
%}}. 
%$$
(4.43) shows (4.30).
%}}

%{\textcolor{blue}{
Finally, we apply the same arguments to the case of the sequence 
$\{f_M\}_{M \in {\mathbb N}}$, $f_M \in {\cal D}({\cal E}_{(\alpha)})$. 
From (4.44) we have 
\begin{equation}
\sup_{M \in {\mathbb N}} {\cal E}_{(\alpha)} (f_{M}, f_{M} )  < \infty.
\end{equation}
%}}
Moreover, 
by (4.15), (4.25) with (4.24) (cf. (4.26)), since
\begin{eqnarray*}
\lefteqn{
A_M \equiv 
\{ {\mathbf x} \, : \, \prod_{i \geq 1}  I_{\mathbb R}(x_i) - \prod_{i \geq 1} \eta_{M,i}(x_i) \ne 0 \} = 
\{ {\mathbf x} \, : \,  \exists i \in {\mathbb N}, \, \, \eta_{M,i}(x_i) \ne 1 \} 
} \\
&&
\subset \{ {\mathbf x} \, :\, \exists i \in {\mathbb N}, \, \, |x_i| > M \gamma_i^{- \frac1p} \beta_i^{- \frac1p} \} = 
\bigcup_{i \geq 1} \{ {\mathbf x} \, :\,  |x_i| > M \gamma_i^{- \frac1p} \beta_i^{- \frac1p} \},
\end{eqnarray*}
for $u \equiv f \cdot \prod_{i \geq 1} I_{{\mathbb R}}(x_i) \in {\cal F}C^{\infty}_0$ and $f_M$ defined by (4.27), we see that 
%$$
\begin{eqnarray*}
\lefteqn{
\| u - f_M \|_{L^2(S; \mu)}  
\leq  \| f \|_{L^\infty} \int_S I_{A_M} (\mathbf x) \, \mu(d {\mathbf x})
%\leq \int_S I_{\{ \prod_{i \geq 1} I_{\mathbb R}(x_i) 
%- \prod_{i \geq 1} \eta_{M,i}(x_i) \ne 0 \}} (\mathbf x) \mu(d {\mathbf x})
}\\
&&
 = \| f \|_{L^\infty}  \, \mu \Big( \bigcup_{i \geq 1} \{ {\mathbf x} \, :\, |x_i| > M \gamma_i^{- \frac1p} {\beta}_i^{- \frac1p} \} \Big) = 
\| f \|_{L^\infty}  \,(1 - \mu(D_M)), 
\end{eqnarray*}
%$$
But, under the condition (4.4), it holds that 
$$
\limsup_{M \to \infty} (1 - \mu(D_M)) = 1 - \liminf_{M \to \infty} \mu(D_M)
\leq 1 - \mu (\liminf_{M \to \infty} D_M)$$
$$  = 1- \mu \big( \cup_{M \in {\mathbb N}} D_M \big) =0.
$$
Thus, $\lim_{M \to \infty} \| f_{M} -  u \|_{L^2(S;\mu)} = 0$, 
and hence 
(4.45) shows that the sequence $\{f_{M}  \}_{M \in {\mathbb N}}$, $f_{M} \in {\cal D}({\cal E}_{(\alpha)})$ satisfies the condition for (4.10) of Lemma 1.  Then, by (4.11), the second assertion of Lemma 1, shows that there exists a subsequence $\{f_{M_l}  \}_{l \in {\mathbb N}}$ of $\{f_M  \}_{M \in {\mathbb N}}$ such that the {\it{Ces{\`a}ro mean}}
\begin{equation}
w_m \equiv \frac{1}m \sum_{l=1}^m (f_{M_l} ) \to u\equiv f \cdot \prod_{i \geq 1} I_{{\mathbb R}}(x_i)  \quad {\mbox{in \quad ${\cal D}({\cal E}_{(\alpha)})$ \quad as \quad $n\to \infty$}}.
\end{equation}
This shows (4.31). 
Since, 
as we have mentioned before, 
from the definition (4.28) and by (4.43), $w_m \in {\tilde{\cal D}}({\cal E}_{(\alpha)}) \subset {\cal D}({\cal E}_{(\alpha)})$ holds, (4.43) and (4.46) verify 
(4.13)and (4.14), respectively.
This complete the proof of Theorem 2-i), 
for the case $S = l^p_{({\beta}_i)}$, $ 1 \leq p < \infty$.\\

{\it{\bf{Proof of (4.13) and (4.14) for $S = l^{\infty}_{(\beta_i)}$.}}} \quad
The same arguments, through which (4.13) and (4.14) are verified for the case where  $S = l^{p}_{(\beta_i)}, \, 1 \leq p <\infty$, can be applied to the case where 
 $S = l^{\infty}_{(\beta_i)}$ as defined by (2.2).  Namely, for  $S = l^{\infty}_{(\beta_i)}$, a definition of the 
compact set $D_M$ corresponding to (4.15) is 
\begin{equation}
D_M \equiv \left\{ {\mathbf x} \in l^{\infty}_{(\beta_i)} \, : \, \, {\beta}_i | x_i|  \leq 
M \cdot {\gamma}_i^{-1},  \,  i \in {\mathbb N} 
 \right\}, \qquad M \in {\mathbb N},
\end{equation}
where $\{\gamma_i \}_{i \in {\mathbb N}}$ is any sequence such that 
$0 < \gamma_1 \leq \gamma_2 \leq \cdots \to \infty$, and the evaluation corresponding to (4.41) is 
\begin{equation*}
{\cal E}^{(i)}_{(\alpha)} (f_{M,k},f_{M,k} )\leq {\cal E}^{(i)}_{(\alpha)} (u,u) + C_{\alpha} \big(
M^{-1} \gamma_i \beta_i \big)^{\alpha +1} \, \| f\|_{L^{\infty}} \cdot \mu \left(
|X_i| > M \gamma_i^{-1} \beta_i^{- 1} \right),
\end{equation*}
\begin{equation}
{\mbox{for any}} \quad 
k \geq n.
\end{equation}
Passing through the same arguments by which (4.46) is derived (cf. (4.45)), under the conditions (4.5) with (4.6), we can prove (4.13) and (4.14) for $S=l^{\infty}_{(\beta_i)}$, and Theorem 2-ii) follows.\\

{\it{\bf{Proof of (4.13) and (4.14) for $S = {\mathbb R}^{\mathbb N}.$}}} \quad
Similar to the above, 
the same arguments, through which (4.13) and (4.14) are verified for the case where  $S = l^{p}_{(\beta_i)}, \, 1 \leq p <\infty$, can be applied to the case where 
 $S = {\mathbb R}^{\mathbb N}$ defined by (2.3).  Namely, for  $S = {\mathbb R}^{\mathbb N}$, a definition of the 
compact set $D_M$ corresponding to (4.15) is 
\begin{equation}
D_M \equiv \left\{ {\mathbf x} \in {\mathbb R}^{\mathbb N} \, : \, \,  | x_i|  \leq 
M \cdot {\gamma}_i,  \quad  i \in {\mathbb N} 
 \right\}, \qquad M \in {\mathbb N},
\end{equation}
where $\{\gamma_i \}_{i \in {\mathbb N}}$ is any sequence such that 
$\gamma_i > 0$, $\forall i \in {\mathbb N}$, and the evaluation corresponding to (4.41) is 
\begin{equation}
{\cal E}^{(i)}_{(\alpha)} (f_{M_k},f_{M_k}  )\leq {\cal E}^{(i)}_{(\alpha)} (u,u) + C_{\alpha} \big(
M^{-1} \gamma_i^{-1} \big)^{\alpha +1} \| f\|_{L^{\infty}}^2 \cdot \mu \left(
|X_i| > M \gamma_i \right),
\,  k \geq n.
\end{equation}
For the probability distribution of the real valued random variable $X_i$, since
\begin{equation}
\mu(|X_i| > M \gamma_i) \leq 1, \qquad \forall i \in {\mathbb N}, \quad \forall M \geq 0, 
\quad \forall \gamma_i \geq 0,
\end{equation}
by taking $\gamma_i \equiv i$,  $i \in {\mathbb N}$, in (4.50), then from (4.50) and (4.51) we see that, for any $f \in C^{\infty}_0({\mathbb R}^n \to {\mathbb R})$, $n \in {\mathbb N}$,
$$
{\cal E}_{(\alpha)}(f_M, f_M ) \leq \, {\cal E}_{(\alpha)}(u,u) + C_{\alpha}  M^{-(\alpha + 1)} \sum_{i =1}^{\infty}  i^{-(\alpha +1)}  < \infty, \qquad \forall M \in {\mathbb N}.
$$
Thus, by passing through the same arguments by which (4.46) is derived (cf. (4.45)), we  prove (4.13) and (4.14) for $S={\mathbb R}^{\mathbb N}$, and Theorem 2-iii) follows.

%{\textcolor{blue}{
Next, let us prove iv).  The point iv) is assured when we  show that $1 \in {\cal D}({\cal E}_{(\alpha)})$. The proof is common to
all the state spaces we are considering, namely $l^p_{(\beta_i)},\, l^{\infty}_{(\beta_i)},\, {\mathbb R}^{\mathbb N}$.
%}}

%{\textcolor{blue}{
Take $\eta \in C_0^{\infty}({\mathbb R} \to {\mathbb R})$ as (4.24), namely let $\eta$ be  such that $\eta (x) \geq 0$, $|\frac{d}{dx} \eta (x)| \leq 1$ for $x \in {\mathbb R}$, and $\eta (x) =1$ for $|x| <1$; $\eta (x) = 0$ for $|x| >3$, and define $u_M(x_1, x_2, \dots) \equiv \eta(x_1\cdot M^{-1}) \prod_{i \geq 2}I_{\mathbb R}(x_i) \in {\cal F}C^{\infty}_0 \subset {\cal D}({\cal E}_{(\alpha)})$ for each $M \in {\mathbb N}$. 
Then, for $0 < \alpha \leq 1$, by (2.8) (cf. (4.44)-(4.35)), 
 for any sufficiently large $M$, since 
$0 \leq \Phi_{\alpha}(u_M.u_M;y_1,y'_1, {\mathbf x} \setminus x_1) \leq 1$, for all $y_1, y'_1 \in {\mathbb R}$ and all ${\mathbf x} \in S$, and for $i \ne 1$, $\Phi_{\alpha}(u_M.u_M;y_i,y'_i, {\mathbf x} \setminus x_i) = 0$ for all $y_i, y'_i \in {\mathbb R}$ and all ${\mathbf x} \in S$,  we see  that 
(using (2.9) and (2.10) and estimates similar to those in (4.33)-(4.35))  $\sup_{M \in {\mathbb N}} {\cal E}_{(\alpha)}(u_M,u_M) < \infty$. 
Since, $\lim_{M \to \infty} u_M({\mathbf x}) = 1 = \prod_{i \geq 1}I_{\mathbb R}(x_i)$ 
pointwise,
and also in $L^2(S; \mu)$,  it follows from the first part of  Lemma 1 that we have $1 \in {\cal D}({\cal E}_{(\alpha)})$.
%}}

This completes the proof of Theorem 2.
%\qquad \qquad \qquad \qquad \qquad \qquad \qquad \qquad 
%\qquad \qquad  \qquad 
\qquad \quad
\qquad  \qquad \qquad  \qquad \qquad \quad 
\bsquare

\bigskip
\bigskip

In the case where $1 < \alpha <2$, for the quasi-regularity of the  Dirichlet form $({\cal E}_{(\alpha)}, {\cal D}({\cal E}_{(\alpha)}))$ considered in  Theorem 1,  the following Theorem 3 holds:
%{\it{The conditions given in Theorem 3 are relatively not easier than the ones in Theorem 2 for $0 < \alpha \leq 1$ (cf., the assumption (2.11) with (2.12)).}}
\begin{theorem} {}\qquad
Let \,  
%{\textcolor{blue}{
$1< \alpha < 2$. 
%}}
 Suppose that the assumption (2.11) with (2.12) hold.
Let  $({\cal E}_{(\alpha)}, {\cal D}({\cal E}_{(\alpha)}))$ be the closed  Markovian symmetric form defined at the beginning of this section through Theorem 1.
Then the following statements hold:\\
i)  \quad In the case where $S= l^p_{(\beta_i)}$, for some $ 1 \leq p < \infty$, as defined by (2.1),  suppose that  
(4.4) holds and that 
there exists a
positive $l^1$ sequence $\{\gamma^{-1}_i \}_{i \in {\mathbb N}}$  and 
 an  $ M_0 \in (0,\infty)$,  such that   
\begin{equation}
\sum_{i =1}^{\infty} \big(\beta^{\frac{1}p}_i {\gamma}_i^{\frac{1}p} \big)^{ \alpha +1} \cdot
\mu \Big(\beta_i^{\frac{1}p} |X_i| > M_0 \cdot {\gamma}_i^{-\frac{1}{p}} \Big) < \infty,
\end{equation}
holds.
Also suppose that 
%for each $M \geq M_0$, suppose that 
%\begin{equation*}
%\lim_{M \to \infty}  M^{- \alpha} 
%\sum_{i=1}^{\infty} L_{M,i} \cdot
%\big(\beta^{\frac1p}_i {\gamma}_i^{\frac{1}p} \big)^{ \alpha}\cdot
%\mu \Big( \beta_i^{\frac{1}p} |X_i| > M\cdot {\gamma}_i^{-\frac{1}{p}} \Big) < \infty,
%\end{equation*}
%and 
\begin{equation}
\sup_{M \geq M_0}
M^{- \alpha} 
\sum_{i=1}^{\infty} L_{M,i} \cdot
\big(\beta^{\frac1p}_i {\gamma}_i^{\frac{1}p} \big)^{ \alpha}\cdot
\mu \Big( \beta_i^{\frac{1}p} |X_i| > M\cdot {\gamma}_i^{-\frac{1}{p}} \Big) < \infty,
\end{equation}
holds, 
where for each $M < \infty$ and  
 $i \in {\mathbb N}$,
 $L_{M,i}$ is the  bound, assumed in assumption (2.12),  for the conditional probability density $\rho$ 
corresponding to the 
compact set  
$$K(M,i) \equiv \big[ -3M\cdot \beta_i^{-\frac1p} \, {\gamma}_i^{-\frac{1}p}, \, 3M\cdot
\beta_i^{-\frac1p} \,  {\gamma}_i^{-\frac{1}p} \big] \subset {\mathbb R}.$$
Then $({\cal E}_{(\alpha)}, {\cal D}({\cal E}_{(\alpha)}))$ is a strictly quasi-regular Dirichlet form.
%taken for $K$ 
 \\
ii) \quad 
 In the case where $S= l^{\infty}_{(\beta_i)}$ as defined by (2.2), suppose that 
(4.6) holds, and that 
there exists a 
sequence $\{ \gamma_i \}_{i \in {\mathbb N}}$ such that $0 < \gamma_1  \leq \gamma_2 \leq \cdots \to \infty$ and 
 an  $0 < M_0 < \infty$, 
% and  (4.6),  and in addition, 
and 
\begin{equation}
 \sum_{i=1}^{\infty} 
\big(\beta_i {\gamma}_i \big)^{ \alpha +1} \cdot
\mu \Big( \beta_i |X_i| > M_0 \cdot \gamma_i^{-1} \Big) < \infty,
\end{equation}
holds.
Also suppose that 
\begin{equation}
\sup_{M \geq M_0}  M^{- \alpha} \sum_{i=1}^{\infty} L_{M,i} \cdot
\big(\beta_i {\gamma}_i \big)^{ \alpha} \cdot
\mu \Big( \beta_i |X_i| > M \cdot  \gamma_i^{-1} \Big) < \infty,
\end{equation}
holds, 
where for each $M < \infty$ and  
 $i \in {\mathbb N}$,
 $L_{M,i}$ is the  bound, assumed in assumption (2.12),  for the conditional probability density $\rho$ 
corresponding to the 
compact set  
$$K(M,i) \equiv \big[ -3M\cdot \beta_i^{-1} \, {\gamma}_i^{-1},  \, 3M\cdot
\beta_i^{- 1} \,  {\gamma}_i^{- 1} \big] \subset {\mathbb R}.$$
 Then $({\cal E}_{(\alpha)}, {\cal D}({\cal E}_{(\alpha)}))$ is a strictly quasi-regular Dirichlet form.
 %where for each  $i \in {\mathbb N}$,
 %$L_{i}$ is the uniform bound for the conditional probability density $\rho$ 
%on the   given compact set  
%$$K_{M,i} \equiv \big[ -6M\cdot \beta_i^{-1} \, \gamma_i^{-1}, 6M\cdot
%\beta_i^{-1} \,  \gamma_i^{-1} \big] \subset {\mathbb R}$$
%taken for $K$ 
%given 
%in the assumption (2.12).
\\
iii) \quad 
 In the case where $S= {\mathbb R}^{\mathbb N}$ as defined by (2.3),  suppose that  there exists a sequence 
$\{ \gamma_i \}_{i \in {\mathbb N}}$ such that $0 < \gamma_i$, $\forall i \in {\mathbb N}$, and 
%an  $0 < M_0 < \infty$,  and the limit 
\begin{equation}
\sup_{M \geq M_0}  M^{-\alpha} \sum_{i = 1}^{\infty} L_{M,i} \cdot \gamma_i^{-{\alpha}} \cdot 
\mu \Big( |X_i| > M \cdot \gamma_i \Big) < \infty,
\end{equation}
holds,  
where for each  
 $i \in {\mathbb N}$,
 $L_{M,i}$ is the  bound, assumed in assumption (2.12),  for the conditional probability density $\rho$ 
corresponding to the 
compact set  
$$K(M,i) \equiv \big[ -6M\cdot \gamma_i, \, \, 6M\cdot
\gamma_i \big] \subset {\mathbb R}.$$
Then 
 $({\cal E}_{(\alpha)}, {\cal D}({\cal E}_{(\alpha)}))$ is a strictly quasi-regular Dirichlet form. 
 %where for each  $i \in {\mathbb N}$,
 %$L_{i}$ is the uniform bound for the conditional probability density $\rho$ 
%given 
%in the assumption (2.12). 
%where for each  $M \in {\mathbb N}$ and $i \in {\mathbb N}$,
% $L_{M,i}$ is the bound for the conditional probability density $\rho$ on the given compact set  
%$$K_{M,i} \equiv \big[ -6M\cdot \gamma_i, \, 6M\cdot
%\gamma_i \big] \subset {\mathbb R}$$
%taken for $K$ in the assumption (2.12).
%\\
%\qquad \qquad \qquad \qquad \qquad \qquad \qquad \qquad \qquad \qquad  \qquad \qquad
%\qquad  \qquad \qquad  \qquad \qquad \quad 
%\bsquare
\end{theorem}

{\bf{Proof of Theorem 3.}} \quad We use the same notations and methods as in the proof of Theorem 2.  To prove  the quasi-regularity we have  to verify (4.12), (4.13) and (4.14) for the case where $1 < \alpha < 2$.  
%{\textcolor{blue}{
To prove 
the strictly part of 
 the quasi-regular Dirichlet form  we have then  to verify 
in addition 
that $1\in {\cal D}({\cal E}_{(\alpha)})$.
%}}

{\bf{Proof of (4.12), (4.13) and (4.14) for $S = l^p{(\beta_i)}$, $1 \leq p < \infty$.}}\\
 Since the compactness of $D_M$ is independent of $\alpha$, passing through the same argument  between (4.15) and (4.31), we can verify (4.12), (4.13) and (4.14) by showing that  (4.30) and (4.31) hold for $1 < \alpha < 2$.

For each $M \in {\mathbb N}$ and $i \in {\mathbb N}$,  
let $a_{M,i}$ and $b_{M,i}$ be the numbers defined by (4.32):
$$a_{M,i} \equiv {\gamma}_i^{-\frac{1}p}, \qquad b_{M,i} \equiv 3{\gamma}_i^{-\frac{1}p}.$$
we note that, by (4.24) and (4.25), the following holds:
\begin{equation}
 \eta_{M,i}(x) - \eta_{M,i}(y) = 0, \qquad 
\forall (x,y) \in  {\cal N}_{M,i},
\end{equation}
where 
$${\cal N}_{N,i} \equiv \{(x,y) \, : \, |x| \leq a_{M,i}, |y| \leq a_{M,i}\} \cup 
\{(x,y) \, : \, |x| > b_{M,i}, |y| > b_{M,i}\}.
$$ 
Since
\begin{eqnarray}
\lefteqn{
  {\mathbb R}^2 \setminus {\cal N}_{M,i}
     }
    \nonumber \\
 &&\subset \Big( \{ (x,y) \, : \, |y| > a_{M,i}, |x| \leq b_{M,i} \}
     \cup \{ (x,y) \, : \, |x| > a_{M,i}, |y| \leq b_{M,i} \} 
% \nonumber \\
%&&{}  \qquad \cup \{ (x,y) \, : \, x > a_{M,i}, |y| \leq b_{M,i} \} 
%\cup \{ (x,y) \, : \, x < -a_{M,i}, |x| \leq b_{M,i} \}
     \Big),
\nonumber
\end{eqnarray}
for each $M \in {\mathbb N}$ and $i \in {\mathbb N}$, (cf. (4.25)), by setting
\begin{equation}
K_{M,i} \equiv \big[ -6M\, {\gamma}_i^{-\frac{1}p} \, \beta_i^{-\frac1p}, 6M\, {\gamma}_i^{-\frac{1}p} \, \beta_i^{-\frac1p} \big] \subset {\mathbb R},
\end{equation}
from (4.57) we have the following bound which is an analogue of (4.33), (4.34) and (4.35):
\begin{eqnarray}
\lefteqn{
\frac{(\eta_{M,i}(y ) - \eta_{M,i}(x))^2}{|y -x |^{2}} I_{|y|   I_{x \ne y} }
\frac{1}{| y - x|^{ \alpha -1}} \cdot \big( I_{K_{M,i}}(x-y) + I_{K_{M,i}^c}(x-y) \big)
} \nonumber \\
&& \leq 
 \big( M^{-1} {\gamma}_i^{\frac{1}p} \beta_i^{\frac1p} \big)^2  I_{x \ne y}
\frac{1}{| y - x|^{ \alpha -1}} \cdot \big( I_{K_{M,i}}(x-y) + I_{K_{M,i}^c}(x-y) \big) 
\nonumber \\
&& \times (I_{|y| > a_{M,i}} I_{|x| \leq b_{M,i}} \, + I_{|x| > a_{M,i}} I_{|y| \leq b_{M,i} }),
\end{eqnarray}
where and in the sequel, to simplify the notations,  we  denote,  e.g., 
$I_{(- \infty, -a_{M,i}) \cup (a_{M,i}, \infty)}(x)$ by $I_{|x| > a_{M,i}}$, and
$I_{[- b_{M,i},\, b_{M,i}]} (y)$ by $I_{|y| \leq b_{M,i}}$, respectively.
%The formulas corresponding to (4.34) and (4.35), respectively,  are 
%\begin{eqnarray}
%\lefteqn{\frac{(\eta_{M,i}(y ) - \eta_{M,i}(x))^2}{|y -x |^{2}} 
%\big(
%I_{-b_{M,i} \leq y < - a_{M,i}} + I_{a_{M,i} < y \leq b_{M,i}} \big) I_{x \ne y}
%} \nonumber \\
%&& \times \frac{1}{| y - x|^{ \alpha -1}} \cdot \big( I_{K_{M,i}}(x-y) + I_{K_{M,i}^c}(x-y) \big)
% \\
%&& \leq 
%2 \big( M^{-1} {\gamma}_i^{\frac{1}p} \beta_i^{\frac1p} \big)^2 
%I_{|y| >a_{M,i}} I_{x \ne y}
%\frac{1}{| y - x|^{ \alpha -1}} \cdot \big( I_{K_{M,i}}(x-y) + I_{K_{M,i}^c}(x-y) \big),
%\nonumber 
%\end{eqnarray}
%and
%\begin{eqnarray}
%\lefteqn{\frac{(\eta_{M,i}(y ) - \eta_{M,i}(x))^2}{|y -x |^{2}} 
%I_{|y| > b_{M,i}}  I_{x \ne y}
%\frac{1}{| y - x|^{ \alpha -1}} \cdot \big( I_{K_{M,i}}(x-y) + I_{K_{M,i}^c}(x-y) \big)
%} \nonumber \\
%&& \leq 
% \big( M^{-1} {\gamma}_i^{\frac{1}p} \beta_i^{\frac1p} \big)^2 
%I_{|y| >a_{M,i}} I_{x \ne y}
%\frac{1}{| y - x|^{ \alpha -1}} \cdot \big( I_{K_{M,i}}(x-y) + I_{K_{M,i}^c}(x-y) \big),
%\nonumber \\
%&&{}
%\end{eqnarray}
%respectively.
Next, 
recall the assumption  (2.11) with (2.12) for a bound of the  conditional probability density $\rho$ corresponding to a given 
compact set.
In the present situation, there exists an $L_{M,i} < \infty$ and for $a.e. y \in {\mathbb R}$ with respect to the Lebesgue measure, such that 
\begin{equation}
\rho(y| {\sigma}^c_i) \cdot I_{[- b_{M,i},\, b_{M,i}]} (y) \leq L_{M,i}, \qquad  \mu-a.e.
\end{equation}
Then, by applying Young's inequality we have that 
(cf. also  (3.16))
 \begin{eqnarray}
\lefteqn{
%{\textcolor{red}{
\sup_{x \in {\mathbb R}} | \int_{\mathbb R} \frac{I_{K_{M,i}}(x-y)}{|x-y|^{\alpha -1}}  \rho(y |{\sigma}^c_i) 
\cdot I_{[- b_{M,i},\, b_{M,i}]} (y) \,
 dy |
%}}
}  \\
% \|\frac{I_{K_{M,i}}}{|\cdot|^{ \alpha -1}} \ast \rho(\cdot \, \big| \, \sigma_{i^c})\|_{L^{\infty}({\mathbb R})} 
&&\leq L_{M,i} \, \| \frac{I_{K_{M,i}}}{|\cdot|^{ \alpha -1}}\|_{L^1({\mathbb R})} 
= \frac{2L_{M,i}}{2- \alpha} \big( 6 M {\gamma}_i^{- \frac{1}p} \beta_i^{- \frac1p} \big)^{2-  \alpha}, 
\nonumber
\quad  \mu-a.e..
\end{eqnarray}
Also by the definition of $K_{M,i}$, we have that
\begin{equation}
\frac{1}{|y_i -y'_i|^{ \alpha -1}} I_{K_{M,i}^c}(y_i -y'_i)  \leq \big( 6^{-1} M^{-1} {\gamma}_i^{\frac{1}p} \beta_i^{ \frac1p} \big)^{ \alpha -1}, \qquad \forall y_i, \, y'_i \in {\mathbb R},
\end{equation}
where  we used (4.58) explicitly.
 For the given $f \in C_0^{\infty}({\mathbb R}^n \to {\mathbb R})$, let us take $M <\infty$ 
sufficiently large in order that (4.37) holds.
Then, by 
%{\textcolor{red}{
(2.13) in Remark 2,
%}}
(4.61) and (4.62),  
from (4.58), (4.59) and (4.60), 
the estimates corresponding to  (4.38), (4.39) and (4.40), for the present 
value of 
$\alpha$,  
are given as follows: \\
for \, $1 \leq i \leq n$, 
\begin{equation}
%{\mbox{for \, $1 \leq i \leq n$}}, \quad 
{\cal E}^{(i)}_{(\alpha)} (f_{M,k}, f_{M,k} ) \leq {\cal E}^{(i)}_{(\alpha)}(u,u), 
\end{equation}
and for \, $n+1 \leq i \leq k$ (below, for the first inequality we use  (4.59), 
for the second and third  inequality we apply (4.61) and (4.62), also, to simplify the notations,  we denote, e.g., 
$I_{[- b_{M,i},\, b_{M,i}]} (y_i)$ by $I_{|y_i| \leq b_{M,i}}$),
\begin{eqnarray}
\lefteqn{
{\cal E}^{(i)}_{(\alpha)}(f_{M,k},f_{M,k}) 
} \nonumber \\
&& =
\int_S (f(x_1, \dots, x_n))^2 \cdot \prod_{j \geq n+1, j \ne i}^k \big( \eta_{M,j}(x_j) \big)^2 
\prod_{l \geq k+1} \big( I_{\mathbb R}(x_l) \big)^2 
\nonumber \\
&& \times
\Big\{
\int_{{\mathbb R}^2} I_{y_i \ne y'_i} 
\frac{1}{|y_i -y'_i|^{2}} ( \eta_{M,i}(y_i) - \eta_{M,i}(y'_i))^2 
\frac{1}{| y_i - y'_i|^{ \alpha -1}} \cdot \big( I_{K_{M,i}}(y_i-y'_i) + I_{K_{M,i}^c}(y_i-y'_i) \big)
\nonumber \\
&& \times
\mu(dy_i\,|\, \sigma_{i^c}) \mu(dy'_i\,|\, \sigma_{i^c}) \Big\} \, \mu(d{\mathbf x})
\nonumber \\
&& \leq 
2 \big( M^{-1} {\gamma}_i^{\frac{1}p} \beta_i^{\frac1p} \big)^2 \int_S (f(x_1, \dots, x_n))^2
\prod_{j \ne i} I_{\mathbb R}(x_j) 
\nonumber \\
&& \times 
\Big\{ \int_{{\mathbb R}^2} I_{|y'_i| > a_{M,i}} 
\frac{1}{| y_i - y'_i|^{ \alpha -1}} \cdot \big( I_{K_{M,i}}(y_i-y'_i) + I_{K_{M,i}^c}(y_i-y'_i) \big)
I_{|y_i| \leq b_{M,i}}
\mu(dy_i \, | \, \sigma_{i^c})
\mu(dy'_i \, | \, \sigma_{i^c}) \Big\} \mu(d{\mathbf x}) 
\nonumber \\
&& \leq 
2 \big( M^{-1} {\gamma}_i^{\frac{1}p} \beta_i^{\frac1p} \big)^2 \, \|f\|_{L^{\infty}}^2
\Big(  \|\frac{I_{K_{M,i}}}{|\cdot|^{ \alpha -1}} \ast 
(\rho(\cdot \, \big| \, \sigma_{i^c}) \cdot I_{|y_i| \leq b_{M,i}})
\|_{L^{\infty}({\mathbb R})} + \big( 6^{-1} M^{-1} {\gamma}_i^{\frac{1}p} \beta_i^{ \frac1p} \big)^{ \alpha -1} \Big)
\nonumber \\
&& \times
\int_S \big\{ \int_{\mathbb R} I_{|y'_i| > a_{M,i}} 
\mu(dy'_i \, | \, \sigma_{i^c}) \big\} \mu(d{\mathbf x}) 
\nonumber \\
&& \leq 
2 \big( M^{-1} {\gamma}_i^{\frac{1}p} \beta_i^{\frac1p} \big)^2 \, 
\Big\{ \frac{2L_{M,i}}{2- \alpha} \big( 6 M {\gamma}_i^{- \frac{1}p} \beta_i^{- \frac1p} \big)^{2-  \alpha} + \big( 6^{-1} M^{-1} {\gamma}_i^{\frac{1}p} \beta_i^{ \frac1p} \big)^{ \alpha -1} \Big\}
\|f\|_{L^{\infty}}^2
\, \mu(|x_i| > a_{M,i}), 
\nonumber \\
&& =2 
\Big\{
\frac{2 \cdot 6^{2- \alpha}L_{M,i} }{2 - \alpha}  \big(M^{-1} {\gamma}_i^{\frac{1}p} \beta_i^{ \frac1p} \big)^{ \alpha} + 6^{1- \alpha} \big( M^{-1} {\gamma}_i^{\frac{1}p} \beta_i^{ \frac1p} \big)^{ \alpha + 1} \Big\} \,
 \|f\|_{L^{\infty}}^2
\, \mu(|x_i| > a_{M,i}), 
%\nonumber 
%\\
%&&{}
%\qquad 
%\qquad \qquad \qquad \qquad 
%{\mbox{for \, $n+1 \leq i \leq k$}};
\end{eqnarray}
also 
for\, $ k+1 \leq i$, 
\begin{equation}
{\cal E}^{(i)}_{(\alpha)} (f_{M,k},f_{M,k} ) = 0, 
%\qquad {\mbox{for\, $ k+1 \leq i$}},
\end{equation}
respectively.
By combining  (4.63), (4.64) and (4.65), from (4.32), the definition of $a_{M,i}$, 
for $f_{M,k}$ defined through $f \in C^{\infty}_0({\mathbb R}^n \to {\mathbb R})$, $n \in {\mathbb N}$, 
we have 
\begin{eqnarray}
\lefteqn{
{\cal E}^{(i)}_{(\alpha)} (f_{M,k},f_{M,k} )
} \nonumber \\
&& \leq {\cal E}^{(i)}_{(\alpha)} (u,u) + 2
\Big\{
\frac{2 \cdot 6^{2- \alpha}L_{M,i}}{1 - \alpha}  \big(M^{-1} {\gamma}_i^{\frac{1}p} \beta_i^{ \frac1p} \big)^{\alpha} + 6^{1- \alpha} \big( M^{-1} {\gamma}_i^{\frac{1}p} \beta_i^{ \frac1p} \big)^{ \alpha + 1} \Big\}
\nonumber \\
&& \times \| f\|_{L^{\infty}}^2 \cdot \mu \left(
|X_i| > M {\gamma}_i^{- \frac{1}{p}} \beta_i^{- \frac1p} \right),
\quad \forall k \geq n.
\end{eqnarray}
By (4.66), for $f_{M,k}$ defined through $f \in C^{\infty}_0({\mathbb R}^n \to {\mathbb R})$, $n \in {\mathbb N}$, in order to make the discussions simple by taking $M$ satisfying both (4.37) and $M \geq M_0$ (cf. (4.52)),
since 
$\mu (|X_i| > M {\gamma}_i^{- \frac{1}{p}} \beta_i^{- \frac1p})$ is a decreasing function of $M >0$, 
we see that if the conditions (4.52)  and (4.53) of Theorem 3-i) are satisfied,  then the following 
 analogue of (4.42) 
 for the present value of  $\alpha$ holds:
\begin{equation}
\sup_{k \geq 1} {\cal E}_{(\alpha)} (f_{M,k}, f_{M,k} ) = \sup_{k \geq 1} \sum_{i \in {\mathbb N}}
{\cal E}^{(i)}_{(\alpha)}(f_{M,k}, f_{M,k}) < \infty.
\end{equation}
By repeating then the same discussions between (4.42) and (4.46) that have been performed for $0 <\alpha \leq 1$,    the proof of Theorem 3-i) is completed 
except for the "strictly" property, that we shall discuss below.

{\bf{Proof of (4.12), (4.13) and (4.14)  for $S = l^{\infty}{(\beta_i)}$, and for $S = {\mathbb R}^{\mathbb N}$.}} \\
The proof for these cases can be done 
in exactly the same way 
%by a completely same manner as 
for the case $S=l^p_{(\beta_i)}$, i.e., by the discussions between (4.57) and (4.67). 
Namely, by changing  $K_{M,i}$ of (4.58) to (cf. (4.47), (4.48))
\begin{equation}
K_{M,i} \equiv \big[-6M\cdot\beta_i^{-1} \gamma_i^{-1}, \, -6M\cdot\beta_i^{-1} \gamma_i^{-1} \big],
\end{equation}
and (cf. (4.49), (4.50))
\begin{equation}
K_{M,i} \equiv \big[-6M\cdot\gamma_i, \, -6M\cdot \gamma_i \big],
\end{equation}
respectively, similarly as in the proof of the corresponding points in Theorem 1,
then through the same discussions as above, we deduce Theorem 3-ii), 3-iii), again except 
for
the "strictly" property.

%{\textcolor{blue}{
Next, let us prove that the quasi-regular Dirichlet forms $({\cal E}_{(\alpha)}, {\cal D}({\cal E}_{(\alpha)}))$  are actually  {\it{strictly}} quasi-regular Dirichlet forms. 
As remarked in the proof of Theorem 2, 
for this it suffices to show that $1 \in {\cal D}({\cal E}_{(\alpha)})$ (see Proposition V-2.15 of [M,R 92]). 
The proof is common for 
every state space $l^p_{(\beta_i)},\, l^{\infty}_{(\beta_i)},\, {\mathbb R}^{\mathbb N}$, and for all $1 < \alpha <2$.
%}}
% {\textcolor{blue}{
Take $\eta \in C_0^{\infty}({\mathbb R} \to {\mathbb R})$ as (4.24), namely let $\eta$ be  such that $\eta (x) \geq 0$, $|\frac{d}{dx} \eta (x)| \leq 1$ for $x \in {\mathbb R}$, and $\eta (x) =1$ for $|x| <1$; $\eta (x) = 0$ for $|x| >3$, and define $u_M(x_1, x_2, \dots) \equiv \eta(x_1\cdot M^{-1}) \prod_{i \geq 2}I_{\mathbb R}(x_i) \in {\cal F}C^{\infty}_0 \subset {\cal D}({\cal E}_{(\alpha)})$ for each $M \in {\mathbb N}$. 
Then, we can set $n=0$ and $k=1$ in (4.64) and (4.65) (although the original argument of (4.64) and (4.65) 
was  performed for $n \in {\mathbb N}$, this extension to $n=0$ is admissible), and we have
\begin{eqnarray*}
{\cal E}^{(1)}_{(\alpha)}(u_M,u_M) 
&\leq& 
\Big\{
\frac{6^{3- \alpha}}{2 - \alpha} L_{1} \big(M^{-1} {\gamma}_1^{\frac{1}p} \beta_1^{ \frac1p} \big)^{ \alpha} + 6^{2- \alpha} \big( M^{-1} {\gamma}_1^{\frac{1}p} \beta_1^{ \frac1p} \big)^{ \alpha + 1} \Big\} \\
&\times&
 \mu(|x_1| > M {\gamma}_1^{-\frac{1}p} \beta_1^{ -\frac1p}), 
\end{eqnarray*}
$$
{\cal E}^{(i)}_{(\alpha)}(u_M,u_M)  = 0, \qquad \forall i \geq 2.$$
By this evaluation, under the assumptions (4.53), (4.55) or (4.56), it holds that 
 $$\sup_{M \in {\mathbb N}} {\cal E}_{(\alpha)}(u_M,u_M) < \infty.$$
Since, $\lim_{M \to \infty} u_M({\mathbf x}) = 1 = \prod_{i \geq 1}I_{\mathbb R}(x_i)$ pointwise,
and also in $L^2(S; \mu)$,  from the first part of  Lemma 1 we get $1 \in {\cal D}({\cal E}_{(\alpha)})$.
%}}

This completes the proof of Theorem 3.
\quad \qquad \qquad \qquad \qquad \qquad \qquad \qquad %\qquad  \qquad \qquad
%\qquad  \qquad \qquad  \qquad \qquad \quad 
\bsquare

\section{Associated Markov processes and a standard procedure of application 
of stochastic quantizations on ${\cal S}'$}

Let $({\cal E}_{(\alpha)}, {\cal D}({\cal E}_{(\alpha)}))$, $0 < \alpha <2$,  be the  
family of 
strictly quasi-regular Dirichlet forms 
on $L^2(S; \mu)$ with a state space $S$ (cf.  (2.1), (2.2), (2.3))
defined by Theorems 2 and 3.  
We shall first apply the general results on strictly quasi-regular Dirichlet forms and associated Markov processes (see Theorem V-2.13 and Proposition V-2.15 of  [M,R 92]) to our case. 
By the strictly quasi-regular Dirichlet form 
  $({\cal E}_{(\alpha)}, {\cal D}({\cal E}_{(\alpha)}))$ 
there exists a properly associated
$S$-valued Hunt process 
\begin{equation}
{\mathbb M} \equiv \Big(\Omega, {\cal F}, (X_t)_{t \geq 0}, (P_{\mathbf x})_{\mathbf x \in 
S_{\triangle}} \Big).
\end{equation}
$\triangle$ is a point adjoined to $S$ as an isolated point of $S_{\triangle} \equiv S \cup \{\triangle \}$. 
Let $(T_t)_{t \geq 0}$ be the strongly continuous contraction semigroup associated with 
$({\cal E}_{(\alpha)}, {\cal D}({\cal E}_{(\alpha)}))$, and $(p_t)_{t \geq 0}$ be the 
corresponding 
transition semigroup of kernels of the Hunt process $(X_t)_{t \geq 0}$, then for any $u \in {\cal F}C^{\infty}_0 
\subset {\cal D}({\cal E}_{(\alpha)})$ the following holds:
\begin{equation}
\frac{d}{dt} \int_{S} \big(p_t u \big)({\mathbf x}) \, \mu(d {\mathbf x}) 
= \frac{d}{dt} \big( T_t u, 1)_{L^2(S;\mu)} = {\cal E}_{(\alpha)}(T_t u, 1) = 0.
\end{equation}
By this, we see that
\begin{equation}
\int_S \big(p_t u\big)({\mathbf x}) \, \mu(d {\mathbf x}) = \int_S u({\mathbf x}) \, \mu(d {\mathbf x}), 
\quad \forall t \geq 0, \quad \forall u \in {\cal F}C^{\infty}_0,
\end {equation}
and hence, 
by the density of ${\cal F}C^{\infty}_0$ in $L^2(S;\mu)$, 
\begin{equation}
\int_S P_{\mathbf x} (X_t \in B) \, \mu(d {\mathbf x}) = \mu(B), \qquad \forall B \in {\cal B}(S), \quad \forall t \geq 0.
\end{equation}
Thus, we have proven the following Theorem 4.
\begin{theorem}
Let $0  < \alpha <2$, and let 
$({\cal E}_{(\alpha)}, {\cal D}({\cal E}_{(\alpha)}))$ be a strictly quasi-regular Dirichlet form on 
$L^2(S; \mu)$ that is defined through Theorem 2 or Theorem 3.
Then to $({\cal E}_{(\alpha)}, {\cal D}({\cal E}_{(\alpha)}))$, there exists a properly associated 
$S$-valued Hunt process (cf. Definitions IV-1.5, 1.8 and 1.13 of [M,R 92] for its
precise definition) ${\mathbb M}$ defined by (5.1), the invariant measure of which is $\mu$ (cf. (5.4)).\\
\quad
\qquad \qquad \qquad \qquad \qquad \qquad \qquad \qquad  \qquad \qquad
\qquad  \qquad \qquad  \qquad \qquad \quad 
\bsquare
\end{theorem}

We introduce below a standard procedure of application of Theorems 1, 2, 3 and 4 to the problem of  {\it{stochastic quantizations}} 
of Euclidean quantum fields, by means of the Hunt processes in Theorem 4. 
Mostly by the term  {\it{stochastic quantization}} one understands methods  to construct diffusion type  Markov processes (with continuous trajectories) which possess  as their invariant measures given probability 
measures $m$ 
associated to a given physical system (e.g., describing a statistical mechanical system or a quantum system, cf. [Parisi,Wu 81] and [A,Ma,R 2015], [A,Ma,R1 93], [A,Ma,R2 93],), with references therein.
Often it is interesting to consider their analogues where ${\mathbb R}^d$ is replaced by the $d$-dimensional 
torus ${\mathbb T}^d$ (since they can be used to approximate the 
 (Euclidean quantum) 
fields on ${\mathbb R}^d$). 
Moreover the original formulation of stochastic quantization can be naturally extended in the sense of asking for Markov processes having $m$ as invariant measure.
Here, we realize such Markov processes by the Hunt processes associated with 
$({\cal E}_{(\alpha)}, {\cal D}({\cal E}_{(\alpha)}))$ given by theorem 4, which are infinite dimensional analogues 
of $\alpha$-stable processes on finite dimensional space (cf. (2.8), (2.9) and (2.10), for the definition of the 
corresponding 
forms ${\cal E}_{(\alpha)}$).
At the end of the present section, examples based on Euclidean free fields, 
or the fields associated with the ${\Phi}^4_2$ and ${\Phi}^4_3$ models are indicated for $0 < \alpha \leq 1$.

\begin{remark} \qquad 
{\it{
Here we only briefly indicate the way these applications are obtained, 
detailed considerations on the individual 
 examples of the  Euclidean quantum field models, e.g., the fields 
constructed from 
 convoluted generalized white noise (cf.  [A,Gottschalk,Wu-J.L. 97]),  and the  free Euclidean field in all ${\mathbb R}^d$, the $P(\phi)_2$ field, 
 the $2$-dimensional fields with trigonometric and exponential potentials,  and 
of 
other random fields, will be carried out in a subsequent paper ([A,Kaga,KawaYaha,Y 2021] part 2: the applications).}}
\end{remark}

Euclidean  ({\it{scalar}}) quantum fields are expressed as random fields on ${\cal S}' \equiv {\cal S}'({\mathbb R}^d \to {\mathbb R})$, or
${\cal S}'({\mathbb T}^d \to {\mathbb R})$, the Schwartz space of real tempered distributions on 
the Euclidean space 
 ${\mathbb R}^d$, or the $d$-dimensional torus ${\mathbb T}^d$, with $d \geq 1$ a given space-time dimension, respectively.  Hence, each Euclidean quantum field is taken as a probability space $\displaystyle{ \big( {\cal S}', 
{\cal B}({\cal S}'), \nu \big)}$, where ${\cal B}({\cal S}')$ is the Borel $\sigma$-field of ${\cal S}'$ and $\nu$ is a Borel probability measure on ${\cal S}'$
(on ${\mathbb R}^d$, or on ${\mathbb T}^d$,
invariant under the 
rigid natural transformations of ${\mathbb R}^d$, or  ${\mathbb T}^d$, respectively,  
that correspond to the Euclidean group). 
In the construction  the standard 
theorem through which  probability measures $\nu$ are determined  is  the Bochner-Minlos theorem 
(cf. e.g., Section 3.2 of [Hida 80]).  

Let us 
 first 
recall the Bochner-Minlos theorem stated in a general framework.  Let $E$ be a nuclear space 
(cf., e.g., Chapters 47-51 of [Tr{\`e}ves 67]). Suppose in particular that $E$ is a countably Hilbert space,  $E$ is characterized by a sequence of Hilbert norms $\|{}\, \,{}\|_n$, $n \in {\mathbb N}\cup \{0\}$ such that 
$\|{}\, \,{}\|_0 < \|{}\, \,{}\|_1 < \cdots <\|{}\, \,{}\|_n < \cdots$.
Let $E_n$ be the completion of $E$ with respect to the norm $\|{}\, \,{}\|_n$, then 
%{\textcolor{blue}{
by definition
%}}
$E = \bigcap_{n \geq 0} E_n$ and $E_0 \supset E_1 \supset \cdots \supset E_n \supset \cdots$. Define 
$$E^{\ast}_n \equiv {\mbox{the dual space of $E_n$, 
%{\textcolor{blue}{
and assume the identification 
%}}
 $E^{\ast}_0 = E_0$ }},$$
then  we have 
$$
E \subset \cdots \subset E_{n+1} \subset E_n \subset \cdots \subset  E_0 = E^{\ast}_0 \subset 
\cdots \subset E^{\ast}_n \subset E^{\ast}_{n+1} \subset \cdots \subset E^{\ast}.
$$
Since $E$ is a nuclear space, for any $m \in {\mathbb N}\cup\{0\}$ there exists an $n 
\in {\mathbb N}\cup \{0\}$, $n > m$, such that the (canonical) injection $T^n_m\, : E_n \to E_m$ is a trace class (nuclear class) positive operator. The Bochner-Minlos theorem is given as follows:
\begin{theorem}{\bf{(Bochner-Minlos Theorem)}} \\
Let $C(\varphi)$, $\varphi \in E$, be a complex valued function on $E$ such that\\
i) \quad \,\, $C(\varphi)$ is continuous with respect to the norm $\| \, \cdot \, \|_{m}$ for some $m \in {\mathbb N} \cup \{0\}$; \\
ii) \quad \, ({\bf{positive definiteness}}) \quad 
for any $k \in {\mathbb N}$,
$$
\sum_{i,j =1}^k {\bar{\alpha}}_i \alpha_j C(
%{\textcolor{blue}{
{{\varphi}}_i 
%}}
- {\varphi}_j) \geq 0, \qquad 
\forall \alpha_i \in {\mathbb C}, \, \, \forall {\varphi}_i \in E, \,\,  i=1,\dots, k;
$$
(where ${\bar{\alpha}}$ means complex conjugate of $\alpha$). \\
iii) \quad ({\bf{normalization}}) \quad  $C(0) = 1$. \\
Then, there exists a unique Borel probability measure $\nu$ on $E^{\ast}$ such that 
$$
C(\varphi) = \int_{E^{\ast}} e^{i <\phi, \varphi>} \nu(d \phi), \qquad \varphi \in E.
$$
Moreover,
 for all $n > m$,
if the (canonical) injection $T^n_m\, : E_n \to E_m$ is a Hilbert-Schmidt operator, then 
the support of $\nu$ is in $E^{\ast}_n$, where 
$<\phi, \varphi> = {}_{E^{\ast}}<\phi, \varphi>_{E}$ is the dualization between $\phi \in E^{\ast}$ and 
$\varphi \in E$.\\
\qquad \qquad \qquad \qquad \qquad \qquad \qquad \qquad \qquad \qquad  \qquad \qquad
\qquad  \qquad \qquad  \qquad 
$\bsquare$
\end{theorem}

\begin{remark} {} \qquad
The assumption on the continuity of $C(\varphi)$ on $E$  given in i) of the above Theorem 5 can be replaced 
by the continuity of $C(\varphi)$ at the {\it{origin}} in  $E$, which is equivalent to i) under 
the assumption that $C(\varphi)$ satisfies ii) and iii) in Theorem 5 (cf. e.g., [It{\^o} K. 76]).
Namely, under the assumption of ii) and iii), the following is equivalent to i): \,
For any $\epsilon >0$ there exists a $\delta >0$ such that 
$$ | C(\varphi) -1 | < \epsilon, \qquad \forall \varphi \in E \quad {\mbox{with}} \quad 
\|\varphi\|_m < \delta.
$$
This can be seen as follows:
Assume that ii) and iii) hold.  For ii), let $k=3$, \,$\alpha_1 = \alpha$, \,$\alpha_2 = -\alpha$, \,$\alpha_3 = \beta$, \,$\varphi_1 = 0$, \,$\varphi_2 = \varphi$ and $\varphi_3 = \psi + \varphi$, 
 (for any $\alpha, \beta \in {\mathbb C}$, and any $\varphi, \psi \in E$)
then by the assumption ii),  the 
positive definiteness of $C$, we have
\begin{eqnarray*}
  \lefteqn{
    {\alpha}{\overline{\alpha}} \cdot ( 2C(0) - C(\varphi) - C(- \varphi))
} \\
 &&+ \alpha {\overline{\beta}} \cdot  (C(- \psi - \varphi) - C(- \psi)) 
 + 
    {\overline{\alpha}} \beta \cdot (C(\psi + \varphi) - C( \psi)) + 
\beta {\overline{\beta}} \cdot C(0) \geq 0.
\end{eqnarray*}
By making use of the fact that $C(- \varphi) = {\overline{C(\varphi)}}$, which follows from ii), and the assumption iii), from the above inequality we have
$$
0 \leq \det \left(
      \begin{array}{cc}
          2- C(\varphi) - {\overline{C(\varphi)}} \, \, &{} \, \, {\overline{C(\psi + \varphi) - C(\psi)}}\\
          C(\psi + \varphi) - C(\psi) \, \,  &{} \, \,1
       \end{array}
     \right).
$$
From this it follows that
$$ | C(\psi + \varphi) - C(\psi)|^2 \leq 2 \, |C(\varphi) - 1|.
$$
\qquad \qquad \qquad \qquad \qquad \qquad \qquad \qquad \qquad \qquad  \qquad \qquad
\qquad  \qquad \qquad  \qquad 
$\bsquare$
\end{remark}

By making use of the support property of $\nu$ by means of the Hilbert-Schmidt 
operators given by Theorem 5, we can present  a framework by which Theorems 1, 2, 3 and 4 
can be applied to the {\it{stochastic quantization}} of  Euclidean quantum fields.

We first define an adequate countably Hilbert nuclear space ${\cal H}_0 \supset {\cal S}({\mathbb R}^d \to {\mathbb R}) \equiv {\cal S}({\mathbb R}^d)$, for a given $d \in {\mathbb N}$ 
 (cf., e.g.,Appendix A.3 of [Hida 80],  Appendix 5 of [Hida,Kuo,Po,Str 93] for the framework on nuclear countable Hilbert space).
\, Let
\begin{equation}
{\cal H}_0 \equiv \Big\{ f \, : \, \|f\|_{{\cal H}_0} = \big((f,f)_{{\cal H}_0} \big)^{\frac12} < \infty, \, \, f : {\mathbb R}^d \to {\mathbb R}, \,\, {\mbox{measurable}} \Big\} \supset  
{\cal S}({\mathbb R}^d),
\end{equation}
where
\begin{equation}
(f,g)_{{\cal H}_0} \equiv (f,g)_{L^2({\mathbb R}^d)} = \int _{{\mathbb R}^d} f(x) g(x) \, dx.
\end{equation}
Let 
\begin{equation}
H \equiv (|x|^2 + 1)^{\frac{d+1}2} (- \Delta +1)^{\frac{d +1}2} (|x|^2 + 1 )^{\frac{d +1}2},
\end{equation}
\begin{equation}
H^{-1} \equiv (|x|^2 + 1)^{- \frac{d+1}2} (- \Delta +1)^{- \frac{d +1}2} (|x|^2 + 1 )^{- 
\frac{d +1}2},
\end{equation}
be the pseudo differential operators on 
${\cal S}'({\mathbb R}^d \to {\mathbb R}) \equiv {\cal S}'({\mathbb R}^d)$
with 
$\Delta$
the $d$-dimensional Laplace operator. 
Note that $H^{-1}$ is strictly positive bounded and symmetric, hence self-adjoint in $L^2({\mathbb R}^d)$.
For each $n \in {\mathbb N}$, define 
\begin{equation}
{\cal H}_n \equiv {\mbox{the completion of ${\cal S}({\mathbb R}^d)$ with respect to 
the norm $\|f \|_n$, \, $f \in {\cal S}({\mathbb R}^d)$}},
\end{equation}
where
$\|f \|^2_n \equiv (f,f)_n$ (in 
%{\textcolor{blue}{
the 
%}}
case where $n=1$, to denote the ${\cal H}_1$ norm we use the exact notation $\| {}\, \, \|_{{\cal H}_1}$,  in order to avoid a confusion with  the notation of some $L^1$ or $l^1$ norms) with 
%{\textcolor{blue}{
the corresponding scalar product
%}}
\begin{equation}
(f,g)_n = (H^n f, H^n g) _{{\cal H}_0}, \qquad f, \, g  \in {\cal S}({\mathbb R}^d).
\end{equation}
%{\textcolor{blue}{
Moreover we define, for $n \in {\mathbb N}$
\begin{equation}
{\cal H}_{-n} \equiv {\mbox{the completion of ${\cal S}({\mathbb R}^d)$ with respect to 
the norm $\|f \|_{-n}$, \, $f \in {\cal S}({\mathbb R}^d)$}},
\end{equation}
%}}
where
$\|f \|^2_{-n} \equiv (f,f)_{-n}$, with 
\begin{equation}
(f,g)_{-n} = ((H^{-1})^n f, (H^{-1})^n g) _{{\cal H}_0}, \qquad f, \, g  \in {\cal S}({\mathbb R}^d).
\end{equation}
Then  obviously, for  $f  \in {\cal S}({\mathbb R}^d)$, 
\begin{equation}
\|f \|_n \leq \|f \|_{n + 1}, \qquad \|f\|_{-n-1} \leq \|f\|_{-n},
\end{equation}
and by taking 
%{\textcolor{blue}{
the
%}}
 inductive limit and setting ${\cal H} = \bigcap_{n \in {\mathbb N}} {\cal H}_n$, 
we have the following inclusions:
\begin{equation}
{\cal H} \subset \cdots \subset {\cal H}_{n + 1} \subset {\cal H}_n \subset \cdots \subset  {\cal H}_0 \subset  \cdots \subset {\cal H}_{-n} \subset {\cal H}_{-n-1} \subset \cdots \subset {\cal H}^{\ast}.
\end{equation}
The (topological)  dual space of ${\cal H}_n$ is ${\cal H}_{-n}$, $n \in {\mathbb N}$.\\
By the operator 
$H^{-1}$ given  by (5.8) on ${\cal S}({\mathbb R}^d)$ we can define,
on each ${\cal H}_n$, $n \in {\mathbb N}$, the bounded symmetric 
(hence self-adjoint) operators  
\begin{equation}
 (H^{-1})^k, \quad k \in {\mathbb N}\cup\{0\}
\end{equation}
(we use the same notations for the operators on ${\cal S}({\mathbb R}^d)$ and on ${\cal H}_n$).
Hence, for the canonical injection 
\begin{equation}
T^{n +k}_n \, : \, {\cal H}_{n +k} \longrightarrow {\cal H}_n, \qquad k, \,n \in {\mathbb N}\cup \{0\},
\end{equation}
it holds that
$$\| T^{n +k}_{n} f \|_{n} = \| (H^{-1})^k f \|_{{\cal H}_0}, \qquad \forall f \in {{\cal H}_{n +k}},$$
where 
by a simple calculation by means of the Fourier transform, and by Young's inequality, we see that for each $n \in {\mathbb N}\cup \{0\}$, $H^{-1}$ on ${\cal H}_n$ is a Hilbert-Schmidt operator and hence 
$({H}^{-1})^2$ on ${\cal H}_n$ is a trace class operator.

Now,  by applying to the strictly positive self-adjoint 
Hilbert-Schmidt (hence compact) 
operator $H^{-1}$, 
 on ${\cal H}_0 = L^2({\mathbb R}^d \to {\mathbb R})$
  the {\it{Hilbert-Schmidt theorem}} (cf., e.g., Theorem VI 16, Theorem VI 22 of [Reed,Simon 80]) 
%{\textcolor{blue}{
we have that 
%}}
there exists an 
 orthonormal base (O.N.B.) $\{{\varphi}_i\}_{i \in {\mathbb N}}$ of ${\cal H}_0$ such that 
\begin{equation}
H^{-1} {\varphi}_i = \lambda_i \, {\varphi}_i, \qquad i \in {\mathbb N},
\end{equation}
where $\{\lambda_i\}_{i \in {\mathbb N}}$ 
%{\textcolor{blue}{
are
%}}
 the corresponding eigenvalues such that  
\begin{equation}
0 < \cdots < \lambda_2 < \lambda_1 \leq 1,
\quad 
{\mbox{which satisfy}} \quad  \,
\sum_{i \in {\mathbb N}} (\lambda_i)^2 < \infty, \quad {\mbox{i.e.,}} \quad 
\{\lambda_i\}_{i \in {\mathbb N}} \in l^2,
\end{equation}
and $\{{\varphi}_i\}_{i \in {\mathbb N}}$ is indexed adequately corresponding to the finite multiplicity of each $\lambda_i$, $i \in {\mathbb N}$.
By the definition (5.9), (5.10), (5.11) and (5.12) (cf. also (5.15)), for each $n \in {\mathbb N} \cup \{0\}$, 
\begin{equation}
\{(\lambda_i)^{n} {\varphi}_i \}_{i \in {\mathbb N}} \quad {\mbox{is an O.N.B. of ${\cal H}_n$ }}
\end{equation}
and 
\begin{equation}
\{(\lambda_i)^{-n} {\varphi}_i \}_{i \in {\mathbb N}} \quad {\mbox{is an O.N.B. of ${\cal H}_{-n}$ }}
\end{equation}
Thus, 
by denoting ${\mathbb Z}$ the set of integers, 
by the Fourier series expansion 
of functions in ${\cal H}_m$, $m \in {\mathbb Z}$ (cf. (5.9)-(5.12)),
such that 
for $f \in {\cal H}_m$, 
%{\textcolor{blue}{
we have
%}}
\begin{equation}
f  = \sum_{i \in {\mathbb N}} a_i (\lambda_i^{m} {\varphi}_i), \quad {\mbox{with}} \quad 
a_i \equiv \big(f, \, (\lambda_i^{m} {\varphi}_i) \big)_{m},
\, \, i \in {\mathbb N}
\end{equation}
%{\textcolor{blue}{
(in particular for $f \in {\cal S}({\mathbb R}^d) \subset {\cal H}_m$, it holds that 
\, $a_i =  {\lambda_i^{-m}} (f, \,
\varphi_i)_{{\cal H}_0}$).
%}}
%{\textcolor{blue}{
Moreover we have 
$$
\sum_{i \in {\mathbb N}} a_i^2 = \|f\|_m^2,
$$
that yields
%}}
 an {\it{isometric isomorphism}}  $\tau_m$ for each $m \in {\mathbb Z}$ such that 
\begin{equation}
\tau_m \, : \, {\cal H}_m \ni f \longmapsto ({\lambda}_1^m a_1, {\lambda}_2^m a_2, \dots) \in l^2_{(\lambda_i^{-2m})},
\end{equation}
where  $l^2_{(\lambda_i^{-2m})}$ is the weighted $l^2$ space defined by (2.1) with $p=2$, 
and  $\beta_i  = \lambda_i^{-2m}$.
Precisely, for 
$f  = \sum_{i \in {\mathbb N}} a_i (\lambda_i^{m} {\varphi}_i) \in {\cal H}_m$ and 
$g = \sum_{i \in {\mathbb N}} b_i (\lambda_i^{m} {\varphi}_i)  \in {\cal H}_m$, with 
$a_i \equiv \big(f, \, (\lambda_i^{m} {\varphi}_i) \big)_{m}$, 
$b_i \equiv \big(g, \, (\lambda_i^{m} {\varphi}_i) \big)_{m}$, $ i \in {\mathbb N}$, 
by $\tau_m$ the following holds (cf. (5.19) and (5.20)):
$$
(f, \,g)_{m} = \sum_{i \in {\mathbb N}} a_i \cdot b_i = 
\sum_{i \in {\mathbb N} } {\lambda_i^{-m}} (\lambda_i^m a_i) \cdot {\lambda_i^{-m}} (\lambda_i^m b_i) 
 = \big( {\tau}_m f, \, {\tau}_m g \big)_{l^2_{(\lambda_i^{-2m})}}.
$$
By the map $\tau_m$ we can identify, in particular, the two systems of Hilbert spaces 
given by (5.23) and (5.24) through the following diagram:
%{\textcolor{blue}{
%\begin{equation}
%{\cal H}_2 \, \, \subset \, \,  \, \, {\cal H}_1 \, \, \subset {\cal H}_0 = L^2({\mathbb R}^d) \subset {\cal H}_{-1} \subset {\cal H}_{-2},
%\end{equation}
%{} \qquad \qquad \qquad \qquad  \qquad   $\parallel$ \qquad \quad $\parallel$
%\qquad \qquad  $\parallel$ \qquad \quad \, \,  $\parallel$ \qquad \, \,  $\parallel$
%\begin{equation}
% l^2_{(\lambda_i^{-4})} \, \, \subset  l^2_{(\lambda_i^{-2})} \, \, \subset \, \, \, \, \, \, \, \, l^2 \, \, \, \, \, \, \, \, \,
%\subset  \, \, \, \, \, \, \,  l^2_{(\lambda_i^2)} \subset  \, \, \, \, \, \, \,  l^2_{(\lambda_i^{4})}.
%\end{equation}
%\newpage
\begin{equation}
\cdots \quad {\cal H}_2 \, \, \subset \, \,  \, \, {\cal H}_1 \, \, \subset {\cal H}_0 = L^2({\mathbb R}^d) \subset {\cal H}_{-1} \subset \, {\cal H}_{-2} \quad \cdots,
\end{equation}
{}\, \, \quad  \qquad \quad
\quad \qquad \qquad \quad   $\parallel$ \qquad \quad $\parallel$
\qquad \qquad  $\parallel$ \qquad \quad \, \,  $\parallel$ \qquad \, \, \, $\parallel$
\begin{equation}
{} \quad 
 \cdots \quad l^2_{(\lambda_i^{-4})} \, \, \subset  l^2_{(\lambda_i^{-2})} \, \, \subset \, \, \, \, \, \, \, \, l^2 \, \, \, \, \, \, \, \, \,
\subset  \, \, \, \, \, \,  l^2_{(\lambda_i^2)} \subset  \, \, \, \, \, \,   l^2_{(\lambda_i^{4})} \quad \cdots.
\end{equation}
%}}
{\bf{Example 1. (The Euclidean free fields)}} \quad Let ${\nu}_0$ be the Euclidean free field 
%{\textcolor{blue}{
probability
%}}
measure on ${\cal S}' \equiv {\cal S}'({\mathbb R}^d)$.
Then,
the (generalized) characteristic function  $C(\varphi)$ of $\nu_0$ in Theorem 5 is given by  
\begin{equation}
C(\varphi) = \exp ({- \frac12 (\varphi, (- \Delta + m^2_0)^{-1} \varphi)_{L^2({\mathbb R}^d)}}),
\qquad {\mbox{for \, $\varphi \in {\cal S}({\mathbb R}^d \to {\mathbb R})$}}.
\end{equation}
Equivalently, $\nu_0$ is a centered Gaussian probability measure on ${\cal S}'$, the covariance of which 
  is given by
\begin{equation*}
\int_{{\cal S}'} <\phi, {\varphi}_1> \cdot  <\phi, {\varphi}_1> \, \nu_0(d \phi) 
= \big(\varphi_1, (- \Delta + m^2_0)^{-1} \varphi_2 \big)_{L^2({\mathbb R}^d)}, 
\end{equation*}
\begin{equation}
{\mbox{for any}} \quad 
 \varphi_1, \, \varphi_2 \in {\cal S}({\mathbb R}^d \to {\mathbb R}),
\end{equation}
where $\Delta$ is the $d$-dimensional Laplace operator and $m_0 >0$ (for $d \geq 3$, we can allow also $m_0 =0$) is a given mass for this  scalar field (the coordinate process to $\nu_0$). By (5.25),  the functional $C(\varphi)$ is continuous with respect to the  norm of the space  ${\cal H}_0 =L^2({\mathbb R}^d)$, and the kernel of 
$(- \Delta + m^2_0)^{-1}$, which is the Fourier inverse transform of $(|\xi|^2 + m^2_0)^{-1}$, $\xi \in {\mathbb R}^d$, is explicitly given by Bessel functions (cf., e.g., section 2-5 of  [Mizohata 73]).
Then, by 
Theorem 5 
and (5.16) 
the support of $\nu_0$ 
can be taken to be 
 in the  Hilbert  spaces  ${\cal H}_{-n}$, $ n \geq 1$ (cf. (5.23) and (5.24)).

Let us apply Theorems 1, 2 and 4 with $p = \frac12$ to this random field. 
For a clarity of the discussion, we start the  consideration from the case where $\alpha= 1$.
 Then, we shall state the corresponding results for the cases where $0 < \alpha <1$.

Now, we take ${\nu}_0$ as a Borel probability measure on ${\cal H}_{-2}$. By (5.22), (5.23) and (5.24), 
by taking $m = -2$, 
$\tau_{-2}$ defines an isometric isomorphism such that 
\begin{equation*}
\tau_{-2} \, : \, {\cal H}_{-2} \ni f \longmapsto ({\lambda}_1^{-2}a_1, {\lambda}_1^{-2}a_2, \dots) \in l^2_{(\lambda_i^{4})},
\end{equation*}
\begin{equation}
{\mbox{with}} \quad 
a_i \equiv (f, \, \lambda_i^{-2} {\varphi}_i)_{-2}, \, \, \, i \in {\mathbb N}.
\end{equation}
Define a probability measure $\mu$ on 
$ l^2_{(\lambda_i^{4})}$ such that 
$$\mu(B) \equiv \nu_0 \circ \tau^{-1}_{-2}(B) \quad {\mbox{for}} \quad B \in {\cal B}(l^2_{(\lambda_i^{4})}).
$$
We set $S = l^2_{(\lambda_i^{4})}$ in Theorems 1, 2 and 4, then it follows that the weight $\beta_i$ satisfies $\beta_i = \lambda_i^4$. We can take ${\gamma_i}^{-\frac12} = \lambda_i$ in Theorem 2-i) 
with $p =2$,  then, from (5.18) we have
\begin{equation}
\sum_{i=1}^{\infty} 
\beta_i \gamma_i \cdot
\mu \Big( \beta_i^{\frac12} |X_i| > M \cdot \gamma_i^{-\frac{1}2}\Big) 
\leq \sum_{i =1}^{\infty} \beta_i \gamma_i = \sum _{i =1}^{\infty} (\lambda_i)^2 < \infty
\end{equation}
(5.28) shows that the condition (4.3) holds 
for $p=2$ and $\alpha =1$.
%{\textcolor{blue}{
Also, as has been mentioned above, 
since $\nu_0({\cal H}_{-n}) = 1$,  for any   $ n \geq 1$, we have 
$$1 = \nu_0({\cal H}_{-1}) = \mu(l^2_{(\lambda_i^2)}) = 
\mu \big( \bigcup_{M \in {\mathbb N}} \{ |X_i| \leq M \beta_i^{- \frac12} \gamma_i^{- \frac12}, \, \forall i \in {\mathbb N} \} \big),
$$
$$
 {\mbox{for $\beta_i = \lambda_i^4$, \, $\gamma_i^{- \frac12} = \lambda_i$}}.
$$
This shows that the condition (4.4) is satisfied for $p=2$.
%}}

Thus, by Theorem 2-i) and Theorem 4,  
for   $\alpha = 1$,
there exists an  $l^2_{((\lambda_i)^{4})}$-valued  Hunt process 
${\mathbb M} \equiv \big(\Omega, {\cal F}, (X_t)_{t \geq 0}, (P_{\mathbf x})_{\mathbf x \in 
S_{\triangle}} \big)$, associated  to the non-local Dirichlet form
$({\cal E}_{(\alpha)}, {\cal D}({\cal E}_{(\alpha)}))$.
We can now define an
${\cal H}_{-2}$-valued 
 process 
$(Y_t)_{t \geq 0}$ 
such that 
$$(Y_t)_{t \geq 0} \equiv \big({\tau}^{-1}_{-2}(X_t) \big)_{t \geq 0}.$$
Equivalently, by (5.27) for $X_t = (X_1(t), X_2(t), \dots) \in l^2_{(\lambda_i^4)}$, $P_{\mathbf x}-a.e.$, by setting $A_i(t)$ such that 
 $A_i(t) \equiv \lambda_i^2 X_i(t)$  (cf. (5.21) and (5.22)),  we see that  
$Y_t$ \, is also given by 
\begin{equation}
Y_t = \sum_{i \in {\mathbb N}} A_i(t) (\lambda_i^{-2} \varphi_i) = \sum_{i \in {\mathbb N}} X_i(t) \varphi_i \in {\cal H}_{-2}, \qquad \forall t \geq 0, \, \, P_{\mathbf x}-a.e..
\end{equation}
By (5.4) and (5.27), $Y_t$ is an ${\cal H}_{-2}$-valued Hunt process that is a {\it{stochastic quantization}} (according to the definition we gave to this term) with respect to the 
non-local Dirichlet form 
$({\tilde{\cal E}}_{(\alpha)}, {\cal D}({\tilde{\cal E}}_{(\alpha)}))$ 
on $L^2({\cal H}_{-2}, \nu_0)$, that is defined through
$({\cal E}_{(\alpha)}, {\cal D}({\cal E}_{(\alpha)}))$, by making use of $\tau_{-2}$. 
This holds for $\alpha = 1$.  

For the cases where $0 < \alpha <1$ and $p = 2$, we can also apply Theorems 1, 2 and 4, and then have the corresponding result to (5.29).  For this purpose we have only to notice that for  $\alpha \in (0, 1)$ if we  take $\nu_0$  as a Borel probability measure on ${\cal H}_{-3}$, and set 
 $S \equiv l^2_{(\lambda_i^6)}$, $\beta_i \equiv \lambda_i^6$, $\gamma_i \equiv \lambda _i$, and define 
\begin{equation*}
\tau_{-3} \, : \, {\cal H}_{-3} \ni f \longmapsto ({\lambda}_1^{-3}a_1, {\lambda}_1^{-3}a_2, \dots) \in l^2_{(\lambda_i^{6})},
\end{equation*}
\begin{equation*}
{\mbox{with}} \quad 
a_i \equiv (f, \, \lambda_i^{-3} {\varphi}_i)_{-3}, \, \, \, i \in {\mathbb N},
\end{equation*}
(cf. (5.22), (5.23), (5.24) and (5.27)),
then 
$$\sum_{i=1}^{\infty} (\beta_i \gamma_i)^{\frac{\alpha +1}2} = \sum_{i=1}^{\infty} (\lambda_i)^{2(\alpha +1)}
 < \infty.
$$
As a consequence, 
for $\alpha \in (0, 1)$, 
we adopt here the same formulation as will be given in the next example, Example 2  (cf. (5.58)-(5.62)), 
we then see that  by  this setting (4.3) and (4.4) also hold 
(cf. (5.28) and the formula given below (5.28)), and have an analogue of (5.29).

The case $1 < \alpha < 2$ requires a separate consideration, see [A,Kaga,Kawa,Yaha,Y 2021]. The diffusion case $\alpha =2$ was already discussed in [A,R 89],
[A,R 91] (and references therein).
\\
\qquad \qquad \qquad \qquad \qquad \qquad \qquad \qquad \qquad  \qquad \qquad
\qquad  \qquad \qquad  \qquad \qquad \quad 
\bsquare

\begin{remark} {}\quad ({\bf{A structural property of Gaussian random fields}}) \quad 
The Euclidean free field considered in Example 1 is a Gaussian field, which possesses a simple mathematical structure as follows:\, 
Let $D(x-y)$ denotes the integral kernel corresponding to the pseudo differential operator 
$(- \Delta + m_0)^{-1}$ on ${\cal S}({\mathbb R}^d \to {\mathbb R})$ considered in Example 1
(cf. (5.26) and the explanation follwing to (5.26)).
Then $(\varphi_1, (- \Delta + m_0)^{-1} \varphi_2)_{L^2({\mathbb R}^d)} 
= 
\int_{{\mathbb R}^d \times {\mathbb R}^d} \varphi_1(x) \, D(x-y)\, \varphi_2(y) \, dx dy$.
$D(x-y)$ is often refered to as the {\it{free propagator}}. For a notational simplicity, let us denote 
$\int_{{\mathbb R}^d \times {\mathbb R}^d} \varphi_i(x) \, D(x-y)\, \varphi_j(y) \, dx dy$, 
$\varphi_i \in {\cal S}({\mathbb R}^d \to {\mathbb R})$, $i=1, \dots, n$,  $n \in {\mathbb N}$,
by $a_{i,j}$. Then it holds that (cf. (5.26)), for $n \in {\mathbb N}$, 
$$
\int_{{\cal S}'} \prod_{i = 1}^{2n} <\phi, {\varphi}_i> \nu_0(d \phi) 
= \sum_{\rm{pairing}} \prod a_{i,j}, \quad 
\int_{{\cal S}'} \prod_{i = 1}^{2n-1} <\phi, {\varphi}_i> \nu_0(d \phi) 
= 0,
$$
where $ \sum_{\rm{pairing}}$ denotes the sum of all the combinations of the {\it{distinguished pairings}} in $\{1, \dots,  2n\}$,
 and $\prod a_{i,j}$ is the product of $a_{i,j}$ characterized by each pairing.
Precisely, each {\it{pairing}} is determined as follows: \, Firstly, as the first pair we choose two indices 
 from $\{1, \dots, 2n\}$,
secondly, 
as the second pair we
choose two indices  from $\{1, \dots, 2n \} \setminus \{{\mbox{{\it{previously chosen indices}}}} \}$, thirdly, 
as the third pair 
 we choose two indices  from $\{1, \dots, 2n \} \setminus \{ {\mbox{{\it{previously chosen indices}}}} \}$, and so on.  
Among all the pairings gotten through the above procedure, 
as a consequence, we have $(2n -1)!!$ number of 
{\it{distinguished pairings}}, where 
$$
(2n-1)!! \equiv (2n-1) \cdot (2n-3) \cdots 3 \cdot 1.$$
Thus, as a special case,  it holds that
$$
\int_{{\cal S}'}  <\phi, {\varphi}>^{2n} \nu_0(d \phi) 
= (2n -1)!! \, ((\varphi, (- \Delta + m_0)^{-1} \varphi)_{L^2({\mathbb R}^d)})^n, \quad \varphi \in 
{\cal S}({\mathbb R}^d \to {\mathbb R}).$$
These properties expressing even moments of mean zero Gaussian random field with sum of products of second moments are characteristic of the fields being Gaussian.
\qquad \qquad \qquad \qquad \qquad \qquad \qquad \qquad \qquad  \qquad \qquad
\qquad  \qquad \qquad  \qquad \qquad \quad 
\bsquare
\end{remark}

{\bf{Example 2. (The Euclidean ${\Phi}^4_1$, ${\Phi}^4_2$ and ${\Phi}^4_3$ fields)}} \quad
As  a rough explanation, the Euclidean $\Phi_d^4$, $d=1,2,3$ fields are  probability measures on ${\cal S}'({\mathbb R}^d \to {\mathbb R})$, 
that describe Euclidean invariant random fields 
  having self interactions with fourth power, and are defined by adding corresponding interaction terms to the Euclidean free fields on ${\cal S}'({\mathbb R}^d \to {\mathbb R})$, 
$d=1,2,3$.  
%{\textcolor{blue}{
(cf., e.g., [Iwata 87] for the case $d=1$; 
%}}
[Simon 74], [Glimm,Jaffe 87] for fields 
in the case $d=2$ 
with general polynomial interactions called $P(\Phi)_2$ (Euclidean) field,  for trigonometric and exponential see also, [A,H-K 76], [A,H-K1 77], [A,H-K2 77], [A,H-K,Zegarlinski 89], [Zegarlinski 84], [A,Y 2002] and references therein;
 for the ${\Phi}^4_3$-field see the references in [Glimm,Jaffe 87]
and [Hairer 2014], [A,Kusuoka-sei 2020], [A.Kusuoka-sei 2021], [Gu,Ho 2019]).  There are several strategies 
 (to cope with the singularities of the fields) 
through which the $\Phi_d^4$, $d=2,3$, Euclidean fields are constructed. 
A standard construction strategy is to start  
 from (continuum) random fields 
on ${\cal S}'({\mathbb R}^d \to {\mathbb R})$, 
 associated with a bounded region in ${\mathbb R}^d$ ( or a torus ${\mathbb T}^d$) and then expanding to ${\mathbb R}^d$, 
or alternative, to 
 start from 
 random 
fields on a lattice $(\epsilon {\mathbb Z})^d$ (with the lattice spacing $\epsilon >0$, which subsequently tends to $0$, i.e., taking a continuum limit).
In both cases, for $d=2,3$, renormalization counter terms are required for a non trivial limit.

The stochastic quantization of $\Phi_d^4$, $d =2,3$, is also not a  trivial problem.  
%{\textcolor{blue}{
In fact
for $d=2$ it was obtained first by the  local Dirichlet form method in the 90's (cf. [A,R 89,90,91]) and in the sense of strong solutions 
in [Da Prato,Debussche] see also [R,Zhu,Zhu 2015].
For $d=3$ it has been an open problem until the publication [Hairer 2014] for the diffusion case on $\Phi^4_3$ Euclidean field on the $3$-dimensional torus ${\mathbb T}^3$.
After this, several publications obtained the stochastic quantization on 
${\mathbb R}^3$ (see [Gu,Ho 2019] and also references in [A,Kusuoka-sei 2020] and 
[A,Kusuoka-sei 2021]).
The methods of the present paper show that a stochastic quantization of $\Phi^4_3$ Euclidean field on ${\mathbb R}^3$ can also be realized by {\it{non-local}} Hunt processes. 
%}}
  In order to understand the difficulty 
 of the stochastic quantization program,
 which is caused by  singularities of 
the  $\Phi^4_d$ Euclidean field measures, we briefly recall above construction procedures.

The first one is as follows: \, Let $\nu_0$ be the Euclidean free field measure on 
${\cal S}' \equiv {\cal S}'({\mathbb R}^d \to {\mathbb R})$, $d = 2,3$, (obtained in  Example 1). For $d=2$, let $:Z^4:$ be the ${\cal S}'$-valued random variable on 
$({\cal S}', {\cal B}({\cal S}'), \nu_0)$ uniquely defined by (see, e.g., [Simon 74], 
and 
for a definition by means of the multiple stochastic integrals, cf., e.g. [A,Y 2002], 
[A,Ferrario,Y 2004]).
also cf. Remark 7)
\begin{eqnarray}
\lefteqn{
\int_{{\cal S}'} {}_{{\cal S}'}<:Z^4:, h>_{{\cal S}} \,
: {}_{{\cal S}'}<\phi, \varphi_1>_{{\cal S}} \cdots  {}_{{\cal S}'}<\phi, \varphi_4>_{{\cal S}}: \nu_0(d \phi) 
}  \\
&&= 
\int_{{\mathbb R}^2} \prod_{j = 1}^4 \Big( \int_{{\mathbb R}^2} (- \Delta + m_0^2)^{-1} 
(x-y_j) \varphi_j(y_j) dy_j \Big) h(x) dx, \quad h, \, \varphi_j \in {\cal S}, \, j=1, \dots, 4,
\nonumber
\end{eqnarray}
where
$h, {\varphi}_j \in {\cal S}, \, j=1,\cdots,4$, 
 the Wick power $: {}_{{\cal S}'}<\phi, \varphi_1>_{{\cal S}} \cdots  {}_{{\cal S}'}<\phi, \varphi_4>_{{\cal S}}:$ is defined through the Hermite polynomials, e.g., a simple case 
$: ({}_{{\cal S}'}<\phi, \varphi>_{{\cal S}})^4: = \sum_{n=0}^2 \frac{4!}{n! (4-2n)!} 
 ({}_{{\cal S}'}<\phi, \varphi>_{{\cal S}})^{4-2n} (- \frac{a}2)^n$ with $a = 
\int_{{\cal S}'} ({}_{{\cal S}'}<\phi, \varphi>_{{\cal S}})^2 \nu_0(d \phi)$.
 For $d=3$, by putting 
 a momentum {\it{cut-off}} ${\kappa}$ and 
an additional counter term (i.e., an additional mass renormalization, diverging when $\kappa$ is removed), 
similar to (5.30), an ${\cal S}'({\mathbb R}^3 \to {\mathbb R})$-valued random variable $:Z^4_{\kappa}:$ can be defined (cf. [Feldman 74], [Feldman,Osterwalder 76]). 

In the case where $d=2$, for bounded regions $\Lambda \subset {\mathbb R}^2$,
 with $I_{\Lambda}$ the corresponding characteristic function, 
 and a {\it{coupling constant}} $\lambda \geq 0$,
\begin{equation}
\nu_{\Lambda}(d \phi)  \equiv \frac{\exp(- \lambda <:Z^4:, I_{\Lambda}>)}{
\int_{{\cal S}'} \exp(- \lambda <:Z^4:, I_{\Lambda}>) \nu_0(d \phi)} \nu_0(d \phi), 
\end{equation}
can be shown to be well defined (see [Ne 66], [Simon 74], [Glimm,Jaffe 87]). 
Moreover, the existence of a weak limit 
 $\nu^{\ast}$ (unique for small values of $\lambda$) 
such that
\begin{equation}
\nu^{\ast} \equiv \lim_{\Lambda \uparrow {\mathbb R}^2} \nu_{\Lambda}
\end{equation}
has been proven 
(cf. e.g., [Simon 74], [Glimm,Jaffe 87],  
 see also e.g.,
[A,R 89, 90, 91]).
The probability measure 
$\nu^{\ast}$ on ${\cal S}'({\mathbb R}^2 \to {\mathbb R})$ defined by (5.32) is known as the $\Phi^4_2$ Euclidean field measure.

%{\textcolor{blue}{
Similarly, in case  where $d=3$, 
 by [Feldman,Osterwalder 76], 
the existence of a weak limit 
point ${\nu}^{\ast}$ 
(for adequately small $\lambda \geq 0$) such that 
\begin{equation}
\nu^{\ast} \equiv  \lim_{\Lambda \uparrow {\mathbb R}^3} \lim_{\kappa \to 1} \nu_{\Lambda, \kappa}, 
\end{equation}
where
\begin{equation}
\nu_{\Lambda, \kappa} \equiv \frac{\exp(- \lambda <:{\tilde Z}^4_{\kappa}:, I_{\Lambda}>)}{
\int_{{\cal S}'} \exp(- \lambda <:{\tilde Z}^4_{\kappa}:, I_{\Lambda}>) \nu_0(d \phi)} \nu_0(d \phi),
\end{equation}
with 
$$
:{\tilde Z}^4_{\kappa}: \equiv :Z^4_{\kappa}: + \lambda a(\kappa) - \lambda^3 b(\kappa), \, \, {\mbox{for suitable constants \,  $a, \, b$ \, depending on $\kappa$}}
$$
(cf. [Feldman 74], [Feldman,Osterwalder 76] and [A,Kusuoka-sei 2020], [A,Kusuoka-sei 2021], 
[Gu,Ho 2019], [Gu,Ho 2019],  also 
(5.37) given below for a lattice formulation) 
has been proven.
The probability measure 
$\nu^{\ast}$ on ${\cal S}'({\mathbb R}^3 \to {\mathbb R})$
defined by (5.33) is known as "the $\Phi^4_3$ Euclidean measure" (the uniqueness of limit points for sufficiently small value of $\lambda \geq 0$ is conjectured).
%}}

%{\textcolor{blue}{
For $d=3$ the power $4$ is a critical point 
 for the existence of corresponding probability measures, 
and the analytic data of the $\Phi^4_3$ Euclidean field measure are quite singular.
The Euclidean invariance is assured by the considerations in [Feldman,Osterwalder 76] (cf. also [Magnen,S{\'e}n{\'e}or 76], [Seiler,Simon 76], [Glimm,Jaffe 87] and references therein). 
For $d=4$ there is no affirmative result on existence of ${\Phi}^4_4$ model, 
%{\textcolor{blue}{
see however, e.g., [Glimm,Jaffe 74], [Fr{\"o}hlich,Park 77], 
[Fr{\"o}hlich], [Aiz1], [Aiz2], [GaR], 
[A 2016], [Gu,Ho 2019], [K.R.It{\^o} 89].
%}}
%}}

The alternative procedure, through which the $\Phi^4_d$ Euclidean field measures are defined, is the following: \, Let $d= 2,3$. For each bounded region $\Lambda \subset 
{\mathbb R}^d$ and a lattice spacing $\epsilon >0$, let 
\begin{equation}
L_{\epsilon, \Lambda} \equiv (\epsilon {\mathbb Z})^d \cap \Lambda,
\end{equation}
 and define a family of real valued random variables $\phi \equiv \{ \phi(x) \, : \, x \in L_{\epsilon, \Lambda} \}$, the probability distribution of which is given by
\begin{equation}
\nu_{\epsilon, \Lambda} \equiv \frac{1}{Z_{\epsilon, \Lambda}} \prod_{x \in L_{\epsilon, \Lambda}} e^{- S_{L_{\epsilon, \Lambda}} (\phi)} d \phi(x),
\end{equation}
where $Z_{\epsilon, \Lambda}$ is the normalizing constant, 
\begin{equation}
S_{L_{\epsilon, \Lambda}} (\phi) \equiv \frac12 \sum_{<x,y>} \epsilon^{d-2} \big(\phi(x) -\phi(y) \big)^2 + \frac12 a_{\epsilon} \sum_{x \in L_{\epsilon, \Lambda}} \epsilon^d \phi^2(x)
+ \frac{\lambda}2 \sum_{x \in L_{\epsilon, \Lambda}} \epsilon^d \phi^4(x),
\end{equation}
with $a_{\epsilon}$ a counter term depending on $\epsilon >0$ and $d=2, 3$, $\lambda \geq 0$ a coupling constant;
  $<x,y>$ denotes the nearest neighbor points in $L_{\epsilon, \Lambda}$. 
In [Brydges,Fr{\"o}hlich,Sokal 83] it is shown, roughly speaking
that, for adequately small $\lambda \geq 0$,  there exists a subsequence 
$\{ \nu_{\epsilon_i, \Lambda_j}\}_{i,j \in {\mathbb N}}$ 
of $\{ \nu_{\epsilon, \Lambda}\}_{\epsilon >0, \Lambda \subset {\mathbb R}^d}$ 
with $\lim_{i \to \infty} \epsilon_i = 0$ and $\lim_{j \to \infty} \Lambda_j = {\mathbb R}^d$, 
and a weak limit 
\begin{equation}
\nu^{\ast} \equiv \lim_{i \to \infty} \lim_{j \to \infty} \nu_{\epsilon_i, \Lambda_j}
\end{equation}
exists in the space of Borel probability measures on ${\cal S}'({\mathbb R}^d \to {\mathbb R})$ by interpreting $\nu_{\epsilon, \Lambda}$ as 
an element in this space,
 for each $\epsilon >0$ and $\Lambda \subset 
{\mathbb R}^d$  (cf., the subsequent precise discussions from (5.39) to (5.56), for the weak convergence). 
For each $d=2,3$, 
$\nu^{\ast}$ given by (5.38) is a $\Phi^4_d$, $d=1,2,3$, Euclidean field measure, defined through a lattice approximation, see [Park 75], [Park 77], [Sokal 82], [Glimm,Jaffe 74].

The present example is  formulated by using the  $\Phi^4_d$ Euclidean field measure 
on ${\cal S}'({\mathbb R}^d \to {\mathbb R})$ constructed through the lattice approximation described above. 
For this purpose, we need to certify the support properties of the measures with more details. 
Precisely, to apply Theorems 1, 2, 3 and 4 
 (cf. Example 1 and Theorem 5) to the example, we have to know that the supports of these measures are in some Hilbert spaces (see (5.48) below).

For each $\epsilon >0$, and bounded region $\Lambda \subset {\mathbb R}^d$, $d=1,2,3$, and for $F(\phi)$ a polynomial in $\{ \phi(x) \, : \, x \in L_{\epsilon, \Lambda} \}$, by (5.35) and (5.36) let 
\begin{equation}
< F >_{\epsilon, \Lambda} \equiv \int_{{\mathbb R}^{N(\epsilon, \Lambda)}} F(\phi) \nu_{\epsilon, \Lambda}(d \phi), 
\end{equation}
where 
$N(\epsilon, \Lambda)$ is the cardinality of $L_{\epsilon,\Lambda}$.
By [Sokal 82] (cf. also section 2 of [Brydges,Fr{\"o}hlich,Sokal 83]), the following limit exists:
\begin{equation}
<F>^{(\epsilon)} \equiv \lim_{\Lambda \uparrow {\mathbb R}^d} < F >_{\epsilon, \Lambda}.
\end{equation}
Also, for each $\epsilon >0$, and $d=1,2,3$, there exists a weak limit 
\begin{equation}
\nu_{\epsilon} \equiv \lim_{\Lambda \uparrow {\mathbb R}^d} \nu_{\epsilon, \Lambda},
\end{equation}
that is a Borel probability measure on ${\cal S}'({\mathbb R}^d \to {\mathbb R})$.  
Let $D^{(\epsilon)}$ be the two point function of the lattice field with the lattice spacing 
$\epsilon >0$ and a given mass $m_0 > 0$, such that 
for $x, y \in (\epsilon {\mathbb Z})^d$  
\begin{equation*}
D^{(\epsilon)} (x-y) = (2 \pi)^{-d} \int_{[- \frac{\pi}{\epsilon}, \frac{\pi}{\epsilon}]^d} 
\big( 2\epsilon^{-2} \sum_{i=1}^d ( 1 - \cos \epsilon k_i) + m_0^2 \big)^{-1} 
e^{i k \cdot (x-y)} dk_1 \cdots dk_d,
\end{equation*}
which is the
 lattice version of the 
corresponding covariance 
operator $(- \Delta + m_0^2)^{-1}$  for the continuous 
Euclidean free field model (cf. (5.25), (5.26) and (5.30)), 
i.e., the so called {\it{free propagator}} 
(cf. Remark 7) 
which is denoted by $C^{\epsilon}(x-y)$ in 
[Brydges,Fr{\"o}hlich,Sokal 83].  By (5.40), define 
\begin{equation}
S^{(\epsilon)} (x-y) \equiv <\phi(x) \cdot \phi(y)>^{(\epsilon)}, \qquad x, y \in (\epsilon {\mathbb Z})^d,
\end{equation}
and 
\begin{equation}
S^{(\epsilon)}_n(x_1, \dots, x_n) \equiv < \prod_{i=1}^n \phi(x_i) >^{(\epsilon)}, \quad x_i \in L_{\epsilon}, \, \, i=1, \dots, n, \, \, n \in {\mathbb N}.
\end{equation}
From Theorem 6.1, Lemma A.1 with the formulas (A.13), (A.16), (A.17)  
and (8.2) in [Brydges,Fr{\"o}hlich,Sokal 83], where  
 $m_0^2$ is taken to be equal $1$ and 
by setting $\alpha =0$,
 we see that there exist  universal constants 
$\lambda_0, K_1, K_2$ such that if $0 \leq \lambda \leq \lambda_0$, then 
we have the inequalities 
\begin{equation}
|\| S^{(\epsilon)} - D^{(\epsilon)} \||^{(\epsilon)} \leq K_1 \lambda^2, \qquad \forall \epsilon >0,
\end{equation}
and 
\begin{equation}
\| D^{(\epsilon)} |_{l^1}^{(\epsilon)} \leq K_2 m_0^{-2}, \qquad \forall \epsilon >0.
\end{equation}
Moreover  the {\it{Gaussian inequality}}
\begin{equation*}
0 \leq S^{(\epsilon)}_{2n} (x_1, \dots, x_{2n}) \leq \sum_{\rm{pairing}} \prod S^{(\epsilon)}(x_i, x_j), 
\end{equation*}
\begin{equation}
S^{(\epsilon)}_{2n-1} (x_1, \dots, x_{2n-1}) = 0, \quad n \in {\mathbb N},
\end{equation}
holds, 
where
the indices $i$ and $j$ move in $\{1, \dots, 2n\}$ and the notation such that $\sum_{{\rm{pairing}}} \prod$ is defined in Remark 7, and
\begin{equation}
|\| f \||^{(\epsilon)} \, \equiv \, \|f \|_{l^1}^{(\epsilon)} + \|f \|_{\infty}^{(\epsilon)} \, \equiv \, \epsilon^d \sum_{x \in (\epsilon {\mathbb Z})^d} |f(x)| + \sup_{x \in (\epsilon {\mathbb Z})^d} |f(x)|.
\end{equation}
{\it{By Remark 7, we note that the inequality in (5.46) holds as the equality for the Euclidean free field model}}. 
By making use of (5.44), (5.45) and (5.46), for $d=1,2,3$, by applying Theorem 5 we
shall prove below 
 that 
the 
supports of the  Borel probability measure $\nu_{\epsilon}$ on ${\cal S}'({\mathbb R}^d \to {\mathbb R})$, for each $\epsilon >0$,  and of a weak limit of 
a subsequence of $\{\nu_{\epsilon} \}_{\epsilon >0}$, denoted by $\nu$, are 
 all 
in the Hilbert space 
${\cal H}_{-2} \subset {\cal S}'({\mathbb R}^d \to {\mathbb R})$ defined by (5.11) 
so that from now on 
\begin{equation}
{\mbox{$\nu_{\epsilon}$ and $\nu$ can be understood as probability measures on ${\cal H}_{-2}$}}.
\end{equation}

In fact, by the Sobolev's embedding theorem and by (5.9) with (5.7) and (5.10), we see that 
there exist some constants $K_3, K_4, K_5$ and 
the following inequality with respect to the norms, and the corresponding {\it{continuous}} embeddings hold:
\begin{eqnarray}
\lefteqn{
\sup_{x \in {\mathbb R}^d} \big| (|x|^2 + 1)^{\frac{d+1}2 } f(x) \big| \leq 
K_3 \big\| (|x|^2 +1)^{\frac{d+1}2} f \big\|_{W^{d +1,2} }
}  \\
&&\leq K_4 \big\| (- \Delta +1)^{\frac{d+1}2} (|x|^2 + 1)^{\frac{d + 1}2} f \big\|_{L^2({\mathbb R}^2)} \leq K_5 \|f \|_{{\cal H}_1},
\quad \forall f \in {\cal S}({\mathbb R}^d \to {\mathbb R}), \nonumber
\end{eqnarray}
\begin{equation*}
{\cal H}_1 \, \hookrightarrow \, (|x|^2 + 1)^{- \frac{d+1}2} W^{d +1, 2} ({\mathbb R}^d) \,  \hookrightarrow \,  (|x|^2 + 1)^{- \frac{d+1}2} 
C_b({\mathbb R}^d \to {\mathbb R}),
\end{equation*}
where $W^{d +1,2} = W^{d +1, 2}({\mathbb R}^d)$ is the Sobolev space, the elements of which are real measurable 
functions having square integrable (with respect to the Lebesgue measure on ${\mathbb R}^d$) partial derivatives (in the sense of distribution) of all orders up to $d+1$, and $C_b({\mathbb R}^d \to {\mathbb R})$ is the space of real valued continuous bounded functions on ${\mathbb R}^d$.  

\begin{remark} \qquad 
 By refining the above discussion, passing through  similar arguments, 
it is possible to get  sharper results than (5.49) 
(cf. [A,Gielerak,Russo 2001] for considerations on the path properties of corresponding Euclidean fields, cf., also references therein),  but (5.49)  is sufficient for the subsequent discussions, in particular, for (5.50) below.
\end{remark}

Then, by  Young's inequality, which is valid for both integrals and sums on lattices (cf., e.g., (A.3) in Appendix of [Brydges,Fr{\"o}hlich,Sokal 83]), from (5.44), (5.45), (5.46), (5.47) and (5.49), we see that (cf., 
(5.41) and (5.42), cf. also (5.26))  there exist constants $K'$, $K$  
such that
\begin{eqnarray}
\lefteqn{
\left| \int_{{\cal S}'} \big( {}_{{\cal S}'}< \phi, \varphi>_{{\cal S}'} \big)^2 \, \nu_{\epsilon}(d \phi) \right|
} 
\nonumber \\
&& = \Big| \epsilon^d \sum_{x \in {(\epsilon {\mathbb Z})^d}} \varphi(x)\, \Big( 
\epsilon^d \sum_{y \in {(\epsilon {\mathbb Z})^d}} 
S^{(\epsilon)} (x-y) \varphi (y) \Big) \Big|
\nonumber \\
&&\leq \Big| \epsilon^d \sum_{x \in {(\epsilon {\mathbb Z})^d}} |\varphi(x)| \, \Big( 
\epsilon^d \sum_{y \in {(\epsilon {\mathbb Z})^d}} 
S^{(\epsilon)} (x-y) |\varphi (y)| \Big) \Big|
\nonumber \\
&& \leq 
K'
\big( \epsilon^d \sum_{x \in {(\epsilon {\mathbb Z})^d}} |\varphi(x)| \big) \cdot ( K_1 \lambda^2 + K_2 m_0^{-2}) \cdot \| \varphi\|_{{\infty}}^{(\epsilon)}
\nonumber \\
&&\leq K' \cdot K_5^2 \big( \epsilon^d \sum_{x \in (\epsilon {\mathbb Z})^d} \| \varphi \|^2_{{\cal H}_1} 
(|x|^2 + 1)^ {{-\frac{d+1}2} } \big) (K_1 \lambda^2 + K_2 m_0^{-2})
\nonumber \\
&&\leq  K \| \varphi \|^2_{{\cal H}_1}, 
\qquad \qquad 
\forall \epsilon >0, \, \, \forall \varphi \in {\cal H}_1,
\end{eqnarray}
where $\| {}\,\, \|_{{\cal H}_1}$ is the ${\cal H}_1$ norm defined by (5.11),  and 
we have used the simplified notations
${\cal S} = {\cal S}({\mathbb R}^d \to {\mathbb R})$, 
${\cal S}' = {\cal S}'({\mathbb R}^d \to {\mathbb R})$.
 In (5.50) for 
the third inequality, we used  Young's inequality to get 
$$
\exists K' \quad {\mbox{such that}} \quad 
\| S^{(\epsilon)} \ast \varphi \|_{\infty}^{(\epsilon)} \leq  K' \|S^{(\epsilon)} \|_{l^1}^{(\epsilon)} \cdot
\|\varphi \|_{\infty}^{(\epsilon)}$$
and for the last but one 
inequality, we used the following consequence from (5.49)  
$$
| \varphi (x) | \leq \big( \sup_{x' \in {\mathbb R}^d} 
 (|x' |^2 + 1)^{  {{ \frac{d +1}2 }} }  |\varphi (x')| \big) \cdot ( |x|^2 + 1)^{- \frac{d +1}2} \leq K_5 \| \varphi \|_{{\cal H}_1} (|x|^2 + 1)^{- \frac{d +1}2},$$ 
$$\forall x \in {\mathbb R}^d,
\quad
\forall \varphi \in {\cal H}_1.
$$
By (5.50), from the {\it{Gaussian inequality}} (5.46) we have
\begin{eqnarray}
\lefteqn{
\left| \int_{{\cal S}'} \big( {}_{{\cal S}'}< \phi, \varphi>_{{\cal S}'} \big)^{2n} \, \nu_{\epsilon}(d \phi) \right|
} 
\nonumber \\
&&= 
  \Big| (\epsilon^d)^{2n}  \sum_{x_1, \dots, x_{2n} \in {(\epsilon {\mathbb Z})^d}} \varphi(x_1) \cdots 
\varphi(x_{2n}) \, S^{(\epsilon)}_{2n}(x_1, \dots, x_{2n}) \Big| 
\nonumber \\
&&\leq (2n-1)!! K^n \| \varphi \|^{2n}_{{\cal H}_1},  \quad
\qquad \forall \epsilon > 0, \, \, \forall \varphi \in {\cal H}_1, \,\, n \in {\mathbb N}.
\end{eqnarray}
To derive the second inequality of (5.51), we used the Gaussian inequality and recalled the last formula given in Remark 7. \,
Now, by (5.51) (cf. (5.41), (5.42), (5.43)), we have (5.52) below, where 
the first equality holds, 
since 
$\lim_{N \to \infty} 
\sum_{  { {n =1} } }^{N} \frac{(-1)^n}{(2n)!}
( {}_{{\cal S}'}< \phi, \varphi>_{{\cal S}'} )^{2n} = e^{i  {}_{{\cal S}'}< \phi, \varphi>_{{\cal S}'}} -1$,
\, $\nu_{\epsilon} \, a.e. \, \phi$, 
   by (5.51) since  
$
 \sum_{  { {n =1} } }^{\infty} \frac{1}{(2n)!} \int_{{\cal S}'} 
( {}_{{\cal S}'}< \phi, \varphi>_{{\cal S}'} )^{2n} \nu_{\epsilon}(d \phi)  < \infty
$, holds, 
by Fubini's Lemma, also the second equality holds; 
the last inequality in (5.52) follows again by applying (5.51): 
\begin{eqnarray}
\lefteqn{
\Big| 
\int_{{\cal S}'} e^{i  {}_{{\cal S}'}< \phi, \varphi>_{{\cal S}'}} \nu_{\epsilon} (d \phi) \, \,
 { {-1}}  
\Big|
}
\nonumber \\
&&= \Big| 
\int_{{\cal S}'} \sum_{  
 { {n =1} } }^{\infty} \frac{(-1)^n}{(2n)!}
\big( {}_{{\cal S}'}< \phi, \varphi>_{{\cal S}'} \big)^{2n} \nu_{\epsilon}(d \phi) \Big| 
\nonumber \\
&&= \Big| \sum_{  { {n =1} } }^{\infty} \frac{(-1)^n}{(2n)!} \int_{{\cal S}'} 
\big( {}_{{\cal S}'}< \phi, \varphi>_{{\cal S}'} \big)^{2n} \nu_{\epsilon}(d \phi) \Big| 
\nonumber \\
&& \leq \sum_{  {{n=1}}  }^{\infty} \frac{(2n-1)!!}{(2n)!} K^n \| \varphi \|^{2n}_{{\cal H}_1} = e^{\frac12 K \|\varphi \|_{{\cal H}_1}^2}  \, \,
{ {-1} }, \quad \forall \epsilon>0, \, \, \forall \varphi \in {\cal H}_1.
\end{eqnarray}

%\medskip
%\begin{remark} 
%\qquad 
%(5.52) is an analogue to
% (12.5.1) in [Glimm,Jaffe 87] holding for  (general) Euclidean $P(\phi)_2$ measures.
%\end{remark}

\medskip
\medskip

Denote, for any $\epsilon >0$, $\varphi \in {\cal S}$  
\begin{equation}
{{C}}^{(\epsilon)} (\varphi) \equiv \int_{{\cal S}'} e^{i  {}_{{\cal S}'}< \phi, \varphi>_{{\cal S}'}} \nu_{\epsilon} (d \phi).
\end{equation}
Since $\nu_{\epsilon}$ is a Borel probability measure on ${\cal S}'$, 
${{C}}^{(\epsilon)}$ satisfies conditions ii) and iii) for $C$ in Theorem 5, the Bochner-Minlos theorem.    
Then, by (5.52) and Remark 6, given after Theorem 5,  the continuity of 
${{C}}^{(\epsilon)}(\cdot)$ on ${\cal H}_1$ follows immediatly. Namely, 
 by Remark 6 (precisely, the last formula in the proof of Remark 6),  and (5.52) the following holds:
$$|C^{(\epsilon)}(\varphi) - C^{(\epsilon)}(\psi)|^2 \leq 2 \cdot | C^{(\epsilon)}(\varphi - \psi) -1|
\leq 
2(
e^{\frac12 K \|\varphi - \psi \|_{{\cal H}_1}^2} -1), 
\quad \forall \varphi, \psi \in {\cal H}_1.$$
We thus see  that 
${{C}}^{(\epsilon)}(\cdot)$ is  continuous with respect to the ${\cal H}_1$ norm $\| {}\,\|_{{\cal H}_1}$
 (by the absolute estimate given by (5.51) and (5.52), it is also possible to show the continuity of $C^{(\epsilon)}$ (see Remark 10, given after this Example)).
Hence, 
 from 
i) of Theorem 5 with its last statement (cf. the arguments between (5.14) and (5.18)), we conclude that the support of $\nu_{\epsilon}$ is in ${\cal H}_{-2}$ for any $\epsilon >0$.  This guarantees (5.48) for $\nu_{\epsilon}$.

To certify (5.48) for a weak limit of $\{ \nu_{\epsilon} \}_{\epsilon > 0}$,  as $\epsilon \downarrow 0$, 
 obtained from 
 (5.50), (5.51) and (5.46), 
 we recall that 
the distribution $S^{(\epsilon)}_k$, $k \in {\mathbb N}$, $\epsilon >0$, satisfies 
\begin{equation}
S^{(\epsilon)}_{2n-1} = 0 \quad {\mbox{and}} \quad 
\big\| S^{(\epsilon)}_{2n} \big\|_{({\cal H}_{-1})^{\otimes 2n}} \leq (2n-1) !! \, K^n, \quad \forall \epsilon >0,
\, \, \forall n \in {\mathbb N},
\end{equation}
Namely, for each $n \in {\mathbb N}$, the family of distributions $\{ {\cal S}^{(\epsilon)}_{2n} \}_{\epsilon > 0}$ forms a bounded set in the Hilbert space $({\cal H}_{-1})^{\otimes 2n}$, the space of $2n$-th tensor powers of ${\cal H}_{-1}$ defined by (5.11). Thus, for each $n \in {\mathbb N}$, we can take a sequence $\{\epsilon_{n,i} \}_{i \in {\mathbb N}}$, with 
$\epsilon_{n,i} >0$ and $\lim_{i \to \infty} \epsilon_{n,i} =0$ such that
$\{ {\cal S}^{(\epsilon_{n,i})}_{2n} \}_{i \in {\mathbb N}}$ converges weakly (as $i \to \infty$) to some $S_{2n} \in ({\cal H}_{-1})^{\otimes 2n}$, that satisfies the same bound as (5.51).  By taking subsequences and using a diagonal argument, we then see that there exists a sequence $\{ \epsilon_i \}_{i \in {\mathbb N}}$ with $\epsilon_i > 0$ and 
$\lim_{i \to \infty} \epsilon_i = 0$ such that $\{ S^{\epsilon_i}_{2n} \}_{i \in {\mathbb N}}$ converges weakly (as $i \to \infty$) to $S_{2n} \in ({\cal H}_{-1})^{\otimes 2n}$  for any $n \in {\mathbb N}$.  By these, we can define the functional 
\begin{equation}
C(\varphi) \equiv \sum_{n =0}^{\infty} \frac{(-1)^n}{(2n)!}  \, \big< S_{2n}, \, {\varphi}^{\otimes 2n} \big>,
\end{equation}
that 
 {{is an absolutely convergent series and}  }
satisfies the same bound as (5.52):
\begin{equation}
|C(\varphi) \, \,  { {-1}}  | \leq e^{\frac12 K \|\varphi \|_{{\cal H}_1}^2}
\, \, 
 {{-1}}, \qquad \forall \varphi \in {\cal H}_1,
\end{equation}
where  $\displaystyle{ < S_{2n}, \, {\varphi}^{\otimes 2n} >= {}_{({\cal H}_{-1})^{\otimes 2n}} \big< S_{2n}, \, {\varphi}^{\otimes 2n} \big>_{({\cal H}_1)^{\otimes 2n}}
}$ is the dualization between $({\cal H}_{-1})^{\otimes 2n}$ and $({\cal H}_{1})^{\otimes 2n}$.  For each $\epsilon >0$, since $C^{(\epsilon)}$ defined by (5.53) satisfies the conditions i), ii), iii) of Theorem 5, by the construction $C(\varphi)$ defined by (5.55) it also satisfies the same conditions. In particular, by (5.56),   $C(\varphi)$ is continuous with respect to the 
${\cal H}_1$ norm 
({\it{precisely, see the explanation given after (5.53)}}, also cf. Remark 10 given after this example) .  Hence, from Theorem 5, we deduce the 
existence of  a Borel probability measure $\nu$ on ${\cal H}_{-2}$ corresponding to $C(\varphi)$ defined by (5.55).
This guarantees (5.48) for $\nu$. 
From now on we fix $\nu$ to be the probability measure as follows:
\begin{equation}
{\mbox{$\nu$ is the  probability measure on ${\cal H}_{-2}$ corresponding to $C(\varphi)$ defined by (5.55)}}.
\end{equation}

We thus arrive 
at a situation analogous to the one in 
%at an analogous situation as in 
Example 1. On the space $L^2({\cal H}_{-2}, \nu)$ with the probability measure $\nu$ defined by (5.57), we can construct an ${\cal H}_{-3}$-valued Hunt process that is a stochastic quantization of the Euclidean $\Phi^4_d$, $d = 2,3$,  field 
%(constructed through the lattice argument) 
 with respect to a non-local Dirichlet form.

Precisely, we repeat the analogous discussions between (5.27) and (5.29) 
(cf. also the susequent discussions for the case where $\alpha \in (0,1)$ given after (5.29)).
As was done in Example 1, 
we interpret ${\nu}$ defined by (5.57) as the Borel probability measure  on 
${\cal H}_{-3}$ (that is  wider space than the original domain ${\cal H}_{-2}$).
 By (5.22), (5.23) and (5.24), 
by taking $m = -3$, 
$\tau_{-3}$ defines an isometric isomorphism such that 
\begin{equation*}
\tau_{-3} \, : \, {\cal H}_{-3} \ni f \longmapsto (\lambda_1^{-3}a_1, \lambda_2^{-3}a_2, \dots) \in l^2_{(\lambda_i^{6})},
%\quad {\mbox{with}} \quad 
%a_i \equiv (f, \, \lambda_i^{-3} {\varphi}_i)_{-3}, \, \, \, i \in {\mathbb N}.
\end{equation*}
\begin{equation}
{\mbox{with}} \qquad 
a_i \equiv (f, \, \lambda_i^{-3} {\varphi}_i)_{-3}, \quad i \in {\mathbb N}.
\end{equation}
Define a probability measure $\mu$ on 
$ l^2_{(\lambda_i^{6})}$ such that 
\begin{equation}
\mu(B) \equiv \nu \circ \tau^{-1}_{-3}(B) \quad {\mbox{for}} \quad B \in {\cal B}(l^2_{(\lambda_i^{6})}).
\end{equation}
We set $S = l^2_{(\lambda_i^{6})}$ in Theorems 1, 2 and 4, then it follows that for the weight $\beta_i$ we have $\beta_i = \lambda_i^6$. We can take ${\gamma_i}^{-\frac12} = \lambda_i$ in Theorem 2-i) 
with $p =2$,  then, from (5.18) we get 
\begin{eqnarray}
\lefteqn{
\sum_{i=1}^{\infty} 
(\beta_i \gamma_i)^{\frac{\alpha +1}{2}} \cdot
\mu \Big( \beta_i^{\frac12} |X_i| > M \cdot \gamma_i^{-\frac12}\Big) 
} \nonumber \\
&\leq& \sum_{i =1}^{\infty} (\beta_i \gamma_i)^{\frac{\alpha +1}{2} } 
= \sum _{i =1}^{\infty} (\lambda_i)^{2(\alpha +1)} < \infty, \qquad \forall \alpha \in (0,1].
%\nonumber 
\end{eqnarray}
(5.60) shows that the condition (4.3) holds.

%{\textcolor{blue}{
Also, as has been mentioned above, 
since $\nu({\cal H}_{-n}) = 1$,  for any   $ n \geq 2$, we have 
$$1 = \nu({\cal H}_{-2}) = \mu(l^2_{(\lambda_i^4)}) = 
\mu \big( \bigcup_{M \in {\mathbb N}} \{ |X_i| \leq M \beta_i^{- \frac12} \gamma_i^{- \frac12}, \, \forall i \in {\mathbb N} \} \big),$$
$$
 {\mbox{for \quad $\beta_i = \lambda_i^6$, \, $\gamma_i^{- \frac12} = \lambda_i$}}.
$$
This shows that the condition (4.4) is satisfied.
%}}

Thus, by Theorem 2-i) and Theorem 4,  
for each  $0 < \alpha \leq 1$,
there exists an  $l^2_{(\lambda_i^{6})}$-valued  Hunt process 
\begin{equation}
{\mathbb M} \equiv \big(\Omega, {\cal F}, (X_t)_{t \geq 0}, (P_{\mathbf x})_{\mathbf x \in 
S_{\triangle}} \big), 
\end{equation}
associated  to the non-local Dirichlet form
$({\cal E}_{(\alpha)}, {\cal D}({\cal E}_{(\alpha)}))$.
We can then define an
${\cal H}_{-3}$-valued 
 process 
$(Y_t)_{t \geq 0}$ 
such that 
$(Y_t)_{t \geq 0} \equiv \big({\tau}^{-1}_{-2}(X_t) \big)_{t \geq 0}.$
Equivalently, by (5.58) for $X_t = (X_1(t), X_2(t), \dots) \in l^2_{(\lambda_i^6)}$, $P_{\mathbf x}-a.e. \, x \in S_{\triangle}$, by setting $A_i(t)$ such that 
$X_i(t) = \lambda_i^{-3} A_i(t)$  (cf. (5.21) and (5.22)),  then 
$Y_t$ \, is given by 
\begin{equation}
Y_t = \sum_{i \in {\mathbb N}} A_i(t) (\lambda_i^{-3} \varphi_i) = \sum_{i \in {\mathbb N}} X_i(t) \varphi_i \in {\cal H}_{-3}, \qquad \forall t \geq 0, \, \, P_{\mathbf x}-a.e.,
\, x \in S_{\triangle}.
\end{equation}
By (5.4) and (5.58), $Y_t$  is an ${\cal H}_{-3}$-valued Hunt process that is a {\it{stochastic quantization}} with respect to the 
non-local Dirichlet form 
$({\tilde{\cal E}}_{(\alpha)}, {\cal D}({\tilde{\cal E}}_{(\alpha)}))$ 
on $L^2({\cal H}_{-3}, \nu)$, that is defined through
$({\cal E}_{(\alpha)}, {\cal D}({\cal E}_{(\alpha)}))$, by making use of $\tau_{-3}$ via (5.59). 
We state the above results as a theorem:
\begin{theorem}
Let $\nu$ be the Euclidean $\Phi^4_d$, $d=2, 3$,  field  measure defined by (5.57). Interpret $\nu$ as a Borel probability measure on ${\cal H}_{-3}$, and let $\mu$ be the Borel probability measure on $l^2_{(\lambda_i^6)}$ that is an image of $\nu$ defined by (5.59). Then, 
for any $0 < \alpha \leq 1$, 
 on $L^2(l^2_{(\lambda_i^6)}; \mu)$,  
 a family of 
non-local quasi-regular Dirichlet form $({\cal E}_{(\alpha)}, {\cal D}({\cal E}_{(\alpha)}))$ can be defined through Theorem 1 and Theorem 2-i), and by Theorem 4,   there exists an $S \equiv l^2_{(\lambda_i^6)}$-valued Hunt process ${\mathbb M} \equiv \big(\Omega, {\cal F}, (X_t)_{t \geq 0}, (P_{\mathbf x})_{\mathbf x \in 
S_{\triangle}} \big)$ properly associated to the Dirichlet form 
 $({\cal E}_{(\alpha)}, {\cal D}({\cal E}_{(\alpha)}))$
(cf. (5.61)). 
Moreover, the stochastic process $(Y_t)_{t \geq 0}$ defined by (5.62) through ${\mathbb M}$ 
is an ${\cal H}_{-3}$-valued Hunt process that is a {\it{stochastic quantization}} with respect to the 
non-local Dirichlet form 
$({\tilde{\cal E}}_{(\alpha)}, {\cal D}({\tilde{\cal E}}_{(\alpha)}))$ 
on $L^2({\cal H}_{-3}, \nu)$, that is defined through
$({\cal E}_{(\alpha)}, {\cal D}({\cal E}_{(\alpha)}))$ by making use of $\tau_{-3}$ defined by (5.58). \, \, 
\bsquare
\end{theorem}

%\medskip
\begin{remark} 
\qquad 
(5.52) is an analogue to
 (12.5.1) in [Glimm,Jaffe 87] holding for  (general) Euclidean $P(\phi)_2$ measures.
\end{remark}

\begin{remark}
\qquad 
By (5.51), since 
$
 \sum_{   {n =0} }^{\infty} \frac{(-1)^n}{(2n)!} \int_{{\cal S}'} 
( {}_{{\cal S}'}< \phi, \varphi>_{{\cal S}'} )^{2n} \nu_{\epsilon}(d \phi) 
$, 
and 
$
C(\varphi) \equiv \sum_{n =0}^{\infty} \frac{(-1)^n}{(2n)!}  \, \big< S_{2n}, \, {\varphi}^{\otimes 2n} \big>,
$
converge 
 absolutely, for $\varphi, \, \psi \in {\cal S}$ we are admitted to perform the following evaluation
(cf. (5.53)), and see, for e.g.,  the continuity of $C^{(\epsilon)}(\varphi)$ with respect to $\varphi \in {\cal H}_1$:
\begin{eqnarray*}
\lefteqn{
|C^{(\epsilon)}(\varphi) - C^{(\epsilon)}(\psi)| = 
 \big| \int_{{\cal S}'} e^{i  {}_{{\cal S}'}< \phi, \varphi>_{{\cal S}'}} \nu_{\epsilon} (d \phi) 
- \int_{{\cal S}'} e^{i  {}_{{\cal S}'}< \phi, \psi>_{{\cal S}'}} \nu_{\epsilon} (d \phi)
\big|
}
\nonumber \\
&&=
\big| \int_{{\cal S}'} e^{i  {}_{{\cal S}'}< \phi, \varphi>_{{\cal S}'}}  
( 1 -  e^{i  {}_{{\cal S}'}< \phi, \psi - \varphi>_{{\cal S}'}}) 
\nu_{\epsilon}(d \phi) 
\big|
\nonumber \\
&& \leq | 
\int_{{\cal S}'} \sum_{  
 { {n =1} } }^{\infty} \frac{(-1)^n}{(2n)!}
\big( {}_{{\cal S}'}< \phi, \psi -\varphi>_{{\cal S}'} \big)^{2n} \nu_{\epsilon}(d \phi) | 
\nonumber \\
%&&= \Big| \sum_{  {\textcolor{red} {n =1} } }^{\infty} \frac{(-1)^n}{(2n)!} \int_{{\cal S}'} 
%\big( {}_{{\cal S}'}< \phi, \varphi>_{{\cal S}'} \big)^{2n} \nu_{\epsilon}(d \phi) \Big| 
%\nonumber \\
&& \leq \sum_{  { {n=1}}  }^{\infty} \frac{(2n-1)!!}{(2n)!} K^n \| \varphi - \psi \|^{2n}_{{\cal H}_1} = e^{\frac12 K \|\varphi - \psi \|_{{\cal H}_1}^2}  \, \,
{ {-1} }, \quad \forall \epsilon>0.
\end{eqnarray*}
In the above, for the inequality of the third line we 
make use of the H{\"o}lder inequality with 
$|e^{i  {}_{{\cal S}'}< \phi, \varphi>_{{\cal S}'}}|  =1$, and 
the  evaluation 
of the fourth line 
holds by (5.51), which follows from the {\it{Gaussian inequality}} (see (5.46)). 

We note that for the Poisson random fields 
the {\it{Gaussian inequality}} does not hold (see [A,Gottschalk,Wu-j.L. 97] for the corresponding detailed discussions).
%also for 
%for the {\it{non-Gaussian}} $\alpha$-stable ($0 < \alpha <2$) type Euclidean random fields, denoted by $\phi$,  even the second 
%moment of $<\phi, \varphi>$, $\varphi \in {\cal S}$ does not exists.
\end{remark}

%{\textcolor{blue}{
A martingale representation of Markov processes properly associated to Dirichlet forms is known as {\it Fukushima decomposition} (for the case where the state spaces $S$ are locally compact  metric spaces cf., e.g., 
  [Fukushima 80], [F,Oshima,Takeda 2011], 
and Theorem 4.3 of [F,Uemura 2012], and for the case where $S$ are general Hausdorff topological spaces cf., e.g.,  Chapter VI of [MR 92] and [A,R{\"u}diger 2003]).
Corresponding to {\it  Fukushima decomposition}, for Theorem 4 and Examples 1, 2 we have the following {\it Remarks} 
{{11 and 12}}.
%}}

\begin{remark} {
%{\textcolor{blue}{
\bf{(Fukushima decomposition)}
%}}
} \\
The subspace  of $S$, in which the $\mu$-tight special standard process and Hunt process  properly associated with the quasi-regular Dirichlet form and strictly quasi-regular Dirichlet form  
on $L^2(S;\mu)$ 
takes values, respectively, can be embedded in some locally compact separable metric space (cf. Definition IV-3.1 and Theorem VI-1.2 of [M,R 93], and also cf. the proof of Theorem 2 of the present paper).  In what follows, the interpretation of Theorem 5.2.2 of [F,Oshima,Takeda 2011] is straightforward, and then Theorem VI-2.5 of [M,R 93] holds.

We use the notions and notations adopted in chapter VI of [M,R 93] (cf. also the notations in  chapter 5 of [F,Oshima,Takeda 2011] and  [A,Ma,R 2015], [A,Ma,R1 93], [A,Ma,R2 93],  [A,R{\"u}diger 2013] and references therein). 

Let 
${\mathbb M} \equiv \big(\Omega, {\cal F}, (X_t)_{t \geq 0}, (P_{\mathbf x})_{\mathbf x \in 
S_{\triangle}} \big)$ 
be the Hunt process defined through Theorem 4. By a direct application of Theorem VI-2.5 of [M,R 93] we see that for $u \in {\cal D}({\cal E}_{(\alpha)})$, there exists a unique martingale additive functional of finite energy (MAF) $M^{[u]}$ and a continuous additive functional of zero energy (CAF's zero energy) $N^{[u]}$ such that 
\begin{equation}
A^{[u]} = M^{[u]} + N^{[u]},
\end{equation}
where 
$$A^{[u]} \equiv (A^{[u]}_t)_{t \geq 0}, \qquad A_t^{[u]} = {\tilde u}(X_t) - {\tilde u}(X_0),
$$
with ${\tilde u}$ an ${\cal E}_{(\alpha)}$-quasi continuous $\mu$-version of $u \in {\cal D}({\cal E}_{(\alpha)})$.  The decomposition formula (5.63) holds for Examples 1 and 2.

In order to consider the martingale problems (cf., e.g., [F,Uemura 2012]) corresponding  
 to 
the decomposition given by 
(5.63), some additional assumptions for the probability measure $\mu$, e.g., a uniform regularity of its density function (cf. (2.11)), are necessary.  Considerations in this direction 
can be found in  [A,Kaga,Kawa,Yaha,Y 2021].  
\end{remark}

\begin{remark} {
%{\textcolor{blue}{
\bf{(Subordination correspondences)}
%}}
} \\
A theory of 
transforming Dirichlet forms and associated symmetric Markov processes by means of the subordinations 
 has been developed both from a  functional analytic and a probabilistic 
point of view  (for the case where the state spaces are locally compact spaces, cf. [Jacob,Schilling 2005], [Schilling 98] 
and references therein, and for the case where the state space are general Hausdorff spaces, cf. [A,R{\"u}diger 2003]).

On $L^2(S;\mu)$ with $\mu$ a Borel probability measure on a general Hausdorff topological space $S$, consider in general a quasi-regular Dirichlet form $({\cal E}, {\cal D}({\cal E}))$.  Let $-L$ be the self-adjoint operator corresponding to $({\cal E}, {\cal D}({\cal E}))$\\ (thus, $Dom(\sqrt{-L}) = {\cal D}({\cal E})$), and $P_{\lambda}$ be the projection valued measure associated to the operator $-L$.  Let $f$ be a Bernstein function on ${\mathbb R}_+$
(cf., e.g., [Schilling 98] 
for its definition and properties), and define a self-adjoint  operator $L^f$ by 
$$L^f \equiv -f(-L) = - \int_0^{\infty} f(\lambda) \, dP_{\lambda}.$$
In [A,R{\"u}diger 2003]  the following mathematical structure of the correspondences between subordinate symmetric processes and the subordinate sub-Markov semi-groups, and hence the associated Dirichlet forms, is introduced (see Theorems 2.7, 2.9, equations (8)(8), (9), Theorem 2.16 and Theorem 3.1 of [A,R{\"u}diger]): \, 
For a quasi-regular Dirichlet form $({\cal E}, {\cal D}({\cal E}))$ as above
\begin{equation}
{\cal E}^f (\phi, \psi) \equiv \left( \sqrt{f(-L)} \phi, \sqrt{f(-L)} \psi \right)_{L^2(S;\mu)},
\, \phi, \psi \in {\cal D}({\cal E}^f) \equiv Dom(\sqrt{f(-L)}),
\end{equation}
defines a non-local quasi-regular Dirichlet form. 
Let $(X^f_t)_{t \geq 0}$ be the $\mu$-tight special standard process properly associated  to the quasi-regular Dirichlet form $({\cal E}^f, {\cal D}({\cal E}^f))$, the corresponding semigroup of which is denoted by 
$(T^f_t)_{t \geq 0}$, 
then $(X^f_t)_{t \geq 0}$ has the same finite dimensional distributions as $(X_{y(t)})_{t \geq 0}$, where $(X_t)_{t \geq 0}$ is the $\mu$-tight special standard process properly associated to the quasi-regular Dirichlet form $({\cal E}, {\cal D}({\cal E}))$ and $(y(t))_{t \geq 0}$ is the increasing L{\'e}vy process defined through the Bernstein function $f$. 

By Theorem VI-2.5 of [M,R 93] (cf. Remark 11), the $\mu$-tight special standard process
$(X^f_t)_{t \geq 0}$ admits the Fukushima decomposition  (5.63).  [A,R{\"u}diger 2003] discusses the Fukushima decomposition and the corresponding martingale problem for the 
process $(X^f_t)_{t \geq 0}$ (see Examples, 1, 2, 3, 4 and Theorem 4.29).

The investigations of the correspondences between the framework of the subordination and the general  framework 
 presented in the present work 
are particularly interesting  and deserve further consideration.
\end{remark}

\section{Future developments.}
Let us add some short comments on future developments corresponding to the present work, 
 possibly encouraging future works by 
  researchers working in  related areas.

The present paper 
is intended to provide an 
 explicit formulation of {\it{non-local}} Dirichlet forms defined on infinite dimensional topological vector spaces.  
Our definitions  
(2.8), (2.9) and (2.10) 
of the Markov symmetric form ${\cal E}_{(\alpha)}$, $0 < \alpha <2$,  
 can be extended in  
several directions. In analogy with  the finite dimensional cases 
 (cf., e.g., [Fukushima 80],  
[A,Song 93], 
[F,Oshima,Takeda 2011], [F,Uemura 2012], [Hoh,Jacob 96], [Masamune,Uemura,Wang 2012],
\,{}\\ \![Schlling,Wang 2014], 
  [Shiozawa,Uemura 2014] and references therein), the kernel 
$\frac{1}{|y_i -y'_i|^{\alpha +1}}$ can  be replaced by more general symmetric ones 
 (and also some 
 non-symmetric ones).

Also the connections between the Hunt processes constructed by the {\it{non-local}} Dirichlet forms on infinite dimensional topological vector spaces and  solutions of  SPDEs   
seem to provide an area of possible extensions of our work 
(cf. e.g., [A,Di Persio,Mastrogiacomo,Smii 2016] and references therein).

The study of contractivity properties of semi-groups with generators associates to Dirichlet forms is a 
particularly 
important subject. E.g., it is well known that  the classical (local) Dirichlet forms, associated to certain Gaussian random fields,
 satisfy  
the logarithmic Sobolev inequality,  which guarantees a spectral gap of the associated self adjoint operator 
and the semi-group corresponding to the operator satisfies the hypercontractivity (for corresponding precise results, cf. for  [Gross 93], [Simon 74],  section X.9 of [Reed,Simon 75], [Aida 2012] and references therein).
 Analogous considerations for the {\it{non-local}} Dirichlet forms associated with  random fields or processes on infinite dimensional 
topological vector spaces is another possible direction for future investigations.

%{\textcolor{blue}{A comment for the co-workers: \,  As an original idea, we may assume that 
%$p= 2$, \, ${\beta}_n = \frac{1}{({n!})^2}$ and $b_n = n!$, and  
%the variable $x_i \in {\mathbb R}$ of the measure $\mu$  corresponds to 
%the random variable $X_i = {\varphi}_i(\phi)$ for ${\varphi}_i \in {\cal S}({\mathbb R}^d)$, with $\phi$ an Euclidean quantum field.}}

                                                                %\begin{figure}
                                                                %\includegraphics{example.eps}
                                                                %\end{figure}
                                                                %\begin{figure*}
                                                                %\includegraphics[width=0.75\textwidth]{example.eps}
                                                                %\end{figure*}
                                                  %\begin{table}
                                                  %\begin{tabular}{lll}
                                                  %\hline\noalign{\smallskip}
                                                  %\end{tabular}
                                                  %\end{table}

\bigskip
\bigskip

{\small{\bf{Dedication}}
 Dedicated to Professor Masatoshi Fukushima for his 88th birth day. Dedicated to  Professor Rahael H{\o}egh-Krohn for his 84th birth anniversary .

\bigskip
 
   {\small{\bf{Acknowledgements}}  
The authors would like to gratefully acknowledge the great hospitality of various institutions.  In particular 
for the  
first named and 
second named authors, 
Giulio Casati and the 
Lake Como school of advanced studies, Complexity and Emergence: Ideas, Methods, with a special attention to Economics and Finance, Italy; for the third and fourth named authors, IAM and HCM at the University of Bonn, Germany; for the fourth named author, SFB 1283 and Bielefeld University, Germany. 
Also, the fourth named author expresses his strong acknowledgements to Professor Michael R{\"o}ckner 
for fruitful discussions.
}


%{\bf{References should be extended and corrected.}}

\end{document}